\documentstyle[11pt,mssymb]{article}

\newcommand{\arrow}{\longrightarrow}
\newcommand{\Z}{{\Bbb Z}}
\newcommand{\C}{{\Bbb C}}
\newcommand{\R}{{\Bbb R}}

\newcommand{\1}{\sqrt{-1}\:}
\newcommand{\g}{{\goth g}}
\newcommand{\inangles}[1]{{\langle #1\rangle}}
\newcommand{\galochka}{\check{\;}}
\newcommand{\restrict}[1]{{|_{{\phantom{|}\!\!}_{#1}}}}
\renewcommand{\c}[1]{{\cal #1}}

\renewcommand{\phi}{\varphi}
\renewcommand{\epsilon}{\varepsilon}
\renewcommand{\geq}{\geqslant}
\renewcommand{\leq}{\leqslant}

\newcommand{\inbfpare}[1]{{%
  \mbox{\tt (}\hspace{-5pt}\mbox{\tt (} #1 %
  \mbox{\tt )}\hspace{-5pt}\mbox{\tt )}%
}}
\newcommand{\comment}[1]{{}}

\def\blacksquare{\hbox{\vrule width 4pt height 4pt depth 0pt}}

\makeatletter

  \renewcommand{\section}{%
     \secdef\SectionNumbered\SectionStarred}

  \def\SectionNumbered[#1]#2{%
      \refstepcounter{section}
      \addcontentsline{toc}{section}
          {\protect\numberline{\thesection}. #1}
      {\vspace{3\baselineskip}\Large\bf
          \thesection.\hspace{1em}#2\par}
      \sectionmark{#1}\vspace{1.5\baselineskip}}

  \newcommand{\SectionStarred}[1]{%
      {\vspace{3\baselineskip}\Large\bf #1\par}
      \sectionmark{#1}\vspace{1.5\baselineskip}}

   \renewcommand{\thesection}{\arabic{section}}
   \@addtoreset{footnote}{section}
   \@addtoreset{equation}{section}

\newcounter{lemma}[section]
\renewcommand{\thelemma}{{Lemma \thesection.\arabic{lemma}}}
\newcommand{\lemma}{%
     \refstepcounter{lemma}
     {\bf \thelemma:\ }}

\newcounter{claim}[section]
\renewcommand{\theclaim}{{Claim \thesection.\arabic{claim}}}
\newcommand{\claim}{%
     \refstepcounter{claim}
     {\bf \theclaim:\ }}

\newcounter{sublemma}[section]
\renewcommand{\thesublemma}{{Sublemma \thesection.\arabic{sublemma}}}
\newcommand{\sublemma}{%
     \refstepcounter{sublemma}
     {\bf \thesublemma:\ }}

\newcounter{corollary}[section]
\renewcommand{\thecorollary}{{Corollary \thesection.\arabic{corollary}}}
\newcommand{\corollary}{%
     \refstepcounter{corollary}
     {\bf \thecorollary:\ }}

\newcounter{theorem}[section]
\renewcommand{\thetheorem}{{Theorem \thesection.\arabic{theorem}}}
\newcommand{\theorem}{%
     \refstepcounter{theorem}
     {\bf \thetheorem:\ }}

\newcounter{proposition}[section]
\renewcommand{\theproposition}
       {{Proposition \thesection.\arabic{proposition}}}
\newcommand{\proposition}{%
     \refstepcounter{proposition}
     {\bf \theproposition:\ }}

\newcounter{definition}[section]
\renewcommand{\thedefinition}
       {{Definition \thesection.\arabic{definition}}}
\newcommand{\definition}{%
     \refstepcounter{definition}
     {\bf \thedefinition:\ }}

\newcounter{example}[section]
\newcounter{problem}[section]
\newcounter{question}[section]
\newcommand{\eqref}[1]{(\ref{#1})}

 \@ifundefined{boxed}{%
     }{}

\newcommand{\ps@verbit}{%
  \renewcommand{\@oddhead}{%
          \scriptsize
          {Period maps and cohomology}
          \hfil\tiny {final version, December 1994, revised May 95}}
  \renewcommand{\@evenhead}{\@oddhead}
  \renewcommand{\@oddfoot}{\hfil\thepage\hfil}
  \renewcommand{\@evenfoot}{\@oddfoot}}

\makeatother

\begin{document}

\begin{center}
\Large\bf Cohomology of compact hyperk\"ahler manifolds.
\end{center}

\centerline{Mikhail Verbitsky,}
\centerline{verbit@math.harvard.edu}

\hfill

\hfill

\hfill

\hfill

\comment{ Let M be a compact simply connected hyperk\"ahler
(or holomorphically symplectic) manifold, \dim H^2(M)=n.
Assume that M is not a product of hyperkaehler manifolds.
We prove that the Lie algebra so(n-3,3) acts by automorphisms
on the cohomology ring H^*(M). Under this action,
the space H^2(M) is isomorphic to the fundamental
representation of so(n-3,3). Let A^r be the subring of
H^*(M) generated by H^2(M). We construct an action of the
Lie algebra so(n-2,4) on the space A, which preserves
A^r. The space A^r is an irreducible representation
of so(n-2,4). This makes it possible to compute
the ring A^r explicitely.}

\hspace{0.2\linewidth}\begin{minipage}[t]{0.7\linewidth}
{\bf Abstract}. Let $M$ be a compact simply connected hyperk\"ahler
(or holomorphically symplectic) manifold, $\dim H^2(M)=n$.
Assume that M is not a product of hyperkaehler manifolds.
We prove that the Lie group $\goth{so}(n-3,3)$ acts by automorphisms
on the cohomology ring $H^*(M)$. Under this action,
the space $H^2(M)$ is isomorphic to the fundamental
representation of $\goth{so}(n-3,3)$. Let $A^r$ be the subring of
$H^*(M)$ generated by $H^2(M)$. We construct an action of the
Lie algebra $\goth{so}(n-2,4)$ on the space $A$, which preserves
$A^r$. The space $A^r$ is an irreducible representation
of $\goth{so}(n-2,4)$. This makes it possible to compute
the ring $A^r$ explicitely.
\end{minipage}

\hfill


\section{Introduction.}\label{_introduction_to_so(b_2,...)_Section_}


Here we give a brief introduction to results of this
paper. The subsequent sections are independent from the
introduction.

The main object of this paper is theory of compact hyperk\"ahler manifolds.
A hyperk\"ahler manifold (see \ref{_hyperkaehler_manifold_Definition_}
for more precise wording) is a Riemannian manifold $M$ equipped
with three complex structures $I$,
$J$ and $K$, such that $I\circ J=-J\circ I=K$ and $M$
is K\"ahler with respect to the complex structures $I$,
$J$ and $K$.

Let $M$ be a complex manifold which admits a hyperk\"ahler
structure. A simple linear-algebraic argument implies that
$M$ is equipped with a holomorphic symplectic form. Calabi-Yau theorem
shows that, conversely, every compact holomorphically symplectic
K\"ahler manifold admits a hyperk\"ahler structure, which is uniquely
defined by these data.

Further on, we do not discriminate between compact holomorphic
symplectic manifolds of K\"ahler type and compact hyperk\"ahler
manifolds.

Let $(M,I)$ be a compact K\"ahler manifold which admits a holomorphic
symplectic form $\Omega$. For simplicity of statements, we assume
in this introduction that $\Omega$ is unique up to a constant.
Denote by $M$ the $C^\infty$-manifold underlying $(M,I)$.
Let $(M,I')$ be another compact holomorphically symplectic
manifold which lies in the same deformation class as $(M,I)$.
Fixing a diffeomorphism of underlying $C^\infty$-manifolds,
we may identify the smooth manifold which underlies $(M,I)$
with that underlying $(M,I')$.

Let $X$ be an arbitrary compact K\"ahler manifold,
$\dim_\C X=n$. We associate with a K\"ahler structure on $X$
so-called Riemann-Hodge pairing

\[
   (\cdot,\cdot): \;\; H^2(X,\R)\times H^2(X,\R)\arrow \R
\]
which is a map associating a number

\begin{equation}\label{_Hodge-Riemann-correct-Equation_}
   \int_X \omega^{n-2}\wedge \eta_1\wedge\eta_2  -
   \frac{n-2}{(n-1)^2}\cdot
   \frac{\int_X \omega^{n-1}\eta_1 \cdot \int_X\omega^{n-1}\eta_2}
   {\int_X \omega^n}
\end{equation}
to classes $\eta_1,\eta_2 \in H^2(X,\R)$, where
$\omega\in H^2(X,\R)$ is the K\"ahler class
(see also \ref{_Hodge_Riema_general_Claim_}).

Consider the positively defined scalar product
$(\cdot,\cdot)_{Metr}$ induced by Riemannian metric
on  the space of harmonic 2-forms, identified by Hodge
with $H^2(X,\R)$.
The pairing \eqref{_Hodge-Riemann-correct-Equation_}
is defined in such a way that on primitive $(1,1)$-forms
with coefficients in $\R$, $(\cdot,\cdot)$ is equal
to $-(\cdot,\cdot)_{Metr}$. Similarly, on the space
\[ \bigg(\R\omega \oplus H^{2,0}(X)\oplus H^{0,2}(X)\bigg) \cap H^2(X,\R), \]
the form $(\cdot,\cdot)$ is equal to $(\cdot,\cdot)_{Metr}$.

\hfill

Let $(\cdot,\cdot)$,
$(\cdot,\cdot)': \; H^2(M,\R)\times H^2(M,\R)\arrow \R$
be the Hodge-Riemann forms associated with $(M,I)$ and $(M,I')$
respectively. The most surprising result of this paper is
following (see \ref{_Hodge_Riemann_independent_Theorem_}
for a different wording of the same statement):

\hfill

\theorem \label{_Riemann_Hodge_unique_in_intro_Theorem_}  
The forms $(\cdot,\cdot)$ and $(\cdot,\cdot)'$
are {\bf proportional}.

\hfill

Taking a K\"ahler class $\omega$ such that
$Vol(M)= 1$, where $Vol(M):= \int_M \omega^{n}$,
we get rid of the ambiguity in the choice of a constant.
If $(M,I)$ and $(M,I')$ both satisfy $Vol(M)= 1$,
then the Hodge-Riemann forms $(\cdot,\cdot)$ and $(\cdot,\cdot)'$
are {\bf equal}.
We call this form {\bf the normalized Hodge-Riemann pairing},
denoted as $(\cdot,\cdot)_{\c H}$. This form is an
invariant of a deformation class of complex manifolds.
One may regard $(\cdot,\cdot)_{\c H}$ as {\bf topological} invariant.
\footnote{In fact, the form
$(\cdot,\cdot):\;H^2(M,\R)\times H^2(M,\R)\arrow \R$ of
\eqref{_Hodge-Riemann-correct-Equation_} is a
{\it topological} invariant associated
with every cohomology class $\omega\in H^2(M,\R)$.
Our proof immediately implies that
for a dense set of $\omega,\omega'\in H^2(M,\R)$,
\ref{_Riemann_Hodge_unique_in_intro_Theorem_} holds.}

We use the form $(\cdot,\cdot)_{\c H}$ to
compute the cohomology algebra $H^*(M)$. We explicitely
compute the subalgebra of $H^*(M)$ generated by $H^2(M)$.
These computations also give the relations between elements
generated by $H^2(M)$ and the rest of $H^*(M)$.

The main ideas of this computation are due to the observations which
involve Torelli theorem and deformation spaces. As follows
from results of Bogomolov
(\cite{_Bogomolov_}, \cite{_Besse:Einst_Manifo_}, \cite{_Todorov_}),
the deformation space of a holomorphically symplectic manifold
is a smooth complex manifold, which is a quotient of a
Stein space by an arithmetic group. From this description
we use only the calculation of dimension of this moduli space.

The map of Kodaira and Spencer

\[ KS:\; T_{[M]}\c M\arrow H^1(TM)
\]
is a linear homomorphism
from the Zariski tangent space $T_{[M]}\c M$ of a moduli space $\c M$ of
deformations of $M$ to the first cohomology of the sheaf of holomorphic
vector fields on $M$. It is proven by Kodaira and Spencer
\cite{_Kodaira_Spencer_} that this
map is an embedding. Results of Bogomolov et al
imply that $\c M$ is smooth and
$KS$ is an isomorphism. This is often called ``Torelli
theorem''. This statement could be translated
to the language of period maps.
Together with \ref{_Riemann_Hodge_unique_in_intro_Theorem_},
this observation implies a nice description of the period map.

By ``period map'' associated with the holomorphic symplectic manifolds
(see Section \ref{_perio_and_forge_Section_}
for details) we mean the following geometrical object.
We define the (coarse, marked) moduli of holomorphic symplectic
manifolds as the space of different complex  holomorphic symplectic
structures on a given $C^\infty$-manifold $M$ up to
diffeomorphisms which act trivially on
$H^*(M)$. Let $Comp$ be a connected component of this moduli space.
For more accurate definition of $Comp$, we refer the reader to
\ref{_Comp_Hyp_Definition_}. For all points $I\in Comp$, we
denote the corresponding complex manifold by $(M,I)$.
{}From the definition of $Comp$ we obtain a canonical
identification of cohomology spaces $H^*(M,I)$ for all $I\in Comp$.

For simplicity of statements, we deal with the simply connected
holomorphically symplectic manifolds $M$ such that $H^{2,0}(M)=\C$
(such manifold are called {\bf simple}).

The period map (Section \ref{_perio_and_forge_Section_})
in this context is a map $P_c:\; Comp \arrow \Bbb P(H^2(M,\C))$
which relies a line $H^{2,0}(M,I)\subset H^2(M,\C)$ to every $I\in Comp$.
Classical results of Bogomolov et al imply that
$P_c$ is an immersion.
Hodge-Riemann relations together with
\ref{_Riemann_Hodge_unique_in_intro_Theorem_}
let one to describe the image of $P_c$ in following
terms. An immediate consequence of Hodge-Riemann relations
is that for all lines $x\in Im(P_c)\subset \Bbb P(H^2(M,\C))$,
and all vectors $l\in x\subset H^2(M,\C)$, we have
$(l,l)_{\c H} =0$. Therefore, $P_c$ maps $Comp$ to a
quadric $C\subset \Bbb P(H^2(M,\C))$, which is defined
by the quadratic form associated with the Riemann-Hodge pairing.
Dimension of $Comp$ is computed from Torelli theorem.
As one can easily check, it is equal to the dimension
of $C$. Since $P_c$ is an immersion, we obtain the following theorem:

\hfill

\theorem 
The period map $P_c:\; Comp \arrow C$ is etale.

\hfill

This is the main
observation used in computations of the cohomology algebra
$H^*(M)$.

\hfill

Let $I\in Comp$. Let $ad I\in End(H^*(M))$ be a linear endomorphism
which maps a $(p,q)$-form $\eta\in H^{p,q}(M)$ to $(p-q)\1 \eta$.
Let $\goth M\subset End(H^*(M))$ be a Lie algebra generated
by the endomorphisms $ad I$ for all $I\in Comp$.
Hodge-Riemann relations imply that the action
of $\goth M$ on $V= H^2(M)$ preserves the scalar
product $(\cdot,\cdot)_{\c H}$. This defines a
Lie algebra homomorphism

\begin{equation} \label{_from_Mum-Tate_to_SO(V)_homom_Equation_}
   \rho:\; \goth M\arrow \goth{so}(V).
\end{equation}
Using period maps and estimation on dimensions of moduli spaces, we
prove that $\rho$ is an {\bf isomorphism}
(\ref{_g_0_computed_Theorem_}, \ref{_g_0_is_Mumf_Tate_Theorem_}).
This statement is an ingredient in the computation of the
algebra $H^*(M)$.

The algebraic structure on $H^*(M)$ is studied using the
general theory of Lefschetz-Frobenius algebras, introduced
in \cite{_Lunts-Loo_}. Lefschetz-Frobenius
algebras are associative graded commutative algebras whose
properties approximate that of cohomology of compact
manifolds which admit K\"ahler structure. We give an
exact definition of this term in Section \ref{_Lefshe_Frob_Section_}
(\ref{_Lefschetz_Frob_alge_Definition_}).
With no loss of generality, reader may think of
these algebras as of cohomology algebras.

With the Lefschetz-Frobenius algebra $A$ we associate
so-called {\bf structure Lie algebra} $\g(A)\subset End(A)$
which acts on $A$ by linear endomorphism. The action of structure
Lie algebra is often sufficient to reconstruct multiplication
on $A$. This algebra is defined using the algebraic version
of the strong Lefschetz theorem. Let $A= \bigoplus\limits^{2d}_{i=0} A_i$
be a graded commutative associative algebra over a field of
characteristic zero. Let $H\in End(A)$ be a linear endomorphism
of $A$ such that for all $\eta \in A_i$, $H(\eta)= (i-d) \eta$.
This endomorphism is usually considered in Hodge theory.

For all $a\in A_2$, denote by $L_a:\;\; A\arrow A$ the linear map which
associates with $x\in A$ the element $ax\in A$. Again, this
operator is a counterpart of the operator $L$ considered in Hodge theory.
The triple $(L_a, H, \Lambda_a)\in End(A)$ is called {\bf a Lefschetz triple}
if

\[ [ L_a, \Lambda_a] = H,\;\; [ H, L_a ] = 2 L_a, \;\;
   [ H, \Lambda_a] = -2 \Lambda_a.
\]
A Lefschetz triple establishes a representation of the
Lie algebra $\goth{sl}(2)$ in the space $A$. For cohomology
algebras, this representation
arises as a part of Lefschetz theory. V. Lunts noticed that
in a Lefschetz triple, the endomorphism $\Lambda_a$ is
uniquely defined by the element $a\in A_2$
(\ref{_Lefshe_tri_unique_Proposition_}). For arbitrary
$a\in A_2$, $a$ is called {\bf of Lefschetz type} if the
Lefschetz triple $(L_a, H, \Lambda_a)$ exists. If $A= H^*(X)$
where $X$ is a compact complex manifold of K\"ahler type,
then all K\"ahler classes $\omega\in H^2(M)$ are elements
of Lefschetz type. On the other hand, not all elements of Lefschetz
type are K\"ahler classes. For instance, if $\omega$ is
of Lefschetz type, then $-\omega$ is also of Lefschetz type,
but $\omega$ and $-\omega$ cannot be K\"ahler classes
simultaneously. As one can easily check (see \cite{_Lunts-Loo_}),
the set $S\subset A_2$ of all elements of Lefschetz type
is Zariski open in $A_2$.

Now we can define the structure Lie algebra $\g(A)$ of $A$.
By definition, $\g(A)\subset End(A)$ is a Lie subalgebra
of $End(A)$ generated by $L_a$, $\Lambda_a$, for all
elements of Lefschetz type $a\in S$.
This Lie algebra is often sufficient to reconstruct the
multiplicative structure on $A$.

Returning to the hyperk\"ahler manifolds, we consider the structure Lie
algebra $\g(A)$ of the ring $A= H^*(M)$, where $M$ is a compact
hyperk\"ahler manifold. It turns out that the structure Lie algebra
$\g(A)$ can be explicitely computed. In particular,
\ref{_g(A)_for_hyperkae_Theorem_} gives us the following theorem:

\hfill

\theorem 
Let $M$ be a compact holomorphically symplectic manifold. Assume for
simplicity\footnote{The structure Lie algebra can be computed
without this trivial assumption, but the statement is
less starightforward. See the discussion after
\ref{_simple_hyperkaehler_mfolds_Definition_} for details.}
that $M$ admits a unique up to a constant holomorphic symplectic form.
Let $n= \dim(H^2(M))$.
Then $\g(A)= \goth{so}(4, n-2)$.

\hfill

This computation, which
also involves the computation of the Lie algebra $\goth M$ of
\eqref{_from_Mum-Tate_to_SO(V)_homom_Equation_},
takes up 4 sections of this paper. However, the main idea of
this computation is easy.

\hfill

Let $M$ be a compact hyperk\"ahler manifold with the complex
structures $I$, $J$, $K$. Consider the K\"ahler forms
$\omega_I$, $\omega_J$, $\omega_K$ associated with these
complex structures. Let $\rho_I:\; \goth{sl}(2)\arrow End(H^*(M))$,
$\rho_J:\; \goth{sl}(2)\arrow End(H^*(M))$,
$\rho_K:\; \goth{sl}(2)\arrow End(H^*(M))$
be the corresponding Lefschetz homomorphisms.
Let $\goth a\subset End(H^*(M))$ be the minimal
Lie subalgebra which contains images of $\rho_I$,
$\rho_J$, $\rho_K$. The algebra $\goth a$ was computed
explicitely in \cite{_so5_on_cohomo_}.

\hfill

\theorem \label{_so_5_Theorem_} 
(\cite{_so5_on_cohomo_})
The Lie algebra $\goth a$ is isomorphic to $\goth{so}(4,1)$.

\hfill

This statement can be regarded as a ``hyperk\"ahler Lefschetz theorem''.
Indeed, its proof parallels the proof of Lefschetz theorem.
One can check that the cohomology classes
$\omega_I$, $\omega_J$, $\omega_K\in H^2(M,\R)$ are orthogonal
with respect to Riemann-Hodge. Let $Hyp$ be the classifying
space of the hyperk\"ahler structures on $M$ (see Section
\ref{_moduli_Section_}). Let
$P_{hyp}:Hyp\arrow H^2(M)\times H^2(M)\times H^2(M)$
be the map which associates with the hyperk\"ahler structure
$\c H= (I, J, K, (\cdot,\cdot))$ the triple
$(\omega_I,\omega_J,\omega_K)$. Then the image of $P_{hyp}$ in
$H^2(M)\times H^2(M)\times H^2(M)$ satisfies

\begin{equation}\label{_image_of_P_hyp_Equation_}
   \forall (x,y,z)\in im P_{hyp}\;\; \bigg |\;\;
   \begin{array}{l}
   (x,y)_{\c H}=(x,z)_{\c H}=(y,z)_{\c H}=0,\\[3mm]
   (x,x)_{\c H}=(y,y)_{\c H}=(z,z)_{\c H},
   \end{array}
\end{equation}
where $(\cdot,\cdot)_{\c H}$ is the Hodge-Riemann pairing
of \eqref{_Hodge-Riemann-correct-Equation_}.
Let $D\subset H^2(M)\times H^2(M)\times H^2(M)$ be the set
defined by the equations \eqref{_image_of_P_hyp_Equation_}.
Using Torelli theorem and Calabi-Yau, we prove the following
statement:

\hfill

\theorem\label{_image_of_P_hyp_Theorem_} 
The image of $P_{hyp}$ is Zariski dense in $D$.

\hfill

\ref{_image_of_P_hyp_Theorem_} shows that all algebraic relations
which are true for \[ (x,y,z)\in P_{hyp}(Hyp) \] are true
for all $(x,y,z)\in D$. Computing the Lie algebra $\goth a$
as in \ref{_so_5_Theorem_}, we obtain a number of relations
between $x,y,z\in H^2(M)$ which hold for all $(x,y,z)\in Im(P_{hyp})$.
Using the density argument, we obtain that these relations
are universally true. This idea leads to the following theorem.

\hfill

\theorem \label{_structure_alge_for_coho_hyperkahe_Theorem_} 
Let $A= H^*(M)$ be a cohomology algebra of a compact
simply connected holomorphically symplectic manifold $M$
with $H^{2,0}(M)\cong \C$. Let $n = \dim H^2(M)$. Then
the structure Lie algebra $\g(A)$ of $A$ is isomorphic
to $\goth{so}(4,n-2)$.

\hfill

It remains to recover the multiplication on $H^*(M)$ from
the structure Lie algebra. This is done as follows. Let
$A = \oplus A_i$ be a Lefschetz-Frobenius algebra,
$\g = \g(A)$ be its structure Lie algebra. Clearly,
$\g$ is graded: $\g=\bigoplus\limits_i\g_{2i}$,
$\g_{2i}(A_j)\subset A_{j+2i}$. Let $k$ denote the
one-dimensional commutative Lie algebra. In the case
of \ref{_structure_alge_for_coho_hyperkahe_Theorem_},
$\g\cong \goth{so}(4,n-2)$, $\g_0 \cong \goth{so}(3,n-3)\oplus k$,
$\dim \g_2=\dim\g_{-2}=n$ and $\dim \g_{2i}=0$ for $|i|>1$
(see \ref{_calculation_of_g(A)_for_minim_Theorem_}).
We say that the Lefschetz-Frobenius algebra $A$ is
{\bf of Jordan type} if $\g_{2i}(A)=0$ for $|i|>1$.
For such algebras, the subspaces $\g_2(A)$,
$\g_{-2}(A)\subset A$ are commutative Lie
subalgebras of $\g(A)$. Let $U_\g$ be the universal
enveloping algebra of $\g= \g(A)$, and $U_{\g_{2}}\subset U_\g$
be the enveloping algebra of the subalgebra $\g_2=\g_2(A)\subset \g$.
Consider the space $A$ as $A$-module. Then, for all
$v\in A$ we have a map $t_v:\; U_{\g}\arrow A$ which
associates $P(v)$ with $P\in U_{\g}$. The Lie algebra
$\g_2$ is commutative, and therefore $U_{\g_2}\cong S^*(\g_2)$.
According to \cite{_Lunts-Loo_}, the natural map
$A_2\arrow \g_2$, $a \arrow L_a$,
is an isomorphism (see \ref{_g_2_is_A_2_Corollary_}
for details). Let $v\in H^0(M)\subset A$ be a unit element of the
ring $A$. Consider the restriction $t$ of $t_v:\; U_\g\arrow A$
to $U_{\g_2}\subset U_\g$. Then $t$ is a map from
$S^* \g_2\cong S^* A_2$ to $A$. Clearly, this
map coinsides with the map $S^* A_2\arrow A$ defined
by the multiplication. This implies that multiplication
by elements from $H^2(M)$ can be recovered from the action
of the structure Lie algebra $\g$.
Using the calculations of
\ref{_structure_alge_for_coho_hyperkahe_Theorem_},
we obtain, in particular, the following theorem
(see Section \ref{_cohomolo_compu_Section_}):

\hfill

\theorem \label{_S^*H^2_is_H^*M_intro-Theorem_} 
Let $M$ be a compact hyperk\"ahler manifold. Let

\[ \bar H^* (M)\subset H^*(M)
\]
be the subalgebra in $H^*(M)$
generated by $H^2(M)$. Let $\dim_\C M=2n$.
Then

\[ \bar H^{2i}(M)\cong S^i H^2(M) \]
for $i\leq n$, and
\[ \bar H^{2i}(M)\cong S^{2n-i} H^2(M) \]
for $i\geq n$.

\hfill

\hfill

\centerline{\Large \bf Contents:}

\hfill

\begin{description}
\item [Section \ref{_introduction_to_so(b_2,...)_Section_}:]\hspace{2mm}
{Introduction.}

\item [Section \ref{hyperk_manif_Section_}:]\hspace{2mm}
{Hyperk\"ahler manifolds.}

\item [Section \ref{_moduli_Section_}:]\hspace{2mm}
{Moduli spaces for hyperk\"ahler and
holomorphically symplectic manifolds.}

\item [Section \ref{_perio_and_forge_Section_}:]\hspace{2mm}
{Periods and forgetful maps.}

\item [Section \ref{_Period_and_Hodge_Riemann_Section_}:]\hspace{2mm}
{Hodge-Riemann relations for the
hyperk\"ahler manifolds and period map.}

\item [Section \ref{_Hodge-Rie_independent_Section_}:]\hspace{2mm}
{The Hodge-Riemann metric on $H^2(M)$ does not depend
on complex structure.}

\item [Section \ref{_Q_c_defini_Section_}:]\hspace{2mm}
{Period map and the space of 2-dimensional
planes in \\ $H^2(M,\R)$.}

\item [Section \ref{_Lefshe_Frob_Section_}:]\hspace{2mm}
{Lefschetz-Frobenius algebras.}

\item [Section \ref{_minimal_Fro_Section_}:]\hspace{2mm}
{The minimal Frobenius algebras and cohomology of compact
K\"ahler surfaces.}

\item [Section \ref{_^dA(V)_Section_}:]\hspace{2mm}
{Representations of $SO(V,+)$ leading to Frobenius algebras.}

\item [Section \ref{_computing_g_for_hyperk_pt-I_Section_}:]\hspace{2mm}
{Computing the structure Lie algebra for the
cohomology of a hyperk\"ahler manifold, part I.}

\item [Section \ref{_compu_g_0_part_1_Section_}:]\hspace{2mm}
{Calculation of a zero graded part of the structure Lie
algebra of the cohomology of a hyperk\"ahler manifold, part I.}

\item [Section \ref{_compu_g_0_part_2_Section_}:]\hspace{2mm}
{Calculation of a zero graded part of the structure Lie
algebra of the cohomology of a hyperk\"ahler manifold, part II.}

\item [Section \ref{_computing_g_for_hyperk_pt-2_Section_}:]\hspace{2mm}
{Computing the structure Lie algebra for the
cohomology of a hyperk\"ahler manifold, part II.}

\item [Section \ref{_cohomolo_compu_Section_}:]\hspace{2mm}
{The structure of the cohomology ring for
compact hyperkaehler manifolds.}

\item [Section \ref{_calcu_dimensi_Section_}:]\hspace{2mm}
{Calculations of dimensions.}

\end{description}

\hfill

\hfill

\begin{itemize}

\item Section \ref{_introduction_to_so(b_2,...)_Section_}
tries to supply motivations and heuristics for the further
study. In the body of this article, we never refer to Section
\ref{_introduction_to_so(b_2,...)_Section_}. Reading
Section \ref{_introduction_to_so(b_2,...)_Section_}
is not necessarily for understanding of this paper.

\item In Section \ref{hyperk_manif_Section_}, we give the definition of a
hyperk\"ahler manifold. We explain the geopetry of quaternionic action
on $H^*(M)$. This section ends with the statement of Calabi-Yau
for compact holomorphically symplectic manifolds. Results
and definitions of this section are well known.

\item Section \ref{_moduli_Section_} begins with a definition
of simple hyperk\"ahler manifolds. This notion stems from the theory
of holonomy groups. Let $M$ be a compact hyperk\"ahler manifold,
$\dim_\R M=4n$. Then the holonomy group of $M$ is naturally
a subgroup of $Sp(n)$. The simple hyperk\"ahler manifold
is a hyperk\"ahler manifold such that its
restricted holonomy group in exactly $Sp(n)$ and
not a proper subgroup of $Sp(n)$. In
\ref{_simple_hyperkaehler_mfolds_Definition_} we give a different,
but equivalent treatment of this notion
(see also \cite{_Beauville_}).Using a formalism
by de Rham and Berger, Bogomolov and Beauville proved that a
hyperk\"ahler manifold $M$ is simple if and only if
there is no finite covering $\tilde M$ of $M$ such that
$\tilde M$ can be represented as a product of two (or more)
non-trivial hyperk\"ahler manifolds. Therefore, it is
usually harmless to assume that a given compact hyperk\"ahler
manifold is simple. For simple hyperk\"ahler
manifolds, $\dim H^{2,0}(M)=1$.

\item After we define simple hyperk\"ahler manifolds, we go for
the marked moduli spaces. We define moduli spaces for complex,
hyperk\"ahler and holomorphically symplectic structures. We
do this in topological, rather than in algebro-geometrical setting.
The most obvious reason for this lopsided treatment is that
moduli of hyperk\"ahler structures don't
have any structures in addition to topology.  Results
and definitions of this section are well known.

\item In Section \ref{_perio_and_forge_Section_},
we define several kinds of period maps. Let $Comp$ be
the marked moduli of complex structures on a simple holomorphically
symplectic manifold. For all complex structures $I\in Comp$,
the space $H^{2,0}(M,I)$ of $(2,0)$-forms is one-dimensional.
With every complex structure $I\in Comp$,
the period maps $P_c:\; Comp \arrow {\Bbb P}(H^2(M, \C))$
associates a line $H^{2,0}(M,I)\subset H^2(M)$
in ${\Bbb P}(H^2(M, \C))$. With every holomorphic symplectic
structure, period map associates a class in $H^2(M,\C)$
represented by a holomorphic symplectic form. Finally, with
every hyperk\"ahler structure, period map $P_{hyp}$ associates
a triple of K\"ahler classes which correspond to the complex
structures $I$, $J$, $K$. We define a number of forgetful maps
(from hyperk\"ahler moduli to complex moduli etc.)
and compare these maps against period mappings.

\item Sections \ref{hyperk_manif_Section_} -
\ref{_perio_and_forge_Section_} contain no new results.
We establish setting for the further study of hyperk\"ahler manifolds.

\item Section \ref{_Period_and_Hodge_Riemann_Section_} answers the
following query. Hodge-Riemann pairing on cohomology satisfies a certain
type of positivity conditions, called {\bf Hodge-Riemann relations}.
What happens with these relations on a hyperk\"ahler manifold?
It turns that for a Hodge-Riemann form $(\cdot,\cdot)$ defined
on $H^2(M)$ and an action of $SU(2)$ on $H^2(M)$ which
comes from quaternions, the following conditions are satisfied:

\begin{description}

\item [{\rm (i)}] Let $H^{inv}\subset H^2(M)$ be the space
of $SU(2)$-invariants, and $H^\bot$ be its orthogonal complement.
Then the Riemann-Hodge pairing is negatively defined on
$H^{inv}$ and positively defined on $H^\bot$.

\item [{\rm (ii)}] The space $H^\bot$ is three-dimensional.

\item  [{\rm (iii)}] The space $H^\bot$ is generated
by the K\"ahler forms associated with the complex
structures $I$, $J$ and $K$.

\end{description}

\item Section \ref{_Hodge-Rie_independent_Section_}
is the crux of the first part of this paper. We prove that
the Hodge-Riemann form \eqref{_Hodge-Riemann-correct-Equation_}
on $H^2(M)$ is independent
(up to a constant) from the K\"ahler structure.
The idea of the proof is the following.

\item The group $SO(3)$ acts on the
space of the hyperk\"ahler structures, replacing the triple $(I, J, K)$
by another orthogonal triple of imaginary quaternions.
As Section \ref{_Period_and_Hodge_Riemann_Section_}
shows, this action does not change the Hodge-Riemann form.
Let $\c H \in Hyp$, and let
$(\omega_I, \omega_J,\omega_K)\in
H^2(M,\R)\times H^2(M,\R)\times H^2(M,\R)$ be the image
of $\c H$ under the action of period map $P_{hyp}$.
By definition, the Hodge-Riemann form depends only on $\omega_I$.
Therefore, we may replace $\c H$ by $\c H'\in Hyp$ such that
$P_{hyp}(\c H')= (\omega_I, \omega_J',\omega_K')$, and
such replacement does not change the Hodge-Riemann
pairing. We show that by iterating such replacements
and action of $SO(3)$, we can connect any two hyperk\"ahler
structure $\c H$ and $\c H'$ satisfying
$Vol_{\c H} (M)= Vol_{\c H'}(M)$, where by $Vol_{\c H}(M)$
we undertstand the volume of $M$ computed with respect to
the Riemannian structure associated with $\c H$.
Since these operations don't change the Riemann-Hodge form,
this form is equal for all hyperk\"ahler structure
$\c H$ of given volume. Section
\ref{_Hodge-Rie_independent_Section_} depends
on Sections
\ref{hyperk_manif_Section_} - \ref{_Period_and_Hodge_Riemann_Section_}.

\item Section \ref{_Q_c_defini_Section_} gives a description of the period
map \[ P_c:\; Comp\arrow \Bbb B(H^2(M, \C))\] in terms of the
manifold $Pl$ of 2-dimensional planes in $H^2(M, \R)$. It turns
out that there exist an etale mapping $Q_c:\; Comp\arrow Pl$.
This is a standard material, covered also in \cite{_Todorov_}.
Our exposition adds a twist to \cite{_Todorov_}, because we use
the normalized Hodge-Riemann form, which was unknown before.
Otherwise, this section depends only on Sections
\ref{hyperk_manif_Section_} - \ref{_perio_and_forge_Section_}.

\item Sections \ref{_Lefshe_Frob_Section_} - \ref{_^dA(V)_Section_}
are completely independed on the preceding sections.

\item In Section \ref{_Lefshe_Frob_Section_}, we give a number
of algebraic definitions. We give an exposition of the theory
of Lefschetz-Frobenius algebras, following \cite{_Lunts-Loo_}.
The aim of this section is a purely algebraic version
of strong Lefschetz theorem.

\item Roughly speaking, Frobenius algebra is a graded algebra
for which an algebraic version of Poincare duality holds.
A typical example of such algebra is an algebra of cohomology
of a compact manifold. Similarly, the Lefschetz-Frobenius algebra
is a Frobenius algebra for which the strong Lefschetz theorem holds -
typically, an algebra of cohomology of a K\"ahler manifold.
In Section \ref{_Lefshe_Frob_Section_} we explain these notions
and define a structure Lie algebra $\g$ of a Lefschetz-Frobenius algebra.
Here we also give a definition of Lefschetz-Frobenius algebras
of Jordan type.

\item Section \ref{_minimal_Fro_Section_} is dedicated to
explicit examples of Lefschetz-Frobenius algebras, called
{\bf minimal Lefhsetz-Frobenius algebras}.
By definition, a minimal Lefschetz-Frobenius algebra
is a Lefschetz-Frobenius algebra $A= A_0 \oplus A_2\oplus A_4$.
The ``Poincare'' form on $A$ defines a bilinear symmetric pairing on
$A_2$. It turns out that this pairing uniquely determines $A$.
Conversely, with every linear space equipped with
non-degenerate symmetric sclalar product
we associate a minimal Lefschetz-Frobenius algebra.
We prove that every minimal Lefschetz-Frobenius algebra
is of Jordan type, and explicitely compute its structure
Lie algebra. For a minimal Lefschetz-Frobenius algebra
$A(V)$ associated with a space $V$, we denote the
corresponding structure Lie algebra by $\goth{so}(V,+)$.
If $V$ is a linear space over $\R$ equipped with a scalar product
of signature $(p,q)$, then $\goth{so}(V,+)\cong \goth{so}(p+1,q+1)$.

\item In Section \ref{_^dA(V)_Section_} we find all
reduced Lefschetz-Frobenius algebras
$A= A_0\oplus A_2 \oplus ... \oplus A_n$  with the structure Lie
algebra $\goth{so}(V,+)$, where $\dim V\geq 3$. By ``reduced'' we
understand the Lefschetz-Frobenius algebras generated by $A_2$.
It turns out that for $n$ even, such algebra
is unique (we denote it by ${}^{\frac{n}{2}}A(V)$),
and for $n$ odd, there is no such algebras.

\item Sections \ref{_Lefshe_Frob_Section_} - \ref{_^dA(V)_Section_}
are purely algebraic, and Sections \ref{hyperk_manif_Section_} -
\ref{_Q_c_defini_Section_} are dealing with geometry. These
sections are mutually independent, and Sections
\ref{_computing_g_for_hyperk_pt-I_Section_} -
\ref{_computing_g_for_hyperk_pt-2_Section_} draw heavily on
both parts, geometrical and algebraic.
In these sections, we compute the structure Lie algebra
of an algebra of cohomology of a simple hyperk\"ahler manifold.
The basic result is that this algebra is isomorphic
to $\goth{so}(V,+)$, where $V$ is the linear space $H^2(M, \R)$
equipped with the normalized Hodge-Riemann pairing.

\item In Section \ref{_computing_g_for_hyperk_pt-I_Section_},
we prove that the Lefschetz-Frobenius algebra $A=H^*(M)$ is of
Jordan type: $\g(A)=\g_{-2}(A)\oplus \g_0(A)\oplus \g_2(A)$.
To prove this we introduce the standard ``density and periods''
argument, which is also used in Sections
\ref{_compu_g_0_part_2_Section_} -
\ref{_computing_g_for_hyperk_pt-2_Section_}.

\item In Section \ref{_compu_g_0_part_1_Section_}, we
construct a map $\g_0(A)\arrow \goth{so}(V)\oplus k$, where
$k$ is a one-dimensional commutative Lie algebra.
We prove that this map is an isomorphism. We also
prove that the Lie subalgebra $\goth M\subset End(V)$
generated by $ad I$ for all $I\in Comp$ is isomorphic
to $\goth{so}(V)$. By definition, $ad I:\; H^i(M)\arrow H^i(M)$
is an endomorphism which maps $\eta\in H^{p,q}_I(M)$
to $(p-q)\1 \eta$. We use the computations related to the
$\goth{so}(5)$-action on $H^*(M)$ (see \cite{_so5_on_cohomo_}).

\item In Section \ref{_compu_g_0_part_2_Section_}, we prove
that the map $u:\; \g_0\arrow \goth{so}(V)\oplus k$, constructed
in Section \ref{_compu_g_0_part_1_Section_}, is an isomorphism.
The proof is computational.

\item In Section \ref{_computing_g_for_hyperk_pt-2_Section_},
we use the results of Sections \ref{_computing_g_for_hyperk_pt-I_Section_}
- \ref{_compu_g_0_part_2_Section_} to finish the computation of
the structure Lie algebra $\g(A)$. This is done by writing down
a linear isomorphism $\g(A)\arrow \goth{so}(4,n-2)$ explicitely.
By computations, we check that this isomorphism is in fact
an isomorphism of Lie algebras.

\item Section \ref{_cohomolo_compu_Section_} is, again, algebraic.
It depends only on Sections
\ref{_Lefshe_Frob_Section_} - \ref{_^dA(V)_Section_}.
In this section, we explicitely compute the graded commutative
algebra ${}^dA(V)$ of Section \ref{_^dA(V)_Section_}.
This computation has the following
geometrical interpretation. Let $A=H^*(M)$ be the algebra of
cohomology of a simple hyperk\"ahler manifold $M$, $\dim_\R M=4d$,
and $A^r\subset A$ be its subalgebra generated by $V= H^2(M)$. The
main theorem of \ref{_^dA(V)_Section_}, together with an isomorphism
$\g(A)\cong\goth{so}(V,+)$ immediately imply
that $A^r\cong {}^dA(V)$. Therefore, by computing
${}^dA(V)$, we compute a big part of the cohomology
algebra of $M$. This way, we obtain a proof of
\ref{_S^*H^2_is_H^*M_intro-Theorem_}. The computation
of ${}^dA(V)$ is based on the classical theory of representations
of $\goth{so}(V)$ and their tensor invariants (\cite{_Weyl_}).

\item The final section applies the result of
\ref{_S^*H^2_is_H^*M_intro-Theorem_} to obtain numerical
lower bounds on Betti and Hodge numbers of a hyperk\"ahler manifold.

\end{itemize}


\section{Hyperk\"ahler manifolds.}\label{hyperk_manif_Section_}


\definition \label{_hyperkaehler_manifold_Definition_} 
(\cite{_Beauville_},
\cite{_Besse:Einst_Manifo_}) A {\bf hyperk\"ahler manifold} is a
Riemannian manifold $M$ endowed with three complex structures $I$, $J$
and $K$, such that the following holds.

\hspace{5mm}   (i)  $M$ is K\"ahler with respect to these structures and

\hspace{5mm}   (ii) $I$, $J$ and $K$, considered as  endomorphisms
of a real tangent bundle, satisfy the relation
$I\circ J=-J\circ I = K$.

\hfill

This means that the hyperk\"ahler manifold has the natural action of
quaternions ${\Bbb H}$ in its real tangent bundle.
Therefore its complex dimension is even.

Let $\mbox{ad}I$, $\mbox{ad}J$ and $\mbox{ad}K$ be the endomorphisms of
the bundles of differential forms over a hyperk\"ahler manifold
$M$ which are defined as follows. Define $\mbox{ad}I$.
Let this operator act as a complex structure operator
$I$ on the bundle of differential 1-forms. We
extend it on $i$-forms for arbitrary $i$ using Leibnitz
formula: $\mbox{ad}I(\alpha\wedge\beta)=\mbox{ad}I(\alpha)\wedge\beta+
\alpha\wedge \mbox{ad}I(\beta)$. Since Leibnitz
formula is true for a commutator in a Lie algebras, one can immediately
obtain the following identities, which follow from the same
identities in ${\Bbb H}$:

\[
   [\mbox{ad}I,\mbox{ad}J]=2\mbox{ad}K;\;
   [\mbox{ad}J,\mbox{ad}K]=2\mbox{ad}I;\;
\]

\[
   [\mbox{ad}K,\mbox{ad}I]=2\mbox{ad}J
\]

Therefore, the operators $\mbox{ad}I,\mbox{ad}J,\mbox{ad}K$
generate a Lie algebra $\goth g_M\cong \goth{su}(2)$ acting on the
bundle of differential forms. We can integrate this
Lie algebra action to the action of a Lie group
$G_M=SU(2)$. In particular, operators $I$, $J$
and $K$, which act on differential forms by the formula
$I(\alpha\wedge\beta)=I(\alpha)\wedge I(\beta)$,
belong to this group.

\hfill

\proposition \label{_there_is_action_of_G_M_Proposition_} 
There is an action of the Lie group $SU(2)$
and Lie algebra $\goth{su}(2)$ on the bundle of differential
forms over a hyperk\"ahler manifold. This action is
parallel, and therefore it commutes with Laplace operator.

{\bf Proof:} Clear.
$\blacksquare$

\hfill

If $M$ is compact, this implies that there is
a canonical $SU(2)$-action on $H^i(M,\R)$ (see
\cite{_so5_on_cohomo_}).

\hfill

Let $M$ be a hyperk\"ahler manifold with a Riemannian form
$\inangles{\cdot,\cdot}$.
Let the form $\omega_I := \inangles{I(\cdot),\cdot}$ be the usual K\"ahler
form  which is closed and parallel
(with respect to the connection). Analogously defined
forms $\omega_J$ and $\omega_K$ are
also closed and parallel.

The simple linear algebraic
consideration (\cite{_Besse:Einst_Manifo_}) shows that \hfill
$\omega_J+\sqrt{-1}\omega_K$ is of
type $(2,0)$ and, being closed, this form is also holomorphic.

\hfill

\definition \label{_canon_holo_symple_form_Definition_}
Let $\Omega:= \omega_J+\sqrt{-1}\omega_K$. This form is called
{\bf the canonical holomorphic symplectic form
of a manifold M}.

\hfill

Let $M$ be a complex manifold which admits a
holomorphic symplectic form $\Omega$. Take the Riemannian
metric $(\cdot,\cdot)$ on $M$, and the corresponding Levi-Civitta
connection. Assume that $\Omega$ is parallel
with respect to the Levi-Civitta connection.
Then the metric $(\cdot,\cdot)$ is hyperkaehler%
\footnote%
{This means that the $(\cdot,\cdot)$ is induced
by a hyperk\"ahler structure on $M$.} (\cite{_Besse:Einst_Manifo_}).

If some {\bf it compact} K\"ahler manifold $M$ admits non-degenerate
holomorphic symplectic form $\Omega$, the Calabi-Yau
(\cite{_Yau:Calabi-Yau_}) theorem
implies that $M$ is hyperk\"ahler
(\ref{_symplectic_=>_hyperkaehler_Proposition_})
This follows from the existence of a K\"ahler
metric on $M$ such that $\Omega$ is parallel under the
Levi-Civitta connection associated with this metric.

\hfill

Let $M$ be a hyperk\"ahler manifold with complex structures
$I$, $J$ and $K$. For any real numbers $a$, $b$, $c$
such that $a^2+b^2+c^2=1$ the operator $L:=aI+bJ+cK$ is also
an almost complex structure: $L^2=-1$.
Clearly, $L$ is parallel with respect to a connection.
This implies that $L$ is a complex structure, and
that $M$ is K\"ahler with respect to $L$.

\hfill

\definition \label{_induced_structures_Definion_}
If $M$ is a  hyperk\"ahler manifold,
the complex structure $L$ is called {\bf induced
by a hyperk\"ahler structure}, if $L=aI+bJ+cK$ for some
real numbers $a,b,c\:|\:a^2+b^2+c^2=1$.

\hfill

\hfill

If $M$ is a hyperk\"ahler manifold and $L$ is induced complex structure,
we will denote $M$, considered as a complex manifold with respect to
$L$, by $(M,L)$ or, sometimes, by $M_L$.

\hfill

Consider the Lie algebra $\goth{g}_M$ generated by ${ad}L$ for all $L$
induced by a hyperk\"ahler structure on $M$. One can easily see
that $\goth{g}_M=\goth{su}(2)$.
The Lie algebra $\goth{g}_M$ is called {\bf isotropy algebra} of $M$, and
corresponding Lie group $G_M$ is called an {\bf isotropy group}
of $M$. By Proposition 1.1, the action of the group is parallel,
and therefore it commutes with the action of Laplace operator on differential
forms. In particular, this implies that the action of the isotropy
group $G_M$ preserves harmonic forms, and therefore this
group canonically acts on cohomology of $M$.

\hfill

\proposition \label{_G_M_invariant_forms_Proposition_} 
Let $\omega$ be a differential form over
a hyperk\"ahler manifold $M$. The form $\omega$ is $G_M$-invariant
if and only if it is of Hodge type $(p,p)$ with respect to all
induced complex structures on $M$.

{\bf Proof:} Assume that $\omega$ is $G_M$-invariant.
This implies that all elements of $\g_M$ act trivially on
$\omega$ and, in particular, that $\mbox{ad}L(\omega)=0$
for any induced complex structure $L$. On the other hand,
$\mbox{ad}L(\omega)=(p-q)\1$ if $\omega$ is of Hodge type $(p,q)$.
Therefore $\omega$ is of Hodge type $(p,p)$ with respect to any
induced complex structure $L$.

Conversely, assume that $\omega$ is of type $(p,p)$ with respect
to all induced $L$. Then $\mbox{ad}L(\omega)=0$ for any induced $L$.
By definition, $\g_M$ is generated by such $\mbox{ad}L$,
and therefore $\g_M$ and $G_M$ act trivially on $\omega$. $\blacksquare$

\hfill

\definition \label{_degree_Definition_}
Let $M$ be a Kaehler manifold and
$\omega\in H^2(M,\R)$ be the Kaehler class of $M$.
Let $dim_\C(M)=n$. Let $\eta\in H^{2i}(M)$.
We define the {\bf degree} of the cohomology class $\eta$
by the formula

\[ deg(\eta):=
   \frac{\int\limits_M \eta\wedge\omega^{n-i}}
   {{\mbox Vol}(M)}.
\]
Clearly, if $\eta$ is pure of Hodge type $(p,q)$ and
$deg(\eta)\neq 0$, then $p=q$.

Let $M$ be a hyperkaehler manifold, and $I$ be an induced complex
structure. Then $(M,I)$ is equipped with the canonical Kaehler metric.
Consider $(M,I)$ as a Kehler manifold. We define {\bf the degree
associated with the induced complex structure $I$}
as the linear homomorphism $deg_I:\; H^{2i}(M, \R)\arrow \R$
which is equal to degree map

\[
   deg:\; H^{2i}((M,I), \R)\arrow \R
\]
defined on the cohomology of the Kaehler manifold $(M,I)$.

\hfill

The following statement follows from a trivial local
computation. The more general form of this claim is proven in
\cite{Verbitsky:Symplectic_I_}.

\hfill

\claim \label{_inv_2-forms_have_zero_degree_Claim_} 
Let $M$ be a hyperkaehler manifold and $\eta\in H^2(M)$
be a $G_M$-invariant hohomology class. Then $deg_I(\eta)=0$
for all induced complex structures $I$.

{\bf Proof:} See Theorem 2.1 of
\cite{Verbitsky:Symplectic_I_} $\;\;\blacksquare$

\hfill

Calabi-Yau theorem provides an elegant way to think of
hyperkaehler manifolds in holomorphic terms. Heuristically speaking,
compact hyperkaehler manifolds are holomorphic manifolds which admit
a holomorphic symplectic form.

\hfill

\definition \label{_holomorphi_symple_Definition_} 
The compact complex manifold $M$ is called
holomorphically symplectic if there is a holomorphic 2-form $\Omega$
over $M$ such that $\Omega^n=\Omega\wedge\Omega\wedge...$ is
a nowhere degenerate section of a canonical class of $M$.
There, $2n=dim_\C(M)$.

Note that we assumed compactness of $M$.%
\footnote{If one wants to define a holomorphic symplectic
structure in a situation when $M$ is not compact,
one should require also the equation $\nabla'\Omega$ to held.
The operator $\nabla':\;\Lambda^{p,0}(M)\arrow\Lambda^{p+1,0}(M)$
is a holomorphic differential defined on differential $(p,0)$-forms
(\cite{_Griffiths_Harris_}).}
One observes that the holomorphically symplectic
manifold has a trivial canonical bundle.
A hyperk\"ahler manifold is holomorphically symplectic
(see Section \ref{hyperk_manif_Section_}). There is a converse proposition:

\hfill

\theorem \label{_symplectic_=>_hyperkaehler_Proposition_}
(\cite{_Beauville_}, \cite{_Besse:Einst_Manifo_})
Let $M$ be a holomorphically
symplectic K\"ahler manifold with the holomorphic symplectic form
$\Omega$, a K\"ahler class
$[\omega]\in H^{1,1}(M)$ and a complex structure $I$. Assume
that

\[ deg ([\Omega]\wedge \bar [\Omega])
   =2 deg([\omega]\wedge[\omega])
\]
Then there exists a {\it unique} hyperk\"ahler structure
$(I,J,K,(\cdot,\cdot))$ over $M$ such that the cohomology
class of the symplectic form $\omega_I=(\cdot,I\cdot)$ is
equal to $[\omega]$ and the canonical symplectic form
$\omega_J+\1\omega_K$ is equal to $\Omega$.

\hfill

\ref{_symplectic_=>_hyperkaehler_Proposition_} immediately
follows from the Calabi-Yau theorem (\cite{_Yau:Calabi-Yau_}).
$\:\;\blacksquare$


\section{Moduli spaces for hyperkaehler and holomorphically symplectic
manifolds.}
\label{_moduli_Section_}


\definition \label{_simple_hyperkaehler_mfolds_Definition_} 
(\cite{_Beauville_}) The connected simply connected compact hyperkaehler
manifold
$M$ is called {\bf simple} if $M$ cannot be represented
as a Cartesian product of two (non-trivial) hyperkaehler manifolds:

\[ M\neq M_1\times M_2, \] where $M_1$, $M_2$ are
hyperkaehler manifolds  such that $dim\; M_1>0$, $dim\; M_2>0$.

\hfill

Bogomolov proved that every compact hyperkaehler manifold has a finite
covering which is a Cartesian product of a compact torus
and simple hyperkaehler manifolds. Even if our results
could be easily carried over for all compact hyperkaehler manifolds,
we restrict ourselves to the case of non-decomposable manifold
to simplify the argument.

\hfill

Let $M$ be a simple  hyperkaehler manifold.
According to Bogomolov's theorem (\cite{_Beauville_}),
for every induced complex structure $I$,
\[dim_\C \bigg(H^{2,0}((M,I))\bigg)=1.\]
This means that the space of holomorphic
symplectic forms on $(M,I)$ is one-dimensional.

{}From now on, we assume that $M$ is a simple compact hyperkaehler
manifold, which is not a torus.

\hfill

The moduli spaces of the hyperkaehler and holomorphically
symplectic manifolds were first studied by Bogomolov
(\cite{_Bogomolov_}). The studies were continued by
Todorov ([Tod]).

\hfill

Let $M_{C^\infty}$ be $M$ considered as a
differential manifold. Let $\mbox{\it Diff}$ be the group of diffeomorphisms
of $M$. Recall that the {\bf hyperkaehler structure} on
$M_{C^\infty}$ was defined as a quadruple
$(I, J, K, (\cdot, \cdot ))$ where

\[ I, J, K \in End(TM_{C^\infty}), \;\; I^2=J^2=K^2=-1 \]
are operators on the tangent bundle $TM_{C^\infty}$
and $(\cdot, \cdot )$ is a Riemannian form. This
quadruple must satisfy certain relations
(\ref{_hyperkaehler_manifold_Definition_}).
Let $\widetilde{\mbox{\it Hyp}}$ be the set of all hyperkaehler
structures on $M_{C^\infty}$. Clearly, the group $\mbox{\it Diff}$ acts
on $\widetilde{\mbox{\it Hyp}}$. The set of all non-isomorphic
hyperkaehler structures on $M_{C^\infty}$ is in bijective
correspondence with the set of orbits of $\mbox{\it Diff}$ on
$\widetilde{\mbox{\it Hyp}}$. However, the geometrical properties of
$\widetilde{\mbox{\it Hyp}}/\mbox{\it Diff}$ are not satisfactory: as a rule,
the natural topology on $\widetilde{\mbox{\it Hyp}}/\mbox{\it Diff}$
is not separable etc. To produce a more geometrical
moduli space of the hyperkaehler structures,
we will refine the space $\widetilde{\mbox{\it Hyp}}/\mbox{\it Diff}$ in
accordance with the general algebro-geometrical formalism
of marked coarse moduli spaces.

\hfill

Let $\underline{\widetilde{Comp}}$ be the set of all integrable complex
structures on  $M_{C^\infty}$. In other words,
$\underline{\widetilde{Comp}}$ is the set of all operators

\[
   I\in End(TM_{C^\infty}),\;\; I^2=-1
\]
such that the almost complex structure defined by $I$ is integrable.
The set of all non-isomorphic complex structures on $M_{C^\infty}$
is in one-to-one correspondence with
$\underline{\widetilde{Comp}}/\mbox{\it Diff}$.

For every $I\in \underline{\widetilde{Comp}}$,
we say that $I$ {\bf admits a hyperkaehler
structure} when there exist a hyperkaehler structure
$(I, J, K, (\cdot,\cdot))$ on $M$. We say that $I$
is {\bf holomorphically symplectic} when the manifold $(M, I)$
admits a non-degenerate holomorphic symplectic form.
The set $\underline{\widetilde{Comp}}$ is endowed with
a natural topology.

\hfill

Let $M$ be a compact hyperkaehler manifold.
Clearly, each of induced complex
structures on $M$ is contained in the same
connected component of $\underline{\widetilde{Comp}}$.
Denote this component by $\underline{\widetilde{Comp}}^\circ $.
\ref{_symplectic_=>_hyperkaehler_Proposition_}
immediately implies the following statement:

\hfill

\corollary \label{_defo_simple_if_it's_Kaeh_Corollary_} 
Let $I\in \underline{\widetilde{Comp}}^\circ$. Then
$(M,I)$ admits the hyperkaehler structure if and only if
$(M,I)$ admits a holomorphic symplectic structure and is
of Kaehler type%
\footnote{We say that a complex manifold $X$ is {\bf of Kaehler type}
if $X$ admits a Kaehler metric.}.

$\blacksquare$

\hfill

Denote the set of all $I\in \underline{\widetilde{Comp}}^\circ$
which admit the hyperkaehler structure by $\widetilde{Comp}^\circ$.

\hfill

The pair $(I,\Omega)$ is called {\bf the holomorphic symplectic
structure on the differential manifold
$M_{C^\infty}$} if $I$ is a complex structure
on $M_{C^\infty}$ and $\Omega$ is a holomorphic symplectic form
over $(M,I)$. Let $\widetilde{Symp}^\circ$ denote the set of
all holomorphic symplectic structures $(I,\Omega)$ on $M_{C^\infty}$
such that $I \in {\widetilde{Comp}}^\circ$.

Let $\mbox{\it Diff}\,^\circ $ be the set of all
$x\in\mbox{\it Diff}$ which act trivially on the cohomology
$H^*(M,\R)$.
Denote by $\widetilde{\mbox{\it Hyp}}^\circ $ the connected component
of $\widetilde{\mbox{\it Hyp}}$ which contains the
initial hyperkaehler structure on $M$.
Let $\mbox{\it Hyp}:=
\widetilde{\mbox{\it Hyp}}^\circ /\mbox{\it Diff}\,^\circ$,
$Symp:=\widetilde{Symp}^\circ/\mbox{\it Diff}\,^\circ$
and $Comp:= \widetilde{Comp}^\circ /\mbox{\it Diff}\,^\circ $.
These spaces are endowed with the natural topology.
Their points can be considered as the classes of
hyperkaehler (respectively, holomorphically symplectic and complex)
structures on $M_{C^\infty}$ up to the action of
$\mbox{\it Diff}\,^\circ $. Slightly abusing the
language, we will refer to these points as to hyperkaehler (resp.,
holomorphically symplectic and
complex) structures. For each $I\in Comp$,
we denote $M$, considered as a complex manifold
with the complex structure $I$, by $(M,I)$.
It is clear that $(M,I)$ is holomorphically symplectic
and admits a hyperkaehler structure for all $I\in Comp$.

\hfill

\definition \label{_Comp_Hyp_Definition_} 
The spaces $\mbox{\it Hyp}$, $Symp$, $Comp$
are called {\bf the coarse moduli spaces
of deformations of the hyperkaehler (respectively,
holomorphically symplectic and complex) structure
 on the marked compact manifold of
hyperkaehler type.}

\hfill

The word {\bf marked} refers to considering
the factorization by $\mbox{\it Diff}\,^\circ $
instead of $\mbox{\it Diff}$. This is roughly equivalent
to fixing the basis in the cohomology $H^*(M)$, hence
``marking''.


\section {Periods and forgetful maps.}\label{_perio_and_forge_Section_}


In assumptions of Section \ref{_moduli_Section_},
let $\mbox{\it Hyp}$, $Symp$, $Comp$ be the moduli spaces
of \ref{_Comp_Hyp_Definition_}.
We define the {\bf period map}

\[
   \tilde P_{hyp}: \;\widetilde{\mbox{\it Hyp}}\arrow H^2(M,\R)\otimes \R^3
\]
as a rule which associates with every hyperkaehler
structure $(I,J,K,(\cdot,\cdot))$ on $M_{C^\infty}$
the triple

\[ ([\omega_I],\: [\omega_J],\: [\omega_K])\in
   H^2(M_{C^\infty},\R)\times H^2(M_{C^\infty},\R)
   \times H^2(M_{C^\infty},\R)
\]
of Kaehler classes corresponding to $I$, $J$ and
$K$ respectively.
By definition, the group $\mbox{\it Diff}\,^\circ $
acts trivially on $H^2(M)$. Therefore, $\tilde P_{hyp}$
descends to a map

\[
   P_{hyp}: \;\mbox{\it Hyp}\arrow H^2(M,\R)\otimes \R^3
\]

Similarly, define the Griffiths'
period map

\[
  P_c:\; Comp\arrow {\Bbb P}^1(H^2(M_{C^\infty},\C))
\]
as a rule which relies the 1-dimensional complex subspace
\[
  H^{2,0}((M_{C^\infty},I))\subset  H^2(M_{C^\infty},\C)
\]
to the complex structure $I$ on $M_{C^\infty}$. Using Dolbeault
spectral sequence, one can easily see that
the subspace $H^{2,0}((M_{C^\infty},I))\subset  H^2(M_{C^\infty},\C)$
is defined independently on the Kaehler metric.

Let $P_s:\; Symp\arrow H^2(M,\C)$
map a pair $(I,\Omega)\in Symp$ to the class
$[\Omega]\in H^2(M,\C)$ which is represented by the
closed 2-form $\Omega$.

There exist a number of natural ``forgetful maps'' between
the spaces $\mbox{\it Hyp}$, $Symp$ and $Comp$.
Here we define some of these maps and find how these
maps relate to period maps.

Let $\c H= (I,J,K, (\cdot,\cdot))\in \mbox{\it Hyp}$ be a hyperkaehler
structure. As in \ref{_canon_holo_symple_form_Definition_},
consider the canonical holomorphic symplectic form
$\Omega:= \omega_J+\1 \omega_K$ associated with $\c H$.
Let $\Phi^{hyp}_s:\; \mbox{\it Hyp}\arrow Symp$ map
$\c H$ to the pair $(I, \Omega)\in Symp$. Let
$\Phi^{hyp}_s:\; \mbox{\it Hyp}\arrow Comp$ map
$\c H$ to $I\in Comp$. Let $\Phi^s_c:\; Symp\arrow Comp$
map $\c S=(I, \Omega)\in Symp$ to $I\in Comp$.

For $h=(x_1,x_2,x_3)\in H^2(M,\R)\times H^2(M,\R)\times H^2(M,\R)$,
let $\pi_i(h)=x_i$.

\hfill

\claim \label{_forgetting-n-periods_Claim_} 
Let $\c H\in Hyp$, $\c S\in Symp$.
Then

(i) $P_s(\Phi^{hyp}_s(\c H))= \pi_2(P_{hyp}(\c H))+
\1\pi_3(P_{hyp}(\c H))$

(ii) The point $P_c(\Phi^{s}_c(\c S))\in {\Bbb P} H^2(M,\C)$
corresponds to a line generated by $P_s(\c S)$.

{\bf Proof:} Clear. $\;\;\blacksquare$

\hfill

According to the general formalism of Kodaira and Kuranishi, $Comp$ is
endowed with a canonical structure of a complex variety.
Using the complex structure on $Comp$, we describe the map
$P_c^s:\; Symp\arrow Comp$ in terms of algebraic geometry.

Let $L$ be a holomorphic vector bundle over a complex variety $X$.
Let $Tot(L)$ be the total space of  $L$.
By definition, $Tot(L)$ is a complex variety which is
smoothly fibered over $X$. Every holomorphic section
$f\in \Gamma_X(L)$ gives a standard holomorphic map
$s_f:\; X\arrow Tot(L)$. Consider the map $s_0:\; X\arrow Tot(L)$
corresponding to the zero section of $L$. This map identifies
$X$ with the closed analytic subspace of $Tot(L)$.
Let $Tot^*(L):= Tot(L)\backslash X$ be the completion
of $Tot(L)$ to $X$.

\hfill

\proposition \label{_Symp_as_a_total_space_Proposition_} 
Let $M$ be a hyperkaehler manifold, $Comp$ and $Symp$ be the
moduli spaces associated with $M$ as in \ref{_Comp_Hyp_Definition_}.
Then there exist a natural holomorphic linear bundle $\tilde \Omega$
on $Comp$ such that the following conditions hold.

(i) There exist a natural homeomorphism
$i:\; Tot^*(\tilde \Omega)\arrow Symp$.

(ii) Let $\pi:\; Tot^*(\tilde \Omega)\arrow Comp$ be the
standard projection. Then the diagram

\[
   \begin{array}{ccccc}
      Tot^*(\tilde \Omega)&\!\!\!\stackrel {i}\arrow \!\!\!& Symp\\[3mm]
      \;\;\;\;\;\searrow\!\!{}^\pi\!\!\!\!\!
      && \!\!\!\swarrow\!\!{}_{\Phi^s_c}\;\;\;\;\; \\[3mm]
     & \!\!\!\!\!Comp &\\
   \end{array}
\]
is commutative.

\hfill

The homeomorphism $i$ defines a complex analytic structure on $Symp$.
Further on, we consider both $Comp$ and $Symp$ as complex analytic
varieties.

\hfill

{\bf Proof of \ref{_Symp_as_a_total_space_Proposition_}:}
Let ${A}:= Comp\times H^2(M, \C)$.
The holomorphic symplectic structure $(I,\Omega)\in Symp$
is uniquely defined by $I\in Comp$ and the cohomology class
$[\Omega]\in H^2(M, \C)$. This defines an injection

\[ j:\; Symp\hookrightarrow {A}, \;\;
   j(I,\Omega)= (I, [\Omega]).
\]
Consider ${A}$ as the total space of a trivial holomorphic
bundle ${A}_b$ with the fiber $H^2(M,\C)$.
We construct $\tilde \Omega$ as a linear subbundle of
${A}_b$, such that its total space coinsides with
$j(Symp)$.

\hfill

Let $U\subset Comp$ be an open set. We say that
{\bf there exists a universal fibration over $U$}
if there exist a smooth complex analytic fibration
$\pi:\;\goth M\arrow U$ such that for all $J\in U$,
the fiber $\pi^{-1}(J)$ is isomorphic to $(M,J)$.

\hfill

\claim \label{_unive_fibra_exi_loca_Claim_} 
For all $I\in Comp$, there exist an open set $U\subset Comp$,
$I\in U$, which admits universal fibration.

{\bf Proof:} This is a consequence of Kodaira-Spencer theory
(see \cite{_Kodaira_Spencer_}). $\;\;\blacksquare$

\hfill

Let $U\subset Comp$ be an open subset which admits a universal fibration
$\goth M\stackrel \pi \arrow U$. Let $\C_{\goth M}$ be the
constant sheaf over $\goth M$, and $\pi_\bullet$ be the sheaf-theoretic
direct image. Let $H^2:=R^2\pi_\bullet \C_{\goth M}$ be the second derived
functor of $\pi_\bullet$ applied to $\C_{\goth M}$. Since $\C_{\goth M}$
is a constructible sheaf, and $\pi$ is a proper morphism,
the sheaf $H^2$ is also constructible.
For every point $I\in U$, the restriction $H^2\restrict{I}$
is isomorphic to $H^2(M,\C)$. Hence, $H^2\restrict{I}$ is a locally
constant sheaf. This sheaf is equipped with a natural flat connection,
known as Gauss-Manin connection. Since $U$ is a
subset in the space of {\it marked} deformations of $M$,
the monodromy of Gauss-Manin connection is trivial. Hence,
the bundle $H^2$ is naturally isomorphic to ${A}_b\restrict{U}$.

Let $F^0\subset F^1\subset F^2={A}_b\restrict{U}$ be the variation of
Hodge structures associated with $\pi$. By definition,
$F^i$ are holomorphic sub-bundles of ${A}_b\restrict{U}$.
For every $I\in U$, we have $F^0\cong H^{2,0}((M,I))$.
Therefore, $F^0$ is a linear sub-bundle of ${A}_b\restrict{U}$.

\hfill

\lemma \label{_hol_bundle_from_Hodge_str_and_j_Lemma_} 
Let $U\subset Comp$ be an open set which admits an universal
fibration $\goth M\arrow U$. Let $F^0\subset {A}_b\restrict{U}$
be the holomorphic linear bundle defined as above.
Let

\[
   Symp(U):= \{ (I,\Omega)\in Symp\;\;|\;\; I\in U\}.
\]
Let $Tot(F^0)\subset U\times H^2(M,\R)$ be the total space
of $F^0$ considered as a subspace in a total space of
${A}_b\restrict{U}$. Then $Tot^*(F^0)$ coinsides with $j(Symp(U))$

{\bf Proof:} Clear. $\;\;\blacksquare$

\hfill

For every open set which admits an universal
fibration $\goth M\arrow U$, we defined the
holomorphic linear bundle $F^0\subset {A}_b\restrict{U}$.
\ref{_hol_bundle_from_Hodge_str_and_j_Lemma_} implies that
$F^0\subset {A}_b\restrict{U}$ is independent from the choice of the
universal fibration, and that the locally defined
sub-bundles $F^0$ can be glued to a globally defined
holomorphic linear sub-bundle in ${A}_b$. Denote this
linear sub-bundle by $\tilde \Omega$. It is clear that
$Tot^*(\tilde \Omega)$ coinsides with $j(Symp)$.
Since $j$ is injective, there exist an inverse homomorphism
$i:\; Tot^*(\tilde \Omega)\arrow Simp$. The condition
(ii) of \ref{_Symp_as_a_total_space_Proposition_} is
obvious. It remains to show that the bijective maps $i$ and
$j= i^{-1}$ are continuous. This is left to the reader
as an exercise.  $\;\;\blacksquare$


\section[Hodge-Riemann relations for the
hyperkaehler manifolds and period map.]
{Hodge-Riemann relations for the \\
hyperkaehler manifolds and period map.}
\label{_Period_and_Hodge_Riemann_Section_}


\subsection{What do we do in this section:} 
With every hyperkaehler manifold $M$, we associate
the action of the group $G_M\cong SU(2)$ on the cohomology
of $M$ (see \ref{_there_is_action_of_G_M_Proposition_}).
Let $P_{hyp}(M):= (\omega_I,\omega_J,\omega_K)$ be the periods
of $M$. We show that the action of $G_M$ may be reconstructed
from the periods. This follows from \ref{_Lambda_dual_to_L_Proposition_}
and \ref{_g_m_from_L_Lambda_Claim_}.

The action of $G_M\cong SU(2)$ on $H^2(M)$ induces a weight
decomposition of $H^2(M)$. Using this decomposition, we obtain
an interesting version of Hodge-Riemann relations
(\ref{_restrictions_of_pairings_to_H^2_Lemma_}). In
particular, we obtain that for every hyperkaehler structure
$(I,J,K,(\cdot,\cdot))$ on $M$, the Hodge-Riemann
pairings associated with the complex structures $I$, $J$ and
$K$ are equal (\ref{_pairings_on_H^2_are_equal_Proposition_}).

\subsection{Hodge-Riemannian pairing.} 

\vspace{3mm}

\hspace{5.5mm}In this subsection we follow \cite{_Weil_}.

\hfill

Let $X$ be a compact Kaehler manifold, and
$\Lambda^*(X)=\oplus \Lambda^{p,q}(X)$
be a space of differential forms equipped with Hodge decomposition.
The Riemannian structure on $X$ equips $\Lambda^*(X)$ with a positively
defined Hermitian metric (see \cite{_Weil_} for correct normalization of
this metric). Integrating the scalar product of two forms over $X$,
we obtain a Hermitian positively defined pairing on the space
of global sections of $\Lambda^*(X)$. Let us identify the cohomology
space of $X$ with the space of harmonic differential forms. This gives
a positively defined Hermitian product on the space of cohomology
$H^*(X)$. We denote it by

\[ (\cdot,\cdot)_{Her}:\; H^i(X,\C)\times H^i(X,\C) \arrow \C. \]

Let $I:\; H^i(X,\C)\times H^i(X,\C) \arrow \C$ map $(x,y)$ to
$(x,\bar y)_{Her}$. Clearly, $I$ is a complex-linear non-degenerate
2-form on $H^i(X,\C)$, which is defined over reals.
Let $A:\;  H^i(X,\C)\times H^i(X,\C) \arrow \C$ map
$x,y \in H^i(X,\C)$ to

\[ \int_X x\wedge y\wedge \omega^{n-i},
\]
where $n=\dim_\C X$, and $i\leq n$. Let $C:\; H^*(X,\C)\arrow H^*(X,\C)$
be the Weil operator, which maps a cohomology class
$\omega\in H^{p,q}(X)\subset H^{p+q}(X)$ to $\1^{p-q} \omega$.
Let $L:\; H^i(X)\arrow H^{i+2}(X)$,
$\Lambda:\; H^i(X)\arrow H^{i-2}(X)$
be the Hodge operators, and $P^i(X)\subset H^i(X)$
be the space of primitive cohomology classes:

\[ P^i(X) = \{ \alpha \in H^i(X) \;\; |\;\; \Lambda(\alpha) \} =0 \]

By Lefshetz theorem,

\[
   H^i(X) = \oplus L^r P^{i-r}(X).
\]
Let $p_r:\; H^i(X) \arrow P^{i-r}(X)$ be a map corresponding to
this decomposition, such that for all $a\in H^*(M)$,

\[ a = \sum_r L^r p_r(a).  \]

The forms $A$ and $I$ are related by the so-called
Hodge-Riemann equation:

\begin{equation}\label{_Hodge_Riemann_general_Equation_}
   (-1)^\frac{(n-i)(n-i-1)}{2} A(a,C b) =
   \sum_r \mu_r\frac{(n-p+r)!}{r!} I(L^r p_r(a), L^r p_r(b)),
\end{equation}
where $\mu_r$ are positive real constants which depend only on
$r$ and dimension of $X$.

Let $\omega\in H^2(X, \R)$ be a Kaehler class of $X$.
We call the form

\begin{equation}\label{_Hodge_Riemann_form_general_Equation_}
  (-1)^\frac{(n-i)(n-i-1)}{2} A(a,C b):
   \; H^i(X,\C)\times H^i(X,\C) \arrow \C
\end{equation}
{\bf the Hodge-Riemann pairing associated with a Kaehler class
$\omega$} and
denote this form by $(\cdot,\cdot)_{\omega}$.

\hfill

\claim \label{_Hodge_Riema_general_Claim_} 

(i) The form $(\cdot,\cdot)_{\omega}$
depends only on the Kaehler class $\omega\in H^2(X,\R)$ of $X$.
In other words, $(\cdot,\cdot)_{\omega}$
would not change if we modify the complex structure or Kaehler
metric, provided that the Kaehler class stays the same.

(ii) The form $(\cdot,\cdot)_{\omega}$ is defined over reals.

(iii) If $X$ is a surface, then the restriction of $(\cdot,\cdot)_{\omega}$
to the primitive cohomology $P^2(X)\subset H^2(X)$
coincides with the intersection form.

(iv) Restriction of $(\cdot,\cdot)_{\omega}$ to $H^2(X)$ can be written
as follows:

\begin{equation} \label{_H_R_to_H^2_formula_Equation_}
   (\eta_1,\eta_2)_{\omega}=\int_X \omega^{n-2}\wedge \eta_1\wedge\eta_2
   - \frac{n-2}{(n-1)^2} \cdot \frac{
   \int_X \omega^{n-1}\eta_1 \cdot \int_X\omega^{n-1}\eta_2}
   {\int_X \omega^n}.
\end{equation}

\hfill

{\bf Proof:} Follows from \eqref{_Hodge_Riemann_general_Equation_}
(see also \cite{_Weil_}). $\;\;\blacksquare$

\hfill

\subsection{Riemann-Hodge relations in hyperkaehler case.} 

Let $M$ be a compact hyperkaehler manifold and $\c H\in Hyp$ be
a hyperkaehler structure on $M$. Let $P_{hyp}(\c H)$ be
denoted by $(x_1, x_2, x_3)$, $x_i\in H^2(M,\R)$.
The hyperkaehler structure $\c H\in Hyp$ defines a Riemannian
metric on $M$. This metric establishes a positively
defined Hermitian scalar product on the space of $\C$-valued
differential forms over $M$. Realizing the cohomology classes
as harmonic forms, we obtain a Hermitian pairing
$(\cdot,\cdot)_{Her}$ on  $H^i(M,\C)$.
Let $\Lambda_{x_i}\!:\; H^i(M)\arrow H^{i-2}(M)$ be
the operator adjoint to $L_{x_i}$ with respect to the
pairing $(\cdot,\cdot)_{Her}$. Clearly,
$\Lambda_{x_i}$ is the Hodge operator associated with the Kaehler
structure on $M$ which is defined by $x_i$ and ${\cal H}$.
Let $\inangles{\cdot,\cdot}_{x_i}$ be the Hodge-Riemann form
\eqref{_Hodge_Riemann_form_general_Equation_} associated
with the Kaehler form $x_i$. Let $L\galochka_{x_i}$ be an
operator adjoint to $L_{x_i}$ with respect to $\inangles{\cdot,\cdot}_{x_i}$.

\hfill

\proposition \label{_Lambda_dual_to_L_Proposition_} 
\ \   $L\galochka_{x_j}=\Lambda_{x_j}$ \ \ \  for $j=1,2,3$.

{\bf Proof:} Fix a choice of $j\in\{1,2,3\}$. For simplicity,
assume that $j=1$. We abbreviate $L_{x_1}$ by $L$, $\Lambda_{x_1}$
by $\Lambda$.
Let $I$ be the complex structure
induced by ${\cal H}$, such that $x_1\in H^2(M,\R)$
is the Kaehler class of $I$. Take the Lefschetz decomposition
(\cite{_Griffiths_Harris_})

\[ H^k(M)=\bigoplus L^i P^{k-2i}(M),\;\;\;
   P^{k-2i}\subset H^{k-2i}(M),
\]
where the space $P^i(M)$ is a space of all primitive
classes:

\[ P^i(M)=
   ker\bigg(\Lambda:\;  H^i(M)\arrow H^{i+2}(M)\bigg).
\]

Let $P^{p,q}(M):= P^{p+q}(M)\cap H^{p,q}(M)$.
It is well known that

\[ P^{i}(M)=\bigoplus\limits_{p+q=i}P^{p,q}(M). \]
(see \cite{_Griffiths_Harris_}, \cite{_Weil_}).
Hodge-Riemann relations (\eqref{_Hodge_Riemann_general_Equation_};
see also \cite{_Weil_})
describe $(\zeta_1,\zeta_2)_{Her}$ in terms of
$\inangles{\zeta_1,\bar\zeta_2}_{x_1}$ and Lefschetz
decomposition. Let $n=dim_\C(M)$.
When

\[ \zeta_1\in L^i P^{p,q}(M), \;\; \zeta_2\in L^{i'} P^{p',q'}(M)
   \;\;\;\mbox{\it and} \;\;\;(i,p,q)\neq (i',p',q'),
\]
both scalar products vanish:

\begin{equation}\label{_Hodge_Riemann_relations_vanishing_Equation}
  (\zeta_1,\zeta_2)_{Her}=
  \inangles{\zeta_1,\bar\zeta_2}_{x_1}=0.
\end{equation}

\hfill

When $\zeta_1\in L^i P^{p,q}(M)$ and
$\zeta_2\in L^i P^{p,q}(M)$, we have

\begin{equation}\label{_Hodge_Riemann_relations_Equation_}
  (\zeta_1,\zeta_2)_{Her}=
  \1^{p-q}(-1)^{\frac{(n-p-q)(n-p-q-1)}{2}}
  \inangles{\zeta_1,\bar\zeta_2}_{x_1}.
\end{equation}

\hfill

The operator $L\galochka_{x_1}$ is adjoint to $L_{x_1}$ with respect
to $\inangles{\cdot,\cdot}_{x_1}$ and $\Lambda_{x_1}$
is adjoint to $L_{x_1}$ with respect to $(\cdot,\cdot)_{Her}$.
Let $\zeta\in  L^i P^{p,q}(M)$. By definition,
$t=L\galochka_{x_1}(\zeta)$ is the element of
$H^{2i-2+p+q}(M)$ such that $\forall \xi \in H^{2i-2+p+q}(M)$
we have

\begin{equation} \label{_L_galochka_definition_Equation_}
   \inangles{t,\xi}_{x_1} =
   \inangles{\zeta,L_{x_1}(\xi)}_{x_1}.
\end{equation}

Using \eqref{_Hodge_Riemann_relations_vanishing_Equation},
we see that $t\in L^{i-1} P^{p,q}(M)$. For
$t\in L^{i-1} P^{p,q}(M)$,
\eqref{_Hodge_Riemann_relations_vanishing_Equation}
shows that if \eqref{_L_galochka_definition_Equation_}
holds for all $\xi\in L^{i-1} P^{p,q}(M)$,
this equation holds for all $\xi \in H^{2i-2+p+q}(M)$.
On the other hand, Hodge-Riemann relations
\eqref{_Hodge_Riemann_relations_Equation_}
imply that for $\xi,t\in  L^{i-1} P^{p,q}(M)$,

\[ (t,\bar \xi)_{Her}=
   \1^{p-q}(-1)^{\frac{(n-p-q)(n-p-q-1)}{2}}
   \inangles{t,\xi}_{x_1} =
\]
\[
   = \1^{p-q}(-1)^{\frac{(n-p-q)(n-p-q-1)}{2}}
   \inangles{\zeta,L_{x_1}(\xi)}_{x_1}=
   (\zeta,L_{x_1}\bar\xi)_{Her}.
\]
Therefore $L\galochka_{x_1}$ is adjoint to
$L_{x_1}$ with respect to $(\cdot,\cdot)_{Her}$.
\ref{_Lambda_dual_to_L_Proposition_} is proven.
$\;\;\blacksquare$

\hfill

In Section \ref{hyperk_manif_Section_},
we defined the action of $G_M\cong SU(2)$
on the cohomology of a hyperkaehler manifold.
For every hyperkaehler structure ${\cal H}\in \mbox{\it Hyp}$,
there is an action of $SU(2)$ on $H^*(M_{C^\infty})$
which is determined by  ${\cal H}$. We proceed to describe this
$SU(2)$-action in terms of the triple

\[ (x_1,x_2,x_3)=P_{hyp}({\cal H})\in
   H^2(M)\oplus H^2(M)\oplus H^2(M).
\]
\ref{_Lambda_dual_to_L_Proposition_} expresses
the Hodge operators $\Lambda_{x_i}$ via the multiplicative
structure on $H^*(M_{C^\infty})$. Let
$\goth a_{\cal H}\subset End(H^*(M,\R)$ be the
Lie algebra generated by $L_{x_i}$, $\Lambda_{x_i}$, $i=1,2,3$.
According to \cite{_so5_on_cohomo_},
$\goth a_{\cal H}\cong \goth{so}(4,1)$.

\hfill

Let $\g_{\cal H}\subset \goth a_{\cal H}$ be the subalgebra of
$\goth a_{\cal H}$ consisting of all elements which
respect the grading on $H^*(M)$ induced by the degree:

\[ \g_{\cal H}:=
   \{ x\in \goth a_{\cal H}\;\; | \;\; x(H^i(M))\subset H^i(M),
   \;\;i=0,1, ... \; 2n.\}
\]

\hfill

\claim \label{_g_m_from_L_Lambda_Claim_} 
The Lie algebra $\g_{\cal H}$ is isomorphic to
$\goth{su}(2)$. Its action coincides with that
of $\g_M$ defined in Section \ref{hyperk_manif_Section_}.

{\bf Proof:} This is Theorem 2 of
\cite{_so5_on_cohomo_}. $\;\;\blacksquare$

\hfill

The forms $(\cdot,\cdot)_{x_i}$, $i=1,2,3$ depend only on the
value of $x_i\in H^2(M,\R)$. Therefore,
\ref{_Lambda_dual_to_L_Proposition_} has the following
interesting consequence:

\hfill

\corollary \label{_so(5)_inde_from_H_Corollary_} %
Let $\c H\in Hyp$ be a hyperkaehler structure
on $M$, and $\g_{\cal H}\cong \goth{su}(2)$,
$\goth a_{\c H}\cong \goth{so}(4,1)$ be the corresponding
 Lie algebras defined as above. Then the action of
$\g_{\cal H}$, $\goth a_{\c H}$ on $H^*(M,\R)$ depends only
on hyperkaehler periods of $\c H$.
In other words, if $\c H_1$,
$\c H_2$ are hyperkaehler structures such that
\[P_{hyp}(\c H_1)= P_{hyp}(\c H_2),\] then action
of $\goth a_{\c H_1}$, $\g_{\c H_1}$ on $H^*(M)$ coinsides
with action of $\goth a_{\c H_2}$, $\g_{\c H_2}$.

$\;\;\blacksquare$

\hfill

\proposition \label{_pairings_on_H^2_are_equal_Proposition_} 
In assumptions of \ref{_Lambda_dual_to_L_Proposition_}, let
$\inangles{\cdot,\cdot}_i$ be the restriction of the
pairing $\inangles{\cdot,\cdot}_{x_i}$ to $H^2(M)$.
Then $\inangles{\cdot,\cdot}_1 =
\inangles{\cdot,\cdot}_2=\inangles{\cdot,\cdot}_3$.

{\bf Proof:} Let $(\cdot,\cdot)$ be the restriction
of $(\cdot,\cdot)_{\cal H}$ to $H^2(M)$.
Let $V$ be the subspace of $H^2(M)$ spanned by
$(x_1,x_2,x_3)$.

Earlier, we defined the action of the Lie algebra
$\g_{\cal H}\cong\goth{su}(2)$ on $H^2(M)$.
Let $H_{inv}$ be the space of all $\g_{\cal H}$-invariant
elements in $H^2(M)$. According to \cite{_so5_on_cohomo_},
the action of $\g_{\cal H}$ on $H^*(M)$ induces the Hodge
decomposition. Namely, for every induced complex structure
$I$ there exist a Cartan subalgebra $\goth h\in \g_{\cal H}$
such that the weight decomposition on $H^*(M)$ induced by
$\goth h$ coincides with the Hodge decomposition

\[ H^i(M)=\bigoplus\limits_{p+q=i}H^{p,q}(M). \]
The space $H^{2,0}(M)$ is one-dimensional
for every induced complex structure. Using the theory
of representations of $\goth{sl}(2)$, one can check that
this implies that $H^2(M)/H_{inv}$ is a simple
3-dimensional representation of $\g_{\cal H}$.
Another trivial calculation shows that
$V$ is a $\g_{\cal H}$-invariant subspace of
$H^2(M)$, and $\g_{\cal H}$ acts on $V$ non-trivially.
Therefore $H^2(M)=H_{inv}\oplus V$.
\ref{_pairings_on_H^2_are_equal_Proposition_}
is implied by the following lemma.

\hfill

\lemma \label{_restrictions_of_pairings_to_H^2_Lemma_} 
Consider the restrictions of $\inangles{\cdot,\cdot}_i$
and $(\cdot,\cdot)$ to $V$ and $H_{inv}$. Then

\begin{equation} \label{_restriction_to_V_Equation_}
  \inangles{x,\bar y} = (x,\bar y)
  \;\;\;\mbox{for} \;\;x,y\in V,
\end{equation}
\begin{equation} \label{_restriction_to_H_inv_Equation_}
   \inangles{x,\bar y}_i =
   -(x,\bar y)\;\;\;\mbox{for}\;\; x,y\in H_{inv}, \;\:
   \mbox{\it\  and finally,}
\end{equation}
\begin{equation} \label{_restriction_to_H_and_V_Equation_}
   \inangles{x,\bar y}_i =
   (x,\bar y)=0\;\;\;\mbox{for} \;\;x\in V,\; y\in H_{inv}.
\end{equation}

\hfill

{\bf Proof:} Let $I$ be an induced complex structure.
Then $V=H^{2,0}\oplus L_I (H^{0,0})\oplus H^{0,2}(M)$.
where the Hodge decomposition is taken with respect
to $I$ and $L_I$ is the Hodge operator of exterrior
multiplication by the Kaehler class of $I$. Then
\eqref{_restriction_to_V_Equation_} immediately follows
from Hodge-Riemann relations
\eqref{_Hodge_Riemann_relations_Equation_}.
By \ref{_inv_2-forms_have_zero_degree_Claim_}
all elements of $H_{inv}$ are primitive.

On the other hand,
$H_{inv}\subset H^{1,1}$ by
\ref{_G_M_invariant_forms_Proposition_}. Therefore
\eqref{_restriction_to_H_inv_Equation_}
follows from
\eqref{_Hodge_Riemann_relations_Equation_},
and \eqref{_restriction_to_H_and_V_Equation_}
follows from \eqref{_Hodge_Riemann_relations_vanishing_Equation}.
\ref{_restrictions_of_pairings_to_H^2_Lemma_}
and consequently \ref{_pairings_on_H^2_are_equal_Proposition_}
are proven.
$\;\;\blacksquare$

\hfill

\corollary \label{_indu_comple_same_HR_Corollary_} 
Let  $\c H= (I, J, K, (\cdot,\cdot))$ be a hyperkaehler
structure on $M$, and $L=a I + bJ + cK$ be an induced
complex structure, $a^2+b^2+c^2=1$. Let $\omega_1\in H^2(M,\R)$
be the Kaehler class associated with the Kaehler manifold
$(M, I)$, and $\omega\in H^2(M,\R)$ be the Kaehler
class associated with the Kaehler manifold $(M, L)$.
Let $(\cdot,\cdot)_\omega$,
$(\cdot,\cdot)_{\omega_1}:\; H^2(M)\times H^2(M)\arrow \C$
be the Hodge-Riemann forms associated with $\omega$, $\omega_1$.
Then $(\cdot,\cdot)_{\omega_1}=(\cdot,\cdot)_\omega$

{\bf Proof:} Follows from \ref{_restrictions_of_pairings_to_H^2_Lemma_}
$\;\;\blacksquare$


\section{The Hodge-Riemann metric on $H^2(M)$ does not depend
on complex structure.}
\label{_Hodge-Rie_independent_Section_}


Let $M_{C^\infty}$ be a compact manifold which
admits a hyperkaehler structure. Let $Hyp$, $Symp$, $Comp$ be the moduli
spaces constructed in Section \ref{_moduli_Section_}.

\hfill

\definition \label{_Hodge_Riemann_asso_w_hyperkeahler_Definition_} 
Let $M$ be a hyperkaehler manifold, $\c H =(I, J, K, (\cdot, \cdot))$
be its hyperkaehler structure and $\omega_1$, $\omega_2$,
$\omega_3\in H^2(M, \R)$ be Kaehler
classes associated with induced complex structures $I$, $J$, $K$.
Consider the Riemann-Hodge pairing
$(\cdot,\cdot)_{\omega_i}:\; H^2(M,\R)\times H^2(M, \R)\arrow \R$, $i=1,2,3$
defined as in \eqref{_Hodge_Riemann_form_general_Equation_}
(see also \eqref{_H_R_to_H^2_formula_Equation_}).
According to \ref{_pairings_on_H^2_are_equal_Proposition_},

\[ (\cdot,\cdot)_{\omega_1}= (\cdot,\cdot)_{\omega_2} =
   (\cdot,\cdot)_{\omega_3}.
\]
Let
\[ \inangles{x,y}_{\c H} := (x,y)_{(vol(M))^{-1/n}\cdot x_i}
\]
where the volume $vol(M)= \int_M x_i^n$ is volume
calculated with respect to the Riemannian metric
$(\cdot,\cdot)$, and $n=\frac{\dim_\R(M)}{2}$.
This pairing is called {\bf the normalized Hodge-Riemann pairing
associated with the hyperkaehler structure $\c H$}.
According to \eqref{_H_R_to_H^2_formula_Equation_},
the normalized Hodge-Riemann pairing
$\inangles{\cdot,\cdot}_{\c H}$ can be expressed by

\begin{equation} \label{_normalized_HR_Equation_}
   (\eta_1,\eta_2)_{\omega}=
   \lambda^{n-2}\int_X \omega^{n-2}\wedge \eta_1\wedge\eta_2
   - \frac{n-2}{(n-1)^2}
   \lambda^{2n-2}\int_X \omega^{n-1}\eta_1 \cdot \int_X\omega^{n-1}\eta_2
\end{equation}
where $\lambda = (vol(M))^{-1/n}$.

\hfill

The main result of this section is the following theorem:

\hfill

\theorem \label{_Hodge_Riemann_independent_Theorem_} 
Let $\c H_1$, $\c H_2\in \mbox{\it Hyp}$ be hyperkaehler structures on
$M_{C^\infty}$. Then
$\inangles{\cdot, \cdot}_{\c H_1} =
   \inangles{\cdot, \cdot}_{\c H_2}$
In other words, the normalized Hodge-Riemann pairing
associated with the point $\c H\in \mbox{\it Hyp}$ does not depend
on the choice of $\c H$ in $Hyp$.

\hfill

{\bf Proof:} The space $\mbox{\it Hyp}$ is endowed with the homogenous
action of the group $SO(3)$ as follows. Let
$(I, J, K, (\cdot, \cdot))\in \mbox{\it Hyp}$ be a hyperkaehler structure.
Consider the complex structures $I$, $J$, $K$ as endomorphisms
of the tangent bundle $TM$. We express this by
$I, J, K\in \Gamma_M(End(TM))$. Consider the
three-dimensional subspace $V\subset \Gamma_M(End (TM))$ generated
by $I$, $J$, $K$. By definition of a hyperkaehler structure,
$V$ is a three-dimensional vector space equipped with a canonical
isomorphism with the space of anti-self-adjoint quaternions.
The space $\c V$ of anti-self-adjoint quaternions is a Lie subalgebra
of the quaternion algebra. The space of sections
$\Gamma_M(End (TM))$ has a canonical algebra structure. By definition,
$V$ is a Lie subalgebra of $\Gamma_M(End(TM))$, and the Lie
algebra structure on $V$ coincides with that on $\c V\subset \Bbb H$.
The Lie algebra $\c V\subset \Bbb H$ is isomorphic to $\goth{so}(3)$.
Consider the adjoint action of $SO(3)$ on $V\cong \c V\cong \goth {so}(3)$.
Let $A\in SO(3)$. By definition of adjoint action,

\begin{equation} \label{_A_of_quate_quate_Equation_}
\begin{array}{l}
   A(I)^2=A(J)^2=A(K)^2, \mbox{ \ and \ } \\
   A(I)\circ A(J)=A(K)=-A(J)\circ A(I)
\end{array}
\end{equation}
The operators $I$, $J$, $K$ are parallel with respect
to the Levi-Civita connection. The operators $A(I)$, $A(J)$, $A(K)$
are linear combinations of $I$, $J$, $K$ with constant coefficients.
Hence, these operators are also parallel. They are orthogonal
by trivial reasons.

\hfill

\claim \label{_Newla_Niere_for_para_Claim_} 
Let $X$ be a Riemannian
manifold equipped with Levi-Civita connection. Let $\c I$ be an
orthogonal almost complex structure which is parallel with respect
to the connection. Then $\c I$ is an integrable (i. e., defines a
complex structure).
Moreover, the Riemannian metric on $X$ is Kaehler.

{\bf Proof:} This follows from Newlander-Nierenberg theorem.
$\;\;\blacksquare$

\hfill

\ref{_Newla_Niere_for_para_Claim_} implies that
$A(I)$, $A(J)$, $A(K)$ are operators of complex structure.
Now, \eqref{_A_of_quate_quate_Equation_} implies that
$(A(I), A(J), A(K),  (\cdot, \cdot))$ is a hyperkaehler structure.
We obtain an action of $SO(3)$ on $Hyp$.

\hfill

\definition \label{_action_SO(3)_on_Hyp_Definition_} 
This action of $SO(3)$ on $\mbox{\it Hyp}$ is called {\bf a standard
action of $SO(3)$ on the space of hyperkaehler structures}.
Two hyperkaehler structures are called {\bf equivalent}
if one can be obtained from another by the standard action
of $SO(3)$.

\hfill

It is easy to describe the action of $SO(3)$ on $\mbox{\it Hyp}$
in terms of the period map:

\hfill

\claim\label{_action_SO(3)_on_Hyp_via_periods_Lemma_} 
Let $\c H \in \mbox{\it Hyp}$ and $A\in SO(3)$. Consider
$P_{hyp}(\c H)$ and $P_{hyp}(A(\c H))$ as elements
of the space

\[ W:= H^2(M,\R)\otimes\R^3
   \cong H^2(M,\R)\oplus H^2(M,\R)\oplus H^2(M,\R).
\]
Then $P_{hyp}(A(\c H))$ is obtained from $P_{hyp}(\c H)$
by applying $Id\otimes A$ to $P_{hyp}(\c H)\in H^2(M,\R)\otimes\R^3$.

{\bf Proof:} Clear. $\;\;\blacksquare$

\hfill

\lemma \label{_Hodge_Riemann_independe_for_equiva_hyperkae_Lemma_} 
In assumptions of \ref{_Hodge_Riemann_independent_Theorem_},
Let $\c H_1$, $\c H_2\in \mbox{\it Hyp}$ be
equivalent hyperkaehler structures on
$M_{C^\infty}$. Then

\[ \inangles{\cdot, \cdot}_{\c H_1} =
   \inangles{\cdot, \cdot}_{\c H_2}
\]

{\bf Proof:} Follows from
\ref{_indu_comple_same_HR_Corollary_}. $\;\;\blacksquare$

\hfill

\definition \label{_admissible_substi_Definition_} 
Let $\c H_1$, $\c H_2\in \mbox{\it Hyp}$. We say that $\c H_2$
{\bf is obtained from $\c H_1$ by an admissible substitution}
if either of the following two conditions hold:

(i) There exists $\lambda\in \R$ such that
$P_3(\c H_1)= \lambda P_3(\c H_2)$.

(ii) $\c H_1$ is equivalent to $\c H_2$.

We say that $\c H_1$ and $\c H_2$ are {\bf well connected}
if $\c H_2$ can be obtained from $\c H_1$ by a sequence of
admissible substitutions. Obviously, this relation is an equivalence relation.

\hfill

\lemma \label{_Hodge_Riemann_independe_for_well_connected_Lemma_} 
In assumptions of \ref{_Hodge_Riemann_independent_Theorem_},
Let $\c H_1$, $\c H_2\in \mbox{\it Hyp}$ be the hyperkaehler structures on
$M_{C^\infty}$. Assume that $\c H_1$ and $\c H_2$ are well connected.
Then

\[   \inangles{\cdot, \cdot}_{\c H_1} =
   \inangles{\cdot, \cdot}_{\c H_2}.
\]

{\bf Proof:} It is sufficient to prove
\ref{_Hodge_Riemann_independe_for_well_connected_Lemma_}
assuming that $\c H_2$ is obtained from $\c H_1$ by admissible
substitution. In other words, we may assume that one
of conditions (i) and (ii) of \ref{_admissible_substi_Definition_}
holds. When (i) holds,
\ref{_Hodge_Riemann_independe_for_well_connected_Lemma_}
follows from
\ref{_Hodge_Riemann_independe_for_equiva_hyperkae_Lemma_}. When (ii)
holds, \ref{_Hodge_Riemann_independe_for_well_connected_Lemma_}
is a direct consequence of \eqref{_normalized_HR_Equation_}
(see also \ref{_Hodge_Riema_general_Claim_} (iv)).
$\;\;\blacksquare$

\hfill

We obtain that \ref{_Hodge_Riemann_independent_Theorem_}
is a consequence of
\ref{_Hodge_Riemann_independe_for_well_connected_Lemma_}
and the following statement:

\hfill

\proposition \label{_hyperk_are_well_connected_Proposition_} 
Let $\c H_1$, $\c H_2\in \mbox{\it Hyp}$ be the hyperkaehler structures on
$M_{C^\infty}$. Then $\c H_1$ and $\c H_2$ are well connected.

{\bf Proof:}

\hfill

\lemma \label{_hyp.st._w/the_same_I_w/conne_Lemma_} 
Let  $\c H=(I, J, K, (\cdot,\cdot))$
and $\c H'=(I', J', K', (\cdot,\cdot)')$ be two hyperkaehler
structures with $I=I'$. Then $\c H$ and $\c H'$ are well
connected.

{\bf Proof:} Since $I=I'$, we have

\[
   P_c(\Phi^{hyp}_c(\c H))=P_c(\Phi^{hyp}_c(\c H')).
\]
By \ref{_forgetting-n-periods_Claim_}, the spaces spanned by
$\inangles{P_2(\c H), P_3(\c H)}$ and $\inangles{P_2(\c H'), P_3(\c H')}$
coinside. Denote $\inangles{P_2(\c H), P_3(\c H)}$ by $W$.
Let $U$ be the group of linear automorphisms
of $W$ which preserve the Hodge-Riemann pairing
$(\cdot,\cdot)_{\c H}$ associated with $\c H$. Using the
basis $W=\inangles{P_2(\c H), P_3(\c H)}$, we may identify
$U$ with $U(1)$. Let $u\in U\cong U(1)$ be represented by
the matrix

\[ u = \bigg(\begin{array}{rr}
                \cos(\alpha) & \sin(\alpha) \\
                -\sin(\alpha) & \cos(\alpha)
             \end{array}
       \bigg).
\]
Let $u(J)= \cos(\alpha) J+ \sin(\alpha) K$ and
$u(K)= \cos(\alpha) K- \sin(\alpha) J$. Checking the
definition of the hyperkaehler structure,
one obtains that $u(\c H):= (I, u(J), u(K), (\cdot, \cdot))$
is a hyperkaehler structure which is equivalent to
$\c H$. By \ref{_action_SO(3)_on_Hyp_via_periods_Lemma_},
$P_2(u(\c H))= u(P_2(\c H))$ and $P_3(u(\c H))= u(P_3(\c H))$.
Choosing a suitable $u$, we can make $P_3(u(\c H))$
proportional to $P_3(\c H')$. For such $u$, $u(\c H)$
is well connected to $\c H'$. Since $u(\c H)$ is equivalent
to $\c H$, we obtain that $\c H$ is well connected to $\c H'$.
$\;\;\blacksquare$

\hfill

\lemma\label{_hyper_w/same_ind_comp_str_well_connect_Corollary_}
Let $\c H_1$, $\c H_2\in \mbox{\it Hyp}$ be the hyperkaehler structures, and
$I\in Comp$ be the complex structure. Assume that $I$ is induced by
$\c H_1$ and $\c H_2$. Then $\c H_1$ is well connected with $\c H_2$.

{\bf Proof:} Clearly,
\ref{_hyper_w/same_ind_comp_str_well_connect_Corollary_}
follows from \ref{_hyp.st._w/the_same_I_w/conne_Lemma_}
and the following statement:

\hfill

\sublemma \label{_induced_compl_str_turn_to_I_Sublemma_} 
Let $(I, J, K, (\cdot,\cdot))=\c H\in \mbox{\it Hyp}$,
$I'\in Comp$ be a complex structure which
is induced by $\c H$. Then $\c H$ is equivalent to a hyperkaehler
structure $\c H' = (I', J', K', (\cdot,\cdot)')$ for some
$J', K', (\cdot,\cdot')$.

{\bf Proof:} Consider the action of $SO(3)$ on the space
$V:=\inangles{I,J,K}\subset \Gamma(End(TM))$ (see
\ref{_action_SO(3)_on_Hyp_Definition_}). Clearly,
$I'$, considered as a section of $\Gamma(End(TM))$, belongs
to $V$. Take a matrix $A\in SO(3)$ which maps
$I\in V$ to $I'$. Then $\c H':=A(\c H)$ satisfies conditions
of \ref{_induced_compl_str_turn_to_I_Sublemma_}.
$\;\;\blacksquare$

\hfill

\definition \label{_well_conne_comple_str_Definition_}
Let $I_1$, $I_2\in Comp$. The complex structures $I_1$, $I_2$
are called {\bf well connected} if there exist well connected
hyperkaehler structures $\c H_1$, $\c H_2$ such that
$\c H_1$ induces $I_1$ and $\c H_2$ induces $I_2$.
By \ref{_hyper_w/same_ind_comp_str_well_connect_Corollary_},
this is an equivalence relation.

\hfill

Let $\omega\in H^2(M,\R)$. Let $Comp^\omega$ be the set
of all $I\in Comp$ such that $\omega$ belongs to the Kaehler
cone of $I$.

\hfill

\claim \label{_Comp^omega_well_conne_Claim_} 
Let $\omega\in H^2(M,\R)$, $I, I'\in Comp^\omega$. Then
the complex structures $I$ and $I'$ are well connected.

{\bf Proof:} Clear. $\;\;\blacksquare$

\hfill

Let $\mbox{\it Kah}$ be the set of all $\omega\in H^2(M,\R)$ such that
$\omega$ is a Kaehler class for some complex structure $I\in Comp$.

\hfill

\definition \label{_well_conne_kah_str_Definition_}
Let $\omega$, $\omega'\in \mbox{\it Kah}$. The classes  $\omega$ and $\omega'$
are called {\bf well connected} if there exist well connected
hyperkaehler structures $\c H$, $\c H'$ such that
$P_1(\c H)= \omega$ and $P_1(\c H')= \omega'$.

\hfill

\lemma \label{_well_conne_Kah_classe_indu_w/c_hype_Lemma_} 
Let $\omega$, $\omega'$ be two well connected classes from $\mbox{\it Kah}$.
Let $\c H$, $\c H'$ be the hyperkaehler structures such that
$P_1(\c H)= \omega$ and $P_1(\c H')= \omega'$. Then
$\c H$, $\c H'$ are well connected.

{\bf Proof:} Consider the spaces $Comp^\omega$ and $Comp^{\omega'}$.
Since $\omega$ is well connected to $\omega'$, there exist
well connected hyperkaehler structures
\[ \c F=\bigg(A, B, C, (\cdot,\cdot)\bigg),\;\; \;
   \c F'=\bigg(A', B', C', (\cdot,\cdot)'\bigg)
\]
such that
$P_1(\c F)= \omega$ and $P_1(\c F')= \omega'$. Therefore,
$A\in Comp^\omega$ is well connected to $A'\in Comp^{\omega'}$.
Let $\c H = (I, J, K, (\cdot,\cdot))$,
$\c H'=(I', J', K', (\cdot,\cdot)')$. By definition,
$I\in Comp^\omega$ and  $I'\in Comp^{\omega'}$. By
\ref{_Comp^omega_well_conne_Claim_}, $I$ is well
connected to $A$ and $I'$ is well connected to $A'$.
Therefore, $I$ is well connected to $I'$. By definition
of well connected complex structures, there exist
well connected hyperkaehler structures $\c G$ and $\c G'$
such that $\c G$ induces $I$ and  $\c G'$ induces $I'$
By \ref{_hyper_w/same_ind_comp_str_well_connect_Corollary_},
$\c G$ is well connected to $\c H$ and $\c G'$ is well
connected to $\c H'$. Since the relation of being well connected
is transitive, $\c H$ is well connected to $\c H'$.
\ref{_well_conne_Kah_classe_indu_w/c_hype_Lemma_}
is proven. $\;\;\blacksquare$

\hfill

To finish the proof of \ref{_hyperk_are_well_connected_Proposition_},
it is sufficient to show that for all $x,y \in \mbox{\it Kah}$, the cohomology
classes $x$ and $y$ are well connected. By
\ref{_symplectic_=>_hyperkaehler_Proposition_},
$\mbox{\it Kah}=P_1(Hyp)$. Since $\mbox{\it Hyp}$ is connected, $\mbox{\it
Kah}$ is also connected.
Therefore \ref{_hyperk_are_well_connected_Proposition_} is
implied by the following lemma:

\hfill

\lemma \label{_C(omega)_open_in_Kah_Lemma_} 
Let $\omega\in \mbox{\it Kah}$, $C(\omega)$ be the set of all
classes $\omega'\in \mbox{\it Kah}$ which are well connected to $\omega$.
Then $C(\omega)$ is open in $\mbox{\it Kah}$.

{\bf Proof:} Since $\mbox{\it Kah}\subset H^2(M,\R)$,
it is sufficient to show that $C(\omega)$ is open in $H^2(M,\R)$.
Since the relation of being well connected is transitive,
it is sufficient to show that $C(\omega)\subset \mbox{\it Kah}$
contains an open neighbourhood of $\omega$ for all $\omega\in Kah$.
Let $I\in Comp^\omega$. Let $\c H$ be a hyperkaehler
structure associated with $I$ and $\omega$
as in \ref{_symplectic_=>_hyperkaehler_Proposition_}.

Let $K(\c H)\subset Kah$ be the set of all $\eta\in H^2(M,\R)$
such that the following condition holds. The hyperkaehler structure
$\c H$ induces a complex structure $L$ such that $\eta\in K(L)$.
As usually, $K(L)$ is the Kaehler cone of $L$.
As the following lemma implies, $K(\c H)\subset C(\omega)$.

\hfill

\sublemma \label{_K(H)_well_conne_to_omega_Sublemma_} 
Let $\omega'\in K(\c H)$. Then $\omega'$ is well connected to $\omega$.

{\bf Proof:} Let $\c H\in Hyp$, $\omega=P_1(\c H)$.
Let $L$ be an induced complex structure such that
$\omega'\in K(L)$. By \ref{_induced_compl_str_turn_to_I_Sublemma_},
there exist a hyperkaehler structure $\c H'=(L, J, K, (\cdot,\cdot))$
which is equivalent to $\c H$. Then, $\c H$ is well connected to $\c H'$.
Let $\c H''$ be the hyperkaehler structure associated with $L$ and $\omega'$
as in \ref{_symplectic_=>_hyperkaehler_Proposition_}. By
\ref{_hyper_w/same_ind_comp_str_well_connect_Corollary_},
$\c H'$ and $\c H''$ are well connected. By definition,
$P_1(\c H'')=\omega'$. Since $\c H$ is well connected to $\c H''$,
$\omega$ is well connected to $\omega'$.
$\;\;\blacksquare$

\hfill

Let $\c H \in Hyp$ and $L$ be an induced complex structure.
As usually, we denote the intersection $H^{1,1}(M,L)\cap H^2(M,\R)$
by $H^{1,1}_L(M,\R)$. Let $\omega\in H^{1,1}_L(M,\R)$.
The hyperkaehler structure $\c H$ induces a Riemannian metric on $M$.
Let $\tilde \omega\in \Lambda^2(M,\R)$ be the harmonic
form which represents the cohomology
class $\omega$. Hodge theory implies that $\tilde\omega$ is
a form of Hodge type (1,1) with respect to the complex structure
$L$. Under these assumptions, we introduce the following
definition.

\hfill

\definition \label{_positive_classes_Definition_} 
We say that the cohomology class $\omega$ is {\bf positive}
with respect to $(\c H,L)$ if the corresponding harmonic
(1,1)-form $\tilde \omega$ is everywhere positively defined.
In other words, $\omega\in H^{1,1}_L(M,\R)$ is {\bf positive}
if the symmetric form

\[ S_p: \; T_pM\times T_pM \arrow \R,\;\;
   S_p(x,y):= \tilde \omega(x, L(y))
\]
is positively defined in every point of $p\in M$.

We denote by $K_{\c H}(L)$
the set of all $\omega\in H^{1,1}_L(M,\R)$
such that $\omega$ is positive with respect to $(\c H,L)$.

\hfill

\claim \label{_positive_form_is_Kaehler_Claim_} 
In assumptions of \ref{_positive_classes_Definition_},
let $\omega\in H^{1,1}_L(M,\R)$ be the two-form which is positive
with respect to $(\c H,L)$. Then $\omega$ is a Kaehler class:
$\omega\in K(L)$.

{\bf Proof:} Clear. $\;\;\blacksquare$

\hfill

Let $\c H\in Hyp$. Denote by $K_c(\c H)$ the set of all
$\omega\in H^2(M,\R)$ such that there exists an induced
complex structure $L$ and $\omega\in K_{\c H}(L)$.
By \ref{_positive_form_is_Kaehler_Claim_},
$K_c(\c H)\subset K(\c H)$. Let $\omega\in C(\omega)$.
This means that $\omega=P_1(\c H)$ for some $\c H\in Hyp$.
Since $K_c(\c H)\subset K(\c H)\subset C(\omega)$, to prove that
$C(\omega)$ is open in $H^2(M,\R)$ it is sufficient to show that
$K_c(\c H)$ contains an open neighbourhood of $\omega$
(we use here the transitiveness of well-connectedness).
Therefor, \ref{_C(omega)_open_in_Kah_Lemma_} is
a consequence of the following statement:

\hfill

\proposition \label{_K(H)_open_in_H^2(M,R)_Sublemma_} 
Let $M$ be a hyperkaehler manifold with the hyperkaehler
structure $\c H$. Then the set $K_c(\c H)\subset H^2(M,\R)$
contains an open neighbourhood of $P_1(\c H)$.

{\bf Proof:} Consider the action of the group of unit quaternions
$G_M\cong SU(2)$ defined as in Section
\ref{hyperk_manif_Section_}. The action of $G_M$ is defined on the
tangent bundle $T(M)$. We naturally extend this action to the tensor
powers of $T(M)$, including $End(T(M))\cong T(M)\otimes T^*(M)$.
Consider the set $R$ of induced complex structures as subset of
the space of sections $\Gamma_M(End(TM))$. An easy local computation
shows that $G_M$ acts transitively on $R\cong S^2$ (see
also \ref{_induced_compl_str_turn_to_I_Sublemma_}). Let $L$ be an
induced complex structure, $\omega\in K_{\c H}(L)$. Let $g\in G_M$,
$L':=g(L)$. Consider the Kaehler form $\omega$ as the section of
$\Lambda^2(TM)\subset T^*(M)\otimes T^*(M)$. Obviously, the 2-form

\[ \inbfpare{\cdot,\cdot}:= g(\omega)(L'(\cdot),\cdot)=
   \omega(g \circ g^{-1}\circ L\circ g (\cdot),g(\cdot))
   = \omega (L(g(\cdot)),g(\cdot))
\]
is symmetric and positively defined. To show that
the Riemannian form $\inbfpare{\cdot,\cdot}$
is Kaehler, we have to prove that the form

\[ \inbfpare{L'(\cdot),\cdot} = - g(\omega)(\cdot,\cdot) \]
is symplectic. Since $G_M$ commutes with Laplacian,
it maps harmonic forms to harmonic ones. Hence, $g(\omega)$ is
a symplectic form. Therefore $\inbfpare{\cdot,\cdot}$ is a
Kaehler metric. This implies that $g(\omega)\in K(L')$.

We proved the following statement:

\hfill

\claim \label{_G_M_acts_transi_on_K(H)_Claim_} 
Let $M$ be a hyperkaehler manifold with the hyperkaehler
structure $\c H$. Consider the action of the group of unit
quaternions $G_M\cong SU(2)$ on $H^2(M,\R)$
(see \ref{_there_is_action_of_G_M_Proposition_}).
Let $g\in G_M$, $L, L'\in R$, $L'=g(L)$. Then
$g:\; H^2(M, \R)\arrow H^2(M,\R)$ induces
an isomorphism from $K_{\c H}(L)\subset H^2(M,\R)$ to
$K_{\c H}(L')\subset H^2(M,\R)$.
$\;\;\blacksquare$

\hfill

The following statement is clear:

\hfill

\claim \label{_K_c(L)_is_open_Claim_} 
In assumptions of \ref{_positive_classes_Definition_},
the set $K_{\c H}(L)$ is open in $H^{1,1}_L(M,\R)$.

$\blacksquare$

\hfill

\ref{_G_M_acts_transi_on_K(H)_Claim_}
establishes a smooth map $\delta:\; G_M\times K_{\c H}(I)\arrow K_c(\c H)$.
Let $b_2:=dim(H^2(M))$. By \ref{_K_c(L)_is_open_Claim_},
$dim_\R K_{\c H}(I)= dim_\C (H^{11}((M, I))$.
Since $M$ is a simple hyperkaehler manifold,
$dim_\C (H^{11}((M, I)) = b_2-2$. According to Section
\ref{hyperk_manif_Section_}, $R\cong S^2$. Therefore,
$dim_\R(R\times K_{\c H}(I))= b_2$. By definition, $G_M$ is
identified with the group of unit quaternions, and $R$
is identified with the set $x\in G_M \ \ | \ \ x^2=-1$.
This identification defines a canonical embedding $R\hookrightarrow G_M$.
Let $\phi:\; R\times K_{\c H}(I)\arrow K_c(\c H)$ be the restriction
of $\delta$ to $R\times K_{\c H}(I)\subset G_M\times K_{\c H}(I)$.
The dimension of $R\times K_{\c H}(I)$ is equal to $b_2$.
Therefore, to prove that $K_c(\c H)$ is open
in $H^2(M,\R)$ it is sufficient to prove the following:

\hfill

\sublemma \label{_K(H)cong_S^2_times_K(I)_Sublemma_} 
The map $\phi:\; R\times K_{\c H}(I)\arrow K_c(\c H)$ is
a diffeomorphism.

{\bf Proof:}
As in \ref{_restrictions_of_pairings_to_H^2_Lemma_},
consider the decomposition $H^2(M,\R)= H_{inv}\oplus V$.
As we have established previously, $V$ is generated
by $P_i(\c H)$, $i=1,2,3$ and $(H_{inv},V)_{\c H}=0$.
For $x\in H^2(M, \R)$, let $\pi_i(x)$ be the orthogonal projection
of $x$ to $H_{inv}$ and $\pi_v(x)$ be the orthogonal projection
of $x$ to $V$. The bilinear form $(\cdot,\cdot)_{\c H}$ is $G_M$-invariant by
\ref{_Hodge_Riemann_independe_for_equiva_hyperkae_Lemma_}.
Therefore, $\pi_i(g(x))=g(\pi_i(x))$.
For every induced complex structure $L$, $\c H$
defines a Kaehler structure on the complex manifold $(M, L)$.
Hence, for every induced complex structure $L$, the hyperkaehler
structure $\c H$ defines a Kaehler form $\omega_L$ and a
degree map $deg_L:\; H^{2i}(M,\R)\arrow \R$. According to
\ref{_inv_2-forms_have_zero_degree_Claim_}, for all
$x\in H_{inv}$ and all induced complex structures $L$,
$deg_L(x)=0$. Therefore for all $x\in K_c(\c H)$, we have
$\pi_v(x)\neq 0$.

Let $y\in K_c(\c H)$. Let $l(y)$ be the line
in the three-dimensional space $V$ generated by $\pi_v(y)$.
The space $V$ is generated by the set of induced complex structures,
which constitute a unit sphere in $V$. Hence, the space
of lines in $V$ is canonically identified with the set of
complex structures up to a sign.
Let $R^{\pm}\cong \R P^2$ be the quotient of $R$ by $\pm 1$.
Let $\theta: \; K_c(\c H)\arrow R^{\pm}$ map $y$ to the
point of $R^{\pm}$ which corresponds to $l(y)$.
Denote the induced complex structures which correspond
to $\theta(y)$ by $L_1$, $L_2$, where $L_2=-L_1$.

Denote the Hodge decomposition
associated with an arbitrary complex structure $L\in Comp$ by $H^{pq}_L$.
According to \ref{_G_M_invariant_forms_Proposition_},
$x\in H^{pp}_L$ if and only if $L(x)=x$.

Obviously,
$L(y)=\pi_i(y)+ L(\pi_v(y))$. Realizing $L$ and $\pi_v(y)$
as quaternions in a usual way, we may check that
$L(\pi_v(y))= L\pi_v(y)L^{-1}$. Since the centralizator
of all elements in $SU(2)$ is one-dimensional, $L\in l$
whenever $L(y)=y$. Therefore for $y\in H^{11}_L$, we
have $L=\pm L_1$. Since $\pi_v$ is an orthogonal projection,

\begin{equation} \label{_pi_v(y)_Equation_}
   \pi_v(y)=\frac{deg_L(y)}{deg_L(\omega_L)}\omega_L,
\end{equation}
where $\omega_L$ is the Kaehler form of $(M,L)$ considered
as an element of $V$, where

\[ deg_L(y):= \int_M \omega_L^{n-1}\wedge y. \]

By definition, $deg_{L_1}(x)=-deg_{L_2}(x)$.
On the other hand, for $x\in K(L)$, we have $deg_L(x)>0$.
Therefore for all $y\in K_c(\c H)$, there exist only
one induced complex structure $L$ such that $y\in K(L)$.
This implies that $\phi$ is a monomorphism. We need to construct the inverse
of $\phi$ and prove that it is smooth. Let
$\Theta^{\pm}: K_c(\c H) \arrow R^{\pm}\times (H_{inv}\oplus \R)$
map $y\in K_c(\c H)$ to the pair

\[ \bigg(\theta (y),\; \pi_i(y)\oplus |deg_{\pi_i(y)}(y)|\bigg) \]
where $deg_{\pi_i(y)}$ is well defined up to a sign. According
to \eqref{_pi_v(y)_Equation_}, up to a sign, one can reconstruct
$\pi_v(y)$ by $\Theta^{\pm}(y)$. Therefore,
$\Theta^{\pm}$ is a double covering. Let
$\rho:\;R^{\pm}\times (H_{inv}\oplus \R)\arrow H^2(M,\R)$ map
$(s,h+t)\in R^{\pm}\times(H_{inv}\oplus \R)$ to

\[ \frac{t}{deg_I(\omega_I)}\omega_I +h. \]

The $I$-degree of $\omega\in K_{\c H}(I)$ is positive.
Therefore \eqref{_pi_v(y)_Equation_} implies that
the map

\[ \phi\circ \Theta^{\pm} \circ \rho:\;
   R\times K_c(\c H)\arrow R^{\pm}\times K_c(\c H)
\]
acts as identity on $K_c(\c H)$ and acts as a double covering on
$R$. Since $\Theta^{\pm}$ is a double covering, $\phi$ is an
open embedding. \ref{_K(H)cong_S^2_times_K(I)_Sublemma_}
and \ref{_C(omega)_open_in_Kah_Lemma_} is proven. The proof of
\ref{_hyperk_are_well_connected_Proposition_} and
consequently \ref{_Hodge_Riemann_independent_Theorem_}
is finished. $\;\;\blacksquare$


\section{Period map and the space of 2-dimensional
planes in $H^2(M,\protect \R)$.} \label{_Q_c_defini_Section_}


There is an alternative way of looking at Griffiths period map
$P_c:\; Comp\arrow {\Bbb P}(H^2(M,\C))$. This enhanced version
of period map is a map from $Comp$ to an open subset in Grassmanian
of all 2-dimensional planes in $H^2(M, \R)$. To define this
map, we remind the reader certain well-known
results from linear algebra.

Let $V_\R$ be an $\R$-linear space endowed with the
non-degenerate symmetric bilinear form $(\cdot,\cdot)_\R$.
Let $V_\C:= V_\R\otimes\C$ be the complexification of
$V_\R$, and $(\cdot,\cdot)_\C$ be the $\C$-linear form on
$V_\C$ obtained as a complexification of $(\cdot,\cdot)_\R$.

In applications, $V_\R= H^2(M,\R)$, $V_\C= H^2(M,\C)$,
and $(\cdot,\cdot)_{\R}$ is the normalized
Hodge-Riemann pairing $(\cdot,\cdot)_{\c H}$.

\hfill

Consider the projectivization ${\Bbb P} V_\C$ as a space
of lines in $V_\C$. For all $x\in V_\C$, let $\bar x$ denote
the complex conjugate to $x$. Let

\[ C:= \{ t\in V_\C \;\; |\;\; \forall x\in t,\, (x,x)_\C=0,
       (x,\bar x)_\C> 0 \}.
\]
Let $\tilde Pl$ be the space of all oriented 2-dimensional
linear subspaces in $V_\R$. Let $Pl\subset \tilde Pl$
be the set of all $L\in \tilde Pl$ such that the restriction
of $(\cdot,\cdot)_\R$ to the 2-dimensional space $L\subset V_R$
is positively defined. Clearly, $Pl$ is open in $\tilde Pl$.

For $t\in {\Bbb P} V_\C$, take $x\in t$, $x\neq 0$.
Let $i_x(t)\subset V_\R$ be the linear span of
$Re(x)$, $Im(x)\in V_\R$. If $i_x(t)$ is two-dimensional,
we consider $i_x(t)$ as the oriented space with the
orientation defined by the basis $(Re(x), Im(x))$.

\hfill

\proposition \label{_from_PV_C_to_Pl_linear-alg_Proposition_} 
Let $t\in C\subset {\Bbb P} V_\C$. Then the space $i_x(t)$
is 2-dimensional and independent on the choice of $x\in t$.
Let $i:\; C\arrow \tilde Pl$ map $t\in C$ to $i_x(t)$.
The image of $i:\; C\arrow \tilde Pl$
coinsides with $Pl\subset \tilde Pl$. Established this
way map $i:\; C\arrow Pl$ is bijective.

{\bf Proof:} Let $t\in C$, $x\in t\subset V_\C$. Let
$y=Re(x), \;z= Im(x)$. Then

\[
   (x,x)_\C= (y,y)_\R-(z,z)_\R + 2\1 (y,z)_\R= 0.
\]
Therefore $(y,y)_\R=(z,z)_\R$ and $(y,z)_\R=0$. On the other hand,

\[
   (x,\bar x)_\C = (y,y)_\R+ (z,z)_\R> 0.
\]
We obtain that

\begin{equation} \label{_y^2=z^2>0_Equation_}
  (y,y)_\R = (z,z)_\R > 0 \;\; \mbox{and}\;\; (y,z)_\R=0
\end{equation}
This implies that the vectors $y$ and $z$ are linearly independent.
For all $c=\lambda e^{\1 \alpha}\in \C$, where $\lambda,\alpha\in \R$,
we have

\[ Re(cx)=\lambda \cos (\alpha) y +\lambda \sin(\alpha) z, \;
   Im(cx)=-\lambda \sin (\alpha) y +\lambda \cos(\alpha) z.
\]
Therefore, $i_x(t)=i_{cx}(t)$. This implies that the map
$i:\; C\arrow \tilde Pl$ is well defined. According to
\eqref{_y^2=z^2>0_Equation_}, $i(C)\subset Pl$.

\hfill

Let us construct the inverse map $j:\; Pl\arrow C$. For
$L\in Pl$, take an oriented orthonormal basis $(y,z)$ in the
Euclidean space $L$. For another orthonormal basis
$(y',z')$ in $L$, we have

\begin{equation} \label{_rotation_on_y,z_Equation_}
    \begin{array}{rrrrr}
       y'& =& \cos (\alpha) y &+&\sin(\alpha) z \\[3mm]
       z'& = & -\sin (\alpha) y &+&\cos(\alpha) z
    \end{array}
\end{equation}
for some $\alpha\in \R$. Let $j(L)\in {\Bbb P}V_\C$
be the line generated by $y+ \1 z\in V_\C$. The equations
\eqref{_rotation_on_y,z_Equation_} imply that $j(L)$ is independent
of the choice of the oriented orthonormal basis $(y,z)$.
Since $(y,z)$ is orthonormal basis, \eqref{_y^2=z^2>0_Equation_}
holds. This equation immediately implies that
$j(L)\in C\subset {\Bbb P}V_\C$.
Finally, it is clear from defintions
that $j\circ i= Id$ and $i\circ j= Id$.
$\;\;\blacksquare$

\hfill

Returning to the hyperkaehler manifolds, consider the space

\[
   \goth C\subset {\Bbb P}(H^2(M,\C))
\]
consisting of lines
$l\in H^2(M,\C)$ such that for all $x\in l$, $(x,x)_{\c H}=0$
and $(x,\bar x)_{\c H}>0$.
Here, as elsewhere,

\[
   (\cdot,\cdot)_{\c H}:\; H^2(M,\C)\times H^2(M,\C)\arrow \C
\]
is a complexification of the normalized Hodge-Riemann pairing.
Hodge-Riemann relations imply that $P_c(Comp)\subset \goth C$.

\ref{_from_PV_C_to_Pl_linear-alg_Proposition_}
establishes a diffeomorphism $\goth i:\; \goth C\arrow Pl$, where
$Pl$ is the space of 2-dimensional subspaces $L\subset H^2(M,\R)$
such that $(\cdot,\cdot)_{\c H}$ is positively defined on $L$.
Let $Q_c:\; Comp\arrow Pl$ be the composition of $P_c$ and $\goth i$.
As results of Bogomolov and Todorov imply (see
\cite{_Bogomolov_}, \cite{_Beauville_}, \cite{_Todorov_}),
the map $Q_c$ is an immersion. Since
$\dim Comp = \dim \goth C= \dim H^2(M) -2$, this map is
etale. It can be described in more straightforward terms
as follows. Let $I\in Comp$. Let $\tilde \Omega$ be a holomorphic
symplectic form over $(M, I)$. Let $\Omega\in H^2(M, \C)$ be the cohomology
class represented by the closed differential form $\tilde \Omega$.
Let $\omega_2:= Re(\Omega)$, $\omega_3:= Im(\Omega)$. Then one can
define $Q_c(I)$ as the linear span of $\omega_2$, $\omega_3$.


\section {Lefschetz-Frobenius algebras.}
\label{_Lefshe_Frob_Section_}

In this section, we give a number of preliminary definitions,
which eventually
lead to a calculation of the cohomology of a compact hyperkaehler
manifold. Some of these definitions are due to V. Lunts (see
\cite{_Lunts-Loo_}).

Further on, by ``algebra'' we understand
an associative algebra with unit.

\hfill

\definition 
Let $A$ be an algebra over a field $k$ and $(\cdot,\cdot)$
be a $k$-valued bilinear form on $A$. The form $(\cdot,\cdot)$
is called {\bf invariant} if for all $a,b,c\in A$,
$(ab,c)=(a,bc)$. A $k$-algebra equipped with an invariant
non-degenerate bilinear form is called {\bf Frobenius algebra}.

\hfill

\definition  
Let $A=\oplus A_i, i=0,...,d$ be a graded supercommutative algebra
over a field $k$ of characteristic zero. Assume that $A$ is equipped
with an invariant bilinear form, such that for all $a\in A_n, b\in A_m$,
the following holds:

\begin{equation} \label{_graded_scalar_pro_Equation_}
\begin{array}{l}
(a,b) = (-1)^{nm}(b,a), \mbox{\ and} \ \ \ \  \ \ \ \  \ \ \ \ \\[2mm]
(a,b)=0 \mbox{\ for\ } n+m\neq d. \ \ \ \  \ \ \ \  \ \ \ \ \\[2mm]
\end{array}
\end{equation}
Then $A$ is called {\bf a graded Frobenius algebra of degree $d$}.
The prototypical example of graded Frobenius algebras is the cohomology
algebra of a compact manifold.

\hfill

Further on, we consider only the graded Frobenius algebras.
For brevity, we sometimes omit the word ``graded''.

Let $A=\oplus A_i, i=0,...,d$ be a graded Frobenius algebra of
degree $d$ over the field $k$. Let $End_k(A)$ be the space of all
$k$-linear endomorphisms of $A$. Let $H\in End_k(A)$ be the endomorphism
which maps $a\in A_i$ to $(2d-i)\cdot a$. This endomorphism is
introduced by Hodge in his study of harmonic forms and Lefschetz
isomorphism. One can check that $H$ is a derivative of an algebra $A$:

\[ H(ab)= H(a) b + a H(b). \]

For all $A\in A$, let $L_a:\; A\arrow A$ map $b\in A$ to $ab$.

\hfill

\definition \label{_Lefschetz_triple_Definition_} 
Let $A=\oplus A_i$ be a graded Frobenius algebra, $a\in A_2$.
Let $L_a$, $H\in End_k(A)$ be as above. The triple of endomorphisms
$L_a$, $H$, $\Lambda_a\in End_k(A)$ is called {\bf a Lefshets
triple} if

\[ [H, L_a] = 2L_a, [H,\Lambda_a] = -2 \Lambda_a,
   [L_a,\Lambda_a] =H.
\]

Clearly, Lefschetz triples correspond to some representations
of the Lie algebra $\goth{sl}(2)$ in $End_k(A)$.
Lefschetz theorem (\cite{_Griffiths_Harris_}) gives examples
of Lefschetz triples for $A=H^*(M)$ and $M$ is a Kaehler manifold.

\hfill

\proposition \label{_Lefshe_tri_unique_Proposition_} 
Let $A$ be a graded Frobenius algebra, $a\in A_2$.
Let $(L_a, H, \Lambda_a)$, $(L_a, H, \Lambda'_a)$
be two Lefschetz triples. Then $\Lambda_a=\Lambda'_a$.
In other words, $\Lambda_a$ is uniquely determined by $a$.

{\bf Proof:} (V. Lunts) Consider the representations $\rho$, $\rho'$
of the Lie algebra $\goth{sl}(2)$ associated with these triples.
Take a basis $(x,y,h)$ in $\goth{sl}(2)$,

\[
   [h,x]=2x, [h,y] =-2y, [x,y]=h,
\]
such that $\rho(x) =\rho'(x) = L_a$, $\rho(h)= \rho'(h)= H$,
$\rho(y)=\Lambda_a$, $\rho'(y)=\Lambda'_a$.
Let $T:= \Lambda_a-\Lambda'_a$. Consider the adjoint
action of $\goth{sl}(2)$ on the space $End(A)$
obtained from $\rho$:

\[
   ad (\rho):\; \goth{sl}(2)\arrow End(End(A))
\]
Clearly, $ad (\rho)(x)(T)=[L_a,\Lambda_a]-[L_a,\Lambda'_a]=0$.
Therefore, $T$ is a highest vector is an $\goth{sl}(2)$-submodule
of $End(A)$ generated by $T$ and $ad (\rho)$. On the other hand,
$ad (\rho)(h)(T)=-2T$, and therefore, the weight of $T$ is $-2$.
This is impossible because $ad (\rho)$ is a finite dimensional
representation of $\goth {sl}(2)$. $\;\;\blacksquare$

\hfill

\definition 
Let $A=\oplus A_i$ be a graded Frobenius algebra, $a\in A_2$.
Then $a$ is called {\bf of Lefschetz type} if a Lefschetz triple
$(L_a, H, \Lambda_a)$ exists.

\hfill

\lemma \label{_a_Lefshe_if_a^i_iso_Lemma_} 
Let $A=A_0\oplus A_1\oplus ... A_{2d}$ be a
graded space. Let $L\in End(A)$ be an endomorphism of
grading 2: $L:\; A_i\arrow A_{i+2}$. Let $H$ act on $A_i$ as the
multiplication by $d-i$, $i=0,1, ... , 2d$.
Then the following conditions are equivalent:

\hfill

(i) There exist an endomoprhism $\Lambda\in End(A)$ of grading
-2 such that the relations

\[ [H,\Lambda] =-2 \Lambda, [H,L] =-2 L,
   [L,\Lambda] = H
\]
hold\footnote{The first two of these relations hold trivially
because of the grading.}.

\hfill

(ii) For all $i=0,1,...,d-1$, the map

\[ L^{d-i}:\; A^i\arrow A^2d-i,
\]
is an isomorphism.

\hfill

{\bf Proof:} Clear (see \cite{_Lunts-Loo_} for details).
$\;\;\blacksquare$

\hfill

When $A$ is a cohomology algebra of a compact
Kaehler manifold $M$, all Kaehler classes
are obviously of Lefschetz type. On the other hand, a class
of Lefschetz type is not necessarily a Kaehler class. For example,
for a Kaehler class $\omega$, the class $-\omega$ is of Lefschetz
type, but $-\omega$ cannot be a Kaehler class by trivial reasons.

\hfill

\definition \label{_Lefschetz_Frob_alge_Definition_} 
Let $A$ be a graded Frobenius algebra. Let $S\subset A_2$
be the set of all elements of Lefschetz type.  The algebra
$A$ is called {\bf a Lefschetz-Frobenius algebra} if the following
conditions hold:

(i) The space $A_0$ is one-dimensional over $k$.

(ii) The set $S$ is Zariski dense in $A_2$.

\hfill

{\bf Example}:
Let $M$ be a compact Kaehler manifold. Then the algebras
$H^*(M)$ and $\oplus H^{p,p}(M)$ are Lefschetz-Frobenius, as
\ref{_Lefshe_Frob_if_a_Lefshe_ele_exists_Proposition_}
implies.

\hfill

By \ref{_a_Lefshe_if_a^i_iso_Lemma_},
the set $S$ of all elements of Lefschetz type
is given by an open condition. Therefore $S$ is
open in $A_2$. Therefore, $A$ is
a Frobenius-Lefschetz algebra if and only if
$S$ is non-empty in $A_2$. We obtained the following
statement:

\hfill

\proposition \label{_Lefshe_Frob_if_a_Lefshe_ele_exists_Proposition_} 
Let $A$ be a graded Frobenius algebra. Assume that $A_2$ contains at
least one element of Lefschetz type. Then $A$ is Lefschetz-Frobenius.

$\blacksquare$

\hfill

For every Lefschetz triple $T$, we define a Lie algebra homomorphism

\[
    \rho_T:\; \goth{sl}(2) \arrow End_k(A)
\]
in an obvious way. For a Lefschetz-Frobenius algebra $A$, let $\g(A)$
be the Lie subalgebra of $End_k(A)$ generated by the images
of $\rho_T$ for all Lefschetz triples $T$. This algebra is our main object
of study. The algebra $\g(A)$ is graded: $\g(A)= \oplus\g_{2i}(A)$,
$g(A_n)\subset A_{n+2i}$ for all $g\in \g_{2i}(A)$. This
algebra is called {\bf the structure Lie algebra of $A$}.

\hfill

\definition \label{_Lefshe_Fro_Definition_} 
Let $A$ be a Lefschetz-Frobenius algebra. Assume that
$\g_{2i}(A)=0$ for $i\neq -1,0,1$:

\[
   \g(A)=\g_{-2}(A)\oplus \g_0(A)\oplus \g_2(A).
\]
Then $A$ is called {\bf a Lefschetz-Frobenius algebra of Jordan type}.
Such $A$ are closely related with Jordan algebras (\cite{_Springer_}).

\hfill

If $A$ is generated by $A_2$ and
$A_0\cong k$, $A$ is called {\bf reduced}. The subalgebra $A^r\subset A$
generated by $A_2$ and $A_0$ is called {\bf reduction of $A$}.
We use the following result.

\hfill

\proposition \label{_Lunts_about_FLJ_Proposition_} 
(\cite{_Lunts-Loo_})
Let $A$ be a Frobenius-Lefschetz algebra, $\g=\oplus \g_{2i}$ be its
structure Lie algebra. Then the following conditions are equivalent:

(i) $\g_2$ is spanned by $L_a$ for all Lefschetz elements $a$,

(ii) $A$ is of Jordan type,

(iii) $[\Lambda_a,\Lambda_b]=0$ for all Lefschetz elements $a,b\in A_2$.

{\bf Proof:} See Proposition 2.6 and Claim 2.6.1 of \cite{_Lunts-Loo_}.
$\;\;\blacksquare$

\hfill

\ref{_Lunts_about_FLJ_Proposition_} immediately implies the following
statement:

\hfill

\corollary \label{_g_2_is_A_2_Corollary_} 
Let $A=\oplus A_i$ be a Lefschetz-Frobenius algebra of Jordan type,
$\g=\g_{-2}\oplus \g_0\oplus \g_2$ be its sturtcure Lie algebra.
Then $\g_2=A_2$.

$\blacksquare$

\hfill

In the future, we often assume that $A$ is reduced.
In this case, the multiplicative structure on $A$ can be
recovered from the $\g(A)$-action. This is done as follows.
Since $[\g_2(A),\g_2(A)]\subset \g_4(A) =0$, the space
$\g_2(A)\subset \g(A)$ is a commutative
subalgebra of $\g(A)$. Consider the corresponding embedding of enveloping
algebras:

\[
   U_{\g_2(A)}\cong S^*(\g_2(A))\hookrightarrow U_{\g(A)}.
\]
Let ${\Bbb I}\in A_0$ be the unit. The representation $\g(A)\arrow End(A)$
induces the canonical map

\[
   U_{\g(A)} \stackrel{\tilde p}\arrow A, \; P\arrow P({\Bbb I}),
\]
where $P\in  U_{\g(A)}$ is an ``polynomial''
over $\g(A)$. Consider the restriction of $\tilde p$ to
$U_{\g_2(A)}\subset  U_{\g(A)}$:

\begin{equation}\label{_p_from_U_g_2_to_A_Equation_}
    p:\;  U_{\g_2(A)} \arrow A.
\end{equation}
Clearly, for all $a\in A_2$, $a$ of Lefschetz type, $L_a\in \g_2(A)$.
Since the set of elements of Lefschetz type is Zariski dense in $A_2$,
we have $L_a\in \g(A)\subset End(A)$ for all $a\in A_2$.
One can easily check that the corresponding map $i:\; A_2\arrow \g_2(A)$
is an isomorphism (see \ref{_Lunts_about_FLJ_Proposition_}).
Therefore, $S^*(\g_2(A))\cong S^*(A_2)$. Applying the isomorphism
$U_{\g_2(A)}\cong S^*(\g_2(A))\cong S^*(A_2)$ to the map
\eqref{_p_from_U_g_2_to_A_Equation_}, we obtain the map

\[ p':\; S^* A_2\arrow A. \]

\hfill

\claim \label{_p_is_multiplication_Claim_} 
The map $p'$ coinsides with the map induced by multiplication
in $A$.

{\bf Proof:} Clear. $\;\;\blacksquare$

\hfill

\ref{_p_is_multiplication_Claim_} implies that the kernel of the map
$p:\; U_{\g_2(A)} \arrow A$ is an ideal in $U_{\g_2(A)}$. This leads
to a more general construction.

\hfill

\definition \label{_multi_associ_with_representa_Definition_} 
Let $\g=\g_{-2}\oplus \g_0\oplus \g_2$ be a graded Lie algebra.
Let $V$ be a representation of $\g$ and $v\in V$ be a vector.
Assume that by applying $\g_2$ to $v$ repeatedly we obtain
the whole space $V$ (i. e., the vector $v$ generates $V$ as a
representation of $\g_2$). Assume that $\g_{-2}(v)=0$,
and that for all $g\in \g_0$, $g(v)$ is proportional to $v$. Let
$p:\; U_{\g_2}\arrow V$ be the map which associates with the polynomial
$P\in U_{\g_2}$ the vector $P(v)\in V$. Clearly, $\ker(p)$ is a left
ideal in $U_{\g_2}$. Since $\g_2$ is commutative, this ideal is two-sided.
Therefore, $V\cong \bigg(U_{\g_2}/ \ker (p)\bigg)$ is equipped with a structure
of commutative algebra. We denote this algebra by $V_{\g, v}$.

\hfill

\claim \label{_g_structure_defines_algebr_Claim_}
Let $A$ be a Lefschetz-Frobenius algebra of Jordan type. Assume that
$A$ is reduced (generated by $A_2$). Consider $A$ as a representation of
$\g(A)\cong \g_{-2}(A)\oplus \g_0(A)\oplus \g_2(A)$. Take the
unity vector ${\Bbb I}\in A_0$. Then the algebra $A_{\g(A),{\Bbb I}}$
coinsides with $A$.

{\bf Proof:} Follows from definitions. $\;\;\blacksquare$

\hfill

\hfill

{\large\bf Appendix. Reduction and the structure Lie algebra.}

\hfill

Let $A$ be a Lefschetz-Frobenius algebra, and $A^r$ be its reduction.
Assume that the restriction $(\cdot,\cdot)_r$ of $(\cdot,\cdot)$ to $A^r$
is non-degenerate. Then $(\cdot,\cdot)_r$
establishes a structure of Frobenius algebra on $A^r$.
We are going to show that $A^r$ is Lefschetz-Frobenius, and relate $\g(A)$
to $\g(A^r)$.

\hfill

\proposition \label {_A^r_Lef-Frob,_pres_by_g(A)_Proposition_} 
Let $A$ be a Lefschetz-Frobenius algebra.
Assume that the restriction of $(\cdot,\cdot)$ to $A^r$
is non-degenerate. Then $A^r$ is also Lefschetz-Frobenius.
Moreover, the action of $\g(A)$ on $A$ preserves the subspace
$A^r\subset  A$.

{\bf Proof:} Let $A^r_\bot$ be the orthogonal complement to $A^r$
in $A$. Since $(\cdot,\cdot)\restrict{A^r}$ is nondegenerate,
$A=A^r\oplus A^r_\bot$.

\hfill

\lemma \label{_A^r^bot_preserved_by_mult_by_A^r_Lemma_} 
Let $a\in A^r$, $b\in A^r_\bot$. Then $ab\in A^r_\bot$.

{\bf Proof:}  It is sufficient to show that for all $c\in A^r$,
$(ab,c)=0$. Since $(\cdot,\cdot)$ is invariant, for all $c\in A$ we have
$(ab,c)=(b, ac)$. Since $ac\in A^r$ and $b\in A^r_\bot$,
$(ab,c)=0$. $\;\;\blacksquare$

\hfill

{}From \ref{_A^r^bot_preserved_by_mult_by_A^r_Lemma_}, we obtain that
the operators $L_a$, $a\in A_2$ preserve the decomposition
$A=A^r\oplus A^r_\bot$. By \ref{_a_Lefshe_if_a^i_iso_Lemma_},
the map $L_a^{d-i}:\; A_i\arrow A_{2d-i}$ is an isomorphism.
Therefore, the restriction of $L_a^{d-i}$ to the $i$-th grading component
$(A^r)_i$ of $A^r$ is an embedding to $(A^r)_{2d-i}$. Since $A^r$
is Frobenius, $dim_k(A^r)_{2d-i}=dim_k(A^r)_i$. Therefore,
the restriction of $L_a^{d-i}$ to $(A^r)_i$
is an isomporphism. Applying \ref{_a_Lefshe_if_a^i_iso_Lemma_}
again, we obtain that for all Lefschetz-type elements
$a\in A_2$, these elements are of Lefschetz type
in $A^r$. Therefore, $A^r$ is a Lefschetz-Frobenius algebra.
It remains to prove that $A^r$ is preserved by $\g(A)$.
Clearly, the generators $L_a$ and $H$ of $\g(A)$ preserve
$A^r$. Therefore, to show that $\g(A)$ preserves $A^r$
it is sufficient to prove that $\Lambda_a$ preserve $A^r$
for all Lefschetz-type elements $a\in A_2$. Let $L=L_a$.
Let $L_r$, $H_r$ be the restrictions of $L$, $H$ to $A^r$ and
$L_{\bot}$ be the restrictions of $L$, $H$ to $A^r_\bot$.
By \ref{_a_Lefshe_if_a^i_iso_Lemma_}, there exists an
endomorphism $\Lambda_\bot:\; A^r_\bot\arrow A^r_\bot$
of grading $-2$ such that $[L_\bot,\Lambda_\bot]= H_\bot$.
Let $\Lambda_r:\; A^r\arrow A^r$ be the endomorphism
of $A^r$ such that $(L_r, H_r, \Lambda_r)$
is a Lefschetz triple. Let $\Lambda_r+\Lambda_\bot$ be an
endomorphism of $A$ such that for all $a=b+c$, $b\in A^r, c\in A_\bot^r$,

\[
   \Lambda_r+\Lambda_\bot(a)= \Lambda_r(b)+\Lambda_\bot(c).
\]
Checking relations, we obtain that $(L, H, \Lambda_r+\Lambda_\bot)$ is
a Lefschetz triple. By \ref{_Lefshe_tri_unique_Proposition_},
$\Lambda_a=\Lambda_r+\Lambda_\bot$. This implies that
$\Lambda_a$ preserves $A^r$. $\;\;\blacksquare$

\hfill

\ref{_A^r_Lef-Frob,_pres_by_g(A)_Proposition_} immediately
implies the following useful statement:

\hfill

\corollary 
Let $A$ be a Lefschetz-Frobenius algebra such that
its reduction $A^r$ is also Frobenius. Then $A^r$
is Lefschetz-Frobenius, and there exists a natural
Lie algebra epimorphism $\g(A)\arrow \g(A^r)$.

$\blacksquare$


\section {The minimal Frobenius algebras and cohomology of compact
Kaehler surfaces.} \label{_minimal_Fro_Section_}


In this section we concentrate on the simplest case of
Frobenius algebras related to Lefschetz theory. Namely,
we analyze the graded Frobenius algebras $A= A_0\oplus A_2\oplus A_4$,
where $dim_k A_0=1$. Such algebras are called minimal.
These algebras are naturally related to the complex surfaces.

\hfill

\definition 
Let $A= \oplus A_i$ be a graded Frobenius algebra.
Assume that $A= A_0\oplus A_2\oplus A_4$, $dim_k A_0=dim_k A_4=1$.
Then $A$ is called {\bf a minimal graded Frobenius algebra}.

\hfill

\proposition \label{_minimal_is_Lefschetz_Proposition_} 
Let $A= A_0\oplus A_2\oplus A_4$ be a minimal graded
Frobenius algebra. Then $A$ is Lefschetz-Frobenius.

{\bf Proof:}
Let $(\cdot,\cdot):\; A\times A\arrow k$ denote the
invariant scalar product on $A$.  The restriction
of $(\cdot,\cdot)$ to $A_2$ is a non-degenerate
bilinear symmetric form (it is non-degenerate
because of grading conditions \eqref{_graded_scalar_pro_Equation_}).
The following statement immediately implies
\ref{_minimal_is_Lefschetz_Proposition_}:

\hfill

\lemma \label{_el-t_with_non_zero_square_Lefschetz_Lemma_} 
Let $A= A_0\oplus A_2\oplus A_4$ be a minimal graded
Frobenius algebra. Let $a\in A_2$ be a vector such that
$(a,a)\neq 0$. Then $a$ is a Lefschetz element.

{\bf Proof:} Let $a^\bot$ be the orthogonal complement of $a$ in $A_2$:

\[
   a^\bot:= \{b\in A_2\; |\; (a,b)=0 \}.
\]
Let ${\Bbb I}\in A_0$ be the unit in $A$.
For all $b\in a^\bot$, $(ab,{\Bbb I})=(a,b)=0$.
Since $A_4$ is one-dimensional and its generator
has non-zero scalar product with ${\Bbb I}$, we have

\begin{equation} \label{_ab=0_for_all_b_in_a^bot_Equation_}
\forall b\in a^\bot,\; \;\; ab=0.
\end{equation}
Let $k_a\subset A_2$ be the one-dimensional space generated
by $a$. Let $A_a:= A_0\oplus ka \oplus A_4$.
Clearly, $A_a$ is a subalgebra of $A$. By
\eqref{_ab=0_for_all_b_in_a^bot_Equation_}, the operator
$L_a$ vanishes on $a^\bot$. Since $H(A_2) =0$, the operator
$H$ also vanishes on $a^\bot\subset A_2$. Therefore
it is sufficient to show that $a$ is a Lefschetz element in the algebra
$A_a$. Since $A_a\cong k[x]/(x^3=0)\cong H^*(\Bbb P^2, k)$, this follows
from Lefschetz theory. \ref{_el-t_with_non_zero_square_Lefschetz_Lemma_}
and consequently, \ref{_minimal_is_Lefschetz_Proposition_},
is proven. We also obtained the following result:

\hfill

\corollary \label{_Lambda_vanish_Corollary_} 
Let $(L_a, H, \Lambda_a)$ be the Lefschetz triple on $A$, where
$A$ is a minimal graded Frobenius algebra. Then
$(L_a, H, \Lambda_a)$ all vanish on $a^\bot$.

$\blacksquare$

\hfill

There is an easy way to construct the minimal graded Frobenius
algebras using spaces with non-degenerate symmetric bilinear forms.
Namely, let $V$ be a linear space over $k$, equipped witha  bilinear form
$(\cdot,\cdot)_V$. Consider the linear space

\[
   A(V):= k{{\Bbb I}}\oplus V\oplus k\Omega,
\]
where $k{\Bbb I}$ and $k\Omega$ are one-dimensional spaces
generated, respectively, by ${\Bbb I}$ and $\Omega$. We introduce a graded
Frobenius algebra structure on $A(V)$ in the following way. The grading
of $V$ is 2, the grading of ${\Bbb I}$ is 0, the grading of $\Omega$ is 4.
The product on $A(V)$ is defined as follows:

\hfill

(i) ${\Bbb I}$ is a unit.

(ii) for $v_1, v_2\in A_2(V)\cong V$, $v_1 v_2= (v_1,v_2)_V \Omega$.

\hfill

It remains to establish the invariant bilinear symmetric form
$(\cdot,\cdot)$ on $A(V)$.

\hfill

(iii) On $A_2(V)\cong V$, $(\cdot,\cdot)$
is equal to $(\cdot,\cdot)_V$.

(iv) The product of ${\Bbb I}$ and $\Omega$ is 1.

\hfill

Together with \eqref{_graded_scalar_pro_Equation_}, relations
(iii) and (iv) define the form $(\cdot,\cdot)$ is a unique way.

One can trivially check that this construction results in a
graded Frobenius algebra. Exactly this algebra appears as the
even cohomology of the compact Kaehler surface $M$, where
$V=H^2(M)$ and $(\cdot,\cdot)_V$ is the intersection form.
In fact, every minimal graded Frobenius algebra can be obtained
this way (\ref{_all_minim_alge_are_associ_Claim_}).

\hfill

\definition \label{_associ_gra_Fro_Definition_}
The graded Frobenius algebra $A(V)$ is called {\bf the minimal graded
Frobenius algebra associated with $V$, $(\cdot,\cdot)_V$}.

\hfill

\claim\label{_all_minim_alge_are_associ_Claim_} 
Let $A= A_0\oplus A_2\oplus A_4$ be the minimal graded
Frobenius algebra. Denote the restriction of
the invariant scalar product to $A_2$ by
$(\cdot,\cdot):\; A_2\times A_2\arrow k$.
Then $A$ is is isomorphic to the minimal graded
Frobenius algebra associated with $A_2$, $(\cdot,\cdot)$.

{\bf Proof:} Clear. $\;\;\blacksquare$

\hfill

\proposition\label{_minimal_alge_are_Jordan_Proposition_}
Let $A= A_0\oplus A_2\oplus A_4$ be a minimal graded Frobenius
algebra. Then $A$ is of Jordan type.

{\bf Proof:} By \ref{_Lunts_about_FLJ_Proposition_}, it is sufficient to
show that for every two elements $a_1, a_2\in A_2$ of Lefschetz
type, $[ \Lambda_{a_1},\Lambda_{a_2}]=0$.  Denote the generators of
$A_0$, $A_4$, by ${\Bbb I}$, $\Omega$, as in
\ref{_associ_gra_Fro_Definition_}. The endomorphism
$[ \Lambda_{a_1},\Lambda_{a_2}]$ has a grading $-4$. Therefore it is
a map from $A_4$ to $A_0$. To prove that
$[ \Lambda_{a_1},\Lambda_{a_2}]=0$ it is sufficient to show
that

\[
   [ \Lambda_{a_1},\Lambda_{a_2}]\:\bigg(\Omega\bigg)=0
\]
\ref{_Lambda_vanish_Corollary_} implies that $\Lambda_{a_i}(\Omega)$
is proportional to $a_i$. An easy calculation in $\goth{sl}(2)$
implies that $\Lambda_{a_i}(\Omega)=-a_i$. Similarly,
$\Lambda_{a_i}(a_j)=(a_i,a_j) {\Bbb I}$. Therefore

\[  [ \Lambda_{a_1},\Lambda_{a_2}]\:\bigg(\Omega\bigg) = (a_1,a_2)\cdot
    {\Bbb I} - (a_2,a_1) \cdot{\Bbb I} =0
\]
$\blacksquare$

\hfill

We proceed to compute the Lie algebra $\g(A)$ associated with
the minimal Frobenius algebra $A=A(V)$, where $V$ is a linear
space equipped with a scalar product. Denote by $\goth{so}(V)$
the Lie algebra of skew-symmetric endomorphisms of $V$. Let
$\goth H$ be the 2-dimensional space over $k$ with the hyperbolic
scalar product. In other words, $\goth H$ has a basis $x, y$
such that $(x,y)=1$, $(x,x)=0$, $(y,y)=0$.

By $\goth{so}(V)\oplus k$ we understand a direct sum of $\goth{so}(V)$
and a trivial Lie algebra of dimension 1.

\hfill

\theorem\label{_calculation_of_g(A)_for_minim_Theorem_} 
Let $V$ be a $k$-linear space equipped with a non-degenerate
scalar product. Let $A=A_0\oplus A_2\oplus A_4$
be the minimal graded Frobenius algebra $A(V)$
constructed by $V$. Take the graded Lie algebra
$\g(A)=\g_{-2}(A)\oplus \g_0(A)\oplus \g_2(A)$
(\ref{_Lefshe_Fro_Definition_}).
Then

\hfill

(i) $\g_0(A)\cong \goth{so}(V)\oplus k$,

(ii) $\g_2(A)\cong \g_{-2}(A)\cong V$,

(ii) $\g(A)\cong \goth{so}(V\oplus \goth H)$.

\hfill

{\bf Proof:} Denote the invariant bilinear form on $A(V)$ by
$(\cdot,\cdot)$. Let $(\cdot,\cdot)'$ another bilinear symmetric form,
defined by

\hfill

$(a,b)'=(a,b)$ if $a,b\in A_2$,

$(a,b)'=-(a,b)$ if $a,b\in A_0\oplus A_4$,

$(a,b)'=0$ if $a\in A_2$ and $b\in A_0\oplus A_4$.

\hfill

Let $A'$ be $A$ equipped with the scalar product $(\cdot,\cdot)'$.
Obviously,

\[ A'\cong V\oplus \goth H. \]

{\bf Step 1:} We are going to show
that $\g(A)\subset \goth{so}(A')$.

By trivial reasons,
$L_a$ and $H$ belong to $\goth{so}(A')$ for all $a\in A_2$.
Let $a\in A_2$ be an element of the Lefschetz type.
To prove that $\Lambda_a\in \goth{so}(A')$, we consider
the decomposition $A' = a^\bot \oplus A_a$ (see the proof
of \ref{_el-t_with_non_zero_square_Lefschetz_Lemma_}).
By \ref{_Lambda_vanish_Corollary_}, $\Lambda_a$ acts trivially
on $a^\bot$. Since the decomposition $A' = a^\bot \oplus A_a$
is orthogonal, it is sufficient to prove that the restriction
of $\Lambda_a$ to $A_a$ is skew-symmetric.

Three-dimensional representation $\rho:\; \goth{sl}(2)\arrow A_a$
obtained from the Lefschetz triple
$(L_a, \Lambda_a, H)$ is naturally isomorphic to the
adjoint representation of $\goth{sl}(2)$. Using this isomorphism,
we obtain that $(\cdot,\cdot)$ is the Killing form of the
Lie algebra $\g(A_a)\cong \goth{sl}(2)$. Therefore,
$\g(A_a)\subset End(A_a)$ consists of skew-symmetric matrices.
This finishes Step 1.

\hfill

{\bf Step 2:} We prove \ref{_calculation_of_g(A)_for_minim_Theorem_} (i).
Since $\g_0(A)$ is the grade-preserving part of $\g(A)$, we have a
homomorphism

\[
   \g_0(A)\stackrel {\mu}\arrow \goth{so}(A_2) \cong \goth{so}(V)
\]
which maps an endomorphism $h\in End(A)$ to its restriction
$h\restrict{A_2}\in End(A_2)$. The kernel of this homomorphism
is the space of all $h\in g_0(A)$ such that $h$ vanishes on $A_2$.
Therefore, $\ker\mu\in \goth{so}(A_0\oplus A_4)$. The algebra
$\goth{so}(A_0\oplus A_4)\cong \goth{so}(\goth H)\cong \goth{so}(1,1)$
is commutative and one-dimensional. Combining $\mu$ and the embedding

\[ \ker\mu\stackrel i\hookrightarrow
   \goth{so}(A_0\oplus A_4)\cong k,
\]
we obtain an embedding

\[
   \g_0(A)\stackrel m\hookrightarrow  \goth{so}(V) \oplus k.
\]
It remains to show that $m$ is a surjection. Consider the Hodge
endomorphism $H\in \g_0(A)\subset End(A)$ introduced a few sentences
before \ref{_Lefschetz_triple_Definition_}. By obvious reasons,
$H\in \ker\mu$.
The map $i:\; \ker \mu \arrow \goth{so}(A_0\oplus A_4)$ is surjective
because $i(H)$ is non-zero, and $\goth{so}(A_0\oplus A_4)$
is one-dimensional. Therefore, to prove that $m$ is surjective,
it is sufficient to show that $\mu:\;\g_0(A)\arrow \goth{so}(A_2)$
is surjective.

Let

\begin{equation} \label{_condi_on_a,b_Equation_}
   a,b\in A_2,\;\; (a,a)\neq, \;(b,b)\neq 0,\;(a,b)=0.
\end{equation}
Let $\inangles{a,b}\subset V$ be the plane generated
by $a$ and $b$, and $\inangles{a,b}^\bot\subset V$ be its orthogonal
completion.

Let $T_{ab}\subset \goth{so}(A_2)$ be the set of all
skew-symmetric endomorphisms which vanish on $\inangles{a,b}^\bot$.

Since $\goth{so}(2)$ is one-dimensional,
for given $a,b\in V$, all elements of $T_{ab}$
are proportional. We notice that the union of all
$T_{ab}$ generates the lie algebra $\goth{so}(A_2)$.
Therefore, to prove that $\mu$ is surjective
it is sufficient to show that $T_{ab}\subset \mu(\g_0(A))$ for all
$a,b$ satisfying conditions \eqref{_condi_on_a,b_Equation_}.
\ref{_Lambda_vanish_Corollary_} implies that
$[L_a,\Lambda_b]\in T_{ab}$. It is easy to check that
this element is non-zero. The proof of
\ref{_calculation_of_g(A)_for_minim_Theorem_} (i) is finished.

\ref{_calculation_of_g(A)_for_minim_Theorem_} (ii) immediately
follows from \ref{_Lunts_about_FLJ_Proposition_}.
\ref{_calculation_of_g(A)_for_minim_Theorem_} (iii) follows fom
the inclusion $\g(A)\subset \goth{so}(A')$ and comparing
dimensions, where dimension of $\g(A)$ is computed via
\ref{_calculation_of_g(A)_for_minim_Theorem_} (i) and
\ref{_calculation_of_g(A)_for_minim_Theorem_} (ii).
$\;\;\blacksquare$

\hfill

We denote the graded Lie algebra $\g(A)$ constructed by $V$ as in
\ref{_calculation_of_g(A)_for_minim_Theorem_} by $\goth{so}(V,+)$.
Clearly, over $\R$, when the symmetric form on $V$ has a signature
$(a,b)$%
\footnote{Of course, this means that
$\goth{so}(V)=\goth{so}(a,b)$}%
,

\[ \goth{so}(V,+)\cong \goth{so}(a+1,b+1). \]


\section[Representations of $SO(V,+)$
leading to Frobenius algebras.]
{Representations of $SO(V,+)$
leading to \\ Frobenius algebras.}
\label{_^dA(V)_Section_}


In this section, we describe all reduced Lefschetz-Frobenius
algebras $A= A_0\oplus A_2 \oplus ... \oplus A_{2d}$ with
$\g(A)\cong \goth{so}(V,+)$. It turns out that such algebras
are uniquely defined by the number $d$, which is {\it even},
whenever $dim V>2$.

Let $V$ be a linear space supplied with a non-degenerate symmetric
bilinear form $(\cdot,\cdot)$. Let $A=A(V)$ be the minimal
Frobenius algebra constructed in Section \ref{_minimal_Fro_Section_}.
Take the tensor product

\[ A^{\otimes d}:=
   \underbrace{A\otimes_k A\otimes_k ... A}_{d\mbox{\ \ times}}.
\]
There is a natural action of $\goth{so}(V,+)$ on $A^{\otimes d}$.
Let ${}^{(d)}A$ be the irreducible $\goth{so}(V,+)$-module generated
by ${}^d{\Bbb I}:= {\Bbb I} \otimes {\Bbb I}\otimes ... {\Bbb I}$,
where ${\Bbb I}\in A$ is the unit. The space ${}^{(d)}A$ is not necessarily a
subalgebra in $A^{\otimes d}$. We introduce a new
algebra structure on ${}^{(d)}A$ which does
not necessarily come from the algebra structure on
$A^{\otimes d}$, but instead comes from $\goth{so}(V,+)$-action
as in \ref{_multi_associ_with_representa_Definition_}.
Denote the graded Lie algebra $\goth{so}(V,+)$ by
$\g= \g_{-2}\oplus \g_0\oplus \g_2$. The algebra
$\g_0$ has an additional decomposition:
$\g_0= \goth{so}(V)\oplus kH$
(\ref{_calculation_of_g(A)_for_minim_Theorem_}).
Clearly, all elements of $\g_{-2}$ vanish on ${}^d{\Bbb I}$,
$\goth{so}(V)\subset \g_0$ vanish on ${}^d{\Bbb I}$ and
$H$ acts on ${}^d{\Bbb I}\in {}^{(d)}A$ as multiplication by $-2d$.
Therefore we may apply
\ref{_multi_associ_with_representa_Definition_}
to the representation ${}^{(d)}A$ and the vector ${}^d{\Bbb I}$.
Denote the resulting algebra ${}^{(d)}A_{\g, {\Bbb I}}$
by ${}^dA(V)$.

\hfill

\theorem \label{_all_alg_with_so_are_^dA_Theorem_} 
Let $\tilde A=A_0\oplus A_2\oplus ... \oplus A_{2n}$ be a
Lefschetz-Frobenius algebra of Jordan type. Assume that the graded
Lie algebra $\g(A)$ is isomorphic to $\goth{so}(V,+)$ and $dim_k V> 2$.
Let $A$ be the subalgebra of $\tilde A$ generated by $A_0$ and $A_2$
(also known as {\bf reduction} of $\tilde A$).
Then

(i) $n$ is even

(ii) $A$ is isomorphic to ${}^{n/2} A(V)$ as a graded algebra%
\footnote{The invariant Frobenius pairing is unique up to a scalar,
which is easy to see.}%
{}.

\hfill

{\bf Remark:} In particlural, this theorem implies that the reduction
$A$ of $\tilde A$ is Frobenius, which is not immediately clear
in general case.

\hfill

{\bf Proof:} Let $\g=\g_{-2}\oplus \g_0\oplus \g_2$ be the graded
Lie algebra $\goth{so}(V,+)$. Denote the unit in $A$ by ${\Bbb I}$.
By definition, $A_0$ is $\g_0$-invariant one-dimensional space.
We have shown that
$\g_0= \goth{so}(V)\oplus kH$. Since $\goth{so}(V)$ is a simple
Lie algebra, $\goth{so}(V)({\Bbb I})=0$. Therefore, to prove
\ref{_all_alg_with_so_are_^dA_Theorem_} (i), it is sufficient
to prove the following lemma.

\hfill

\lemma \label{_repres_so(V,+)_even-weight_Lemma_} 
Let $\rho:\;\g\arrow End(M)$ be a simple representation
of $\g\cong \goth{so}(V,+)$, where $dim_k(V)>2$. Let
${\Bbb I}\in M$ be a vector such that $\g_{-2}({\Bbb I})=0$,
$\goth{so}(V)({\Bbb I})=0$,%
\footnote{As elsewhere, we use the decomposition
$\g_0\cong \goth{so}(V)\oplus kH$ provided by
\ref{_calculation_of_g(A)_for_minim_Theorem_}.}
and $H({\Bbb I})=-2n\cdot{\Bbb I}$. Then $n$ is even.

{\bf Proof:}

{\bf Step 1:} We show that $k\neq 1$.

Assume the contrary.
Consider the decomposition $M\cong M_{-1}\oplus M_1$, given
by the weights of $H$.  Since ${\Bbb I}$ is a highest
weight vector for some root system in $\g$
(see \ref{_Cartan_exists_in_so(V,+)_Claim_}),
the corresponding weight space is one-dimensional.
This implies that $dim_k(M_{-1})=1$. There exists an automorphism
of $\goth{so}(V,+)$ which maps $H$ to $-H$ (an easy check).
Therefore, $dim_k(M_{-1})=dim_k(M_1)=1$. In other words,
the representation $M$ is two-dimensional. Consider $M$
as a representation of $\goth{so}(V)\subset \goth{so}(V,+)$.
Since ${\Bbb I}$ is invariant with respect to $\goth{so}(V)$,
the space $M$ is decomposed into a sum of two
one-dimensional $\goth{so}(V)$-invariant subspaces.
For $dim_k(V)>2$, the Lie algebra $\goth{so}(V)$
is simple. Therefore its one-dimensional representations
are trivial. We obtain that $\goth{so}(V)(M)=0$.
Since $\g=\goth{so}(V,+)$ is a simple Lie algebra,
the homomorphism $\rho:\; \g\arrow End(M)$
cannot have proper non-zero kernel. Therefore,
$\rho(\g)=0$.

{\bf Step 2:} We use the following statement, which is easy to check.

\hfill

\claim \label{_Cartan_exists_in_so(V,+)_Claim_} 
Let $\goth{so}(V)\oplus kH \stackrel i \hookrightarrow  \goth{so}(V,+)$
be the embedding provided by \ref{_calculation_of_g(A)_for_minim_Theorem_}
(i). Then there exists a Cartan subalgebra $\goth h'\subset \goth{so}(V)$
such that $i(\goth h'\oplus kH)$ is a Cartan subalgebra in $\goth{so}(V,+)$.

$\blacksquare$

\hfill

Take a Cartan subalgebra $\goth h:= i(\goth h'\oplus kH)\subset \g$
provided by \ref{_Cartan_exists_in_so(V,+)_Claim_}. Clearly,
the linear map

\[ -H\galochka:= -(H,\cdot), \;\;-H\galochka:\; \goth h\arrow k\]
is a root. Taking a root system $\alpha_1, ... ,\alpha_m$ in
$\goth h'\subset \goth{so}(V)$, we obtain that
$\alpha_1, ... ,\alpha_m, -H\galochka$ is a root system in
$\g$. In this root system, ${\Bbb I}\in M$ is a highest weight
vector of the representation $M$. It is known that the
set\footnote{Here, $\goth h\galochka$ means $\goth h$ dual}
$W\subset \goth h\galochka$ of possible
weights of the highest weight vector coinsides with the intersection
of a weight lattice $L$ in $\goth h\galochka$ and a Weyl chamber.
In particular, $W\subset\goth h\galochka$ is an abelian semi-group
with group structure induced from $\goth h\galochka$.

The weight of ${\Bbb I}$ corresponding to the root system
$\alpha_1, ... ,\alpha_m, -H\galochka$ is $(0,0 ,..., 2n)$.
Let $W_0\subset W$ be the set of all
$(0,0 ,..., 2n)\subset W\subset \goth h\galochka$ which correspond to
representations $M$ satisfying conditions of
\ref{_repres_so(V,+)_even-weight_Lemma_}.
Clearly,

\[ W_0= \bigg \{ (x_1, .... , x_m)\in
   W \;\;|\;\; x_1=...=x_{m-1}=0\bigg \}
\]
By Step 1, $n\neq 1$. Since Weyl chamber is invariant with
respect to homotheties, the semigroup $W_0$ is isomorphic
to $\Z_{\geq 0}$. Therefore, $n$ is never odd.
\ref{_repres_so(V,+)_even-weight_Lemma_} is proven.
$\;\;\blacksquare$

\hfill

To prove \ref{_all_alg_with_so_are_^dA_Theorem_},
we notice that the simple representation of a reductive Lie algebra
is uniquely determined by its highest weight.
The weights of representations of $\goth{so}(V,+)$
corresponding to Lefschetz-Frobenius algebras are computed
a few lines above. In particular, we obtained that for
some root system in $\g$, the highest weight of $A$
is $(0,0,0... , 2n)$. Simple  representations with a given
highest weight are isomorphic. Therefore,
as a representation of $\g$, $A\cong {}^dA(V)$.
By \ref{_g_structure_defines_algebr_Claim_},
the action of $\g$ detrermines multiplication in $A$ uniquely.
This finishes the proof of \ref{_all_alg_with_so_are_^dA_Theorem_}.
$\;\;\blacksquare$


\section[Computing the structure Lie algebra for the
cohomology of a hyperkaehler manifold,
part \ I.]{Computing the structure Lie algebra for the
cohomology of a hyperkaehler manifold, \\part \ I.}
\label{_computing_g_for_hyperk_pt-I_Section_}


Let $M$ be a compact Kaehler manifold. According to
\ref{_Lefshe_Frob_if_a_Lefshe_ele_exists_Proposition_},
the ring $A:=H^*(M)$ is Lefschetz-Frobenius. The aim of this
section is to compute $\g(A)$ in the case when
$M$ is a simple compact hyperkaehler manifold. The answer is
hinted at by the following statement, which is proven in
\cite{_so5_on_cohomo_}:

\hfill

\proposition \label{_subalg_in_g_genera_by_three_Kae_forms_Proposition_}
Let $\c H\in Hyp$ be a hyperkaehler structure on $M$. Consider the
Kaehler classes

\[ \omega_I= P_1(\c H),\;\; \omega_J=P_2(\c H),\;\;
   \omega_K= P_3(\c H), \;\;\omega_I, \omega_J,\omega_K\in H^2(M,\R).
\]
Cohomology classes $\omega_I, \omega_J,\omega_K$ are Lefschetz elements in
the graded Frobenius algebra $A:= H^*(M,\R)$. Consider the graded
subalgebra $\g(\c H)$ in $\g(A)$ generated by
$L_{\omega_I},L_{\omega_J},L_{\omega_K}$,
$\Lambda_{\omega_I}, \Lambda_{\omega_J},\Lambda_{\omega_K}$ and
$H$. Then $\g(\c H)$ is isomorphic to $\g(\R^3)$, where
$\g(\R^3)\cong \goth{so}(4,1)$ is the structure Lie algebra
of the minimal graded Lefschetz-Frobenius algebra corresponding
to the linear space $\R^3$ with positively defined scalar product.
In particular, the graded algebra
$\g(\c H)$ is independent from $\c H$ and $M$.

{\bf Proof:} This statement is proven in \cite{_so5_on_cohomo_}.
It is based on the commutation relations in $\g(\c H)$ given as follows.
Denote $P_i(\c H)$ by $\omega_i$, $i=1,2,3$. Denote $L_{\omega_i}$ by
$L_i$, and $\Lambda_{\omega_i}$ by $\Lambda_i$. Let
$K_{ij}:= [L_i,\Lambda_j]$, $i\neq j$.
Then the following relations  are true:

\begin{equation} \label{_so5_relations_Equation_}
\begin{array}{l}
{}[ L_i, L_j] = [\Lambda_i,\Lambda_j] =0;\;\; \\[2mm]
[L_i,\Lambda_i]= H;\;\; [H, L_i] = 2 L_i;
\;\; [H, \Lambda_i]=-2\Lambda_i \\[2mm]{}
K_{ij}=-K_{ji}, [K_{ij}, K_{jk}]=2 K_{ik}, [K_{ij}, H] =0 \\[2mm]{}
[K_{ij} L_j]=2 L_i;\;\; [K_{ij} \Lambda_j] = 2 \Lambda_i \\[2mm]{}
[K_{ij}, L_k] = [K_{ij}, \Lambda_k] =0\;\; (k\neq i,j)\\[2mm]
\end{array}
\end{equation}

 $\;\;\blacksquare$

\hfill

Let $M$ be a simple hyperkaehler manifold, $A=H^*(M,\R)$ be its
cohomology ring, equipped with the invariant pairing provided
by the Poincare duality. We consider $A=\oplus A_i=\oplus H^i(M,\R)$
as a graded Frobenius algebra over $\R$.
In Section \ref{_Hodge-Rie_independent_Section_},
we defined the normalized Hodge-Riemann pairing $(\cdot,\cdot)_{\c H}$
on $H^2(M,\R)$. Let $V$ be a linear space
$H^2(M)$ equipped with this pairing.
In Section \ref{_minimal_Fro_Section_}, we defined a graded
Frobenius algebra $\g(V)$, also denoted by $\goth{so}(V,+)$.
By definition, over $\R$ the algebra $\goth{so}(V,+)$ is isomorphic to
$\goth{so}(m+1,n+1)$, where $(m,n)$ is the signature of $V$.

\hfill

\theorem\label{_g(A)_for_hyperkae_Theorem_} 
In this notation, $\g(A)$ is isomorphic to $\g(V)$.

\hfill

\ref{_g(A)_for_hyperkae_Theorem_} is the main result of this paper.
It is proven in the subsequent sections. The present section is
dedicated to proving that the Lefschetz-Frobenius
algebra $A$ is of Jordan type. This
is a crucial step in proof of \ref{_g(A)_for_hyperkae_Theorem_}.

\hfill

{\bf Remarks:} For a hyperkaehler
manifold with $dim \; H^2(M) =3$ (minimal dimension which
is not obviously impossible), \ref{_g(A)_for_hyperkae_Theorem_}
is equivalent to \ref{_subalg_in_g_genera_by_three_Kae_forms_Proposition_}.
For a hyperkaehler manifold with $dim\; H^2(M) =n$,
the Riemann-Hodge metric on $H^2(M,\R)$ has a signature $(3,n-3)$.
This means that $\g(V)$ is isomorphic to $\goth{so}(4,n-2)$.

\hfill

\proposition \label{_H^*_hyp_Jordan_type_Proposition_} 
Let $A$ be the Lefschetz-Frobenius
algebra of cohomology of a simple compact hyperkaehler
manifold. Then $A$ is of Jordan type.

{\bf Proof:} According to \ref{_Lunts_about_FLJ_Proposition_},
it is sufficient to show that for all $a,b \in A_2$,
$a$, $b$ of Lefschetz type, $[\Lambda_a,\Lambda_b]=0$.
Let $a,b\in A_2$, $\c H\in Hyp$, $\c H = (I,J,K, (\cdot,\cdot))$.
We write $a\bullet_{\c H}b$ if there exist complex structures
$I_1, I_2$ which are induced by $\c H$ such that $a, b$ are
Kaehler classes corresponding to $I_1$, $I_2$ and the metric
$(\cdot,\cdot)$.\footnote{the metric $(\cdot,\cdot)$ is Kaehler
with respect to the complex structures $I_1$, $I_2$, as the
definition of induced complex structures implies.}
Clearly, if $a\bullet_{\c H} b$, then $a$ and $b$
are of Lefschetz type.

\hfill

\lemma \label{_a_bullet_H_b=>_Lambdas_commute_Lemma_} 
Let $\alpha,\beta\in A_2$, $\c H\in Hyp$, $\alpha\bullet_{\c H}\beta$.
Then $[\Lambda_\alpha,\Lambda_\beta]=0$.

{\bf Proof:} Let $\g(\c H)$ be the graded Lie algebra defined
in \ref{_subalg_in_g_genera_by_three_Kae_forms_Proposition_}.
Since $\g(\c H)=\g(\R^3)$,
 the $-4$-th component of $\g(\c H)$ vanishes: $\g_{-4}(\c H)=0$.
Therefore for all $\lambda,\mu\in \g_{-2}(\c H)$, we have
$[\lambda,\mu]=0$. Therefore, to prove
\ref{_a_bullet_H_b=>_Lambdas_commute_Lemma_} it is sufficient
to show that $\Lambda_\alpha,\Lambda_\beta\in \g_{-2}(\c H)$.
This is implied by the following lemma.

\hfill

\lemma\label{_Kaehle_cla_indu_by_H_in_g(H)_Lemma_} 
Let $\c H\in Hyp$, $\c H=(I,J,K,(\cdot,\cdot))$
and $I'$ be a complex structure induced
by $\c H$. Let $\omega\in A_2$ be the Kaehler class corresponding
to $I'$ and the metric $(\cdot,\cdot)$. Then
$\Lambda_\omega \in \g_{-2}(\c H)$.

{\bf Proof:} Let $\omega_i:=P_i(\c H)$, $i=1,2,3$.
By definition of induced complex structures,
$I'= a I+b J+ c K$, where $a,b,c\in \R$, $a^2+b^2+c^2=1$.
Since $\omega(x,y)= (x, I'y)$, we have
$\omega=a\omega_1+b\omega_2+c\omega_3$. Let

\[ \Lambda:= a\Lambda_{\omega_1}+
   b\Lambda_{\omega_2} + c\Lambda_{\omega_3}.
\]
Using \eqref{_so5_relations_Equation_}, we see that
$[L_\omega,\Lambda]=H$. The other relations defining Lefschetz
triples checked automatically, we obtain that
$(L_\omega, H, \Lambda)$ is a Lefschetz triple. Therefore,
$\Lambda=\Lambda_\omega\in\g_2(\c H)$. This proves
\ref{_a_bullet_H_b=>_Lambdas_commute_Lemma_}
and \ref{_Kaehle_cla_indu_by_H_in_g(H)_Lemma_}.
$\;\;\blacksquare$

\hfill

Let $a,b\in A_2$. We write $a\bullet b$ if there exist $\c H\in Hyp$,
$\lambda\in \R$, $\lambda\neq 0$ such that $a\bullet_{\c H}\lambda b$.
Clearly, $\Lambda_{\lambda b} = \lambda^{-1}\Lambda_b$.
Therefore, \ref{_a_bullet_H_b=>_Lambdas_commute_Lemma_}
implies the following statement.

\hfill

\claim \label{_a_bullet_b_Lambdas_commute_Claim_} 
Let $a,b\in A_2$, $a\bullet b$. Then $[\Lambda_a,\Lambda_b]=0$.

$\blacksquare$

\hfill

The relation $a\bullet b$ is prominent further on in this paper.
We could have given an alternative definition of $a\bullet b$ as follows.

\hfill

\lemma \label{_a_bullet_b_if_hyperk_exists_Lemma_} 
Let $a,b \in A_2$. The following conditions are equivalent:

(i) $a \bullet b$

(ii) There exists $\c H\in Hyp$ such that $a$ and $b$
can be expressed as linear combinations of $P_i(\c H)$,
$i=1,2,3$.

{\bf Proof:} Clear from \ref{_action_SO(3)_on_Hyp_via_periods_Lemma_}.
$\;\;\blacksquare$

\hfill

{\bf Remark:} Let $a\bullet b$, $a,b\neq 0$. Then both $a$ and $b$ are
elements of Lefschetz type.

\hfill

Let $S\in A_2$ be a set of all elements of Lefschetz type.
By \ref{_a_Lefshe_if_a^i_iso_Lemma_}, $S$ is Zariski open in $A_2$.
Let $\nu:\; S\arrow End(A)$ map $a\in S$ to $\Lambda_a\in End(A)$.
An easy linear-algebraic check implies that $\nu$ is an algebraic
map. Therefore the map $\eta:\; S\times S\arrow End(A)$,
$a,b \arrow [\Lambda_a,\Lambda_b]$ is also an algebraic map.
To prove \ref{_H^*_hyp_Jordan_type_Proposition_}
it is sufficient to show that $\eta$ is identically zero,
as shown by \ref{_Lunts_about_FLJ_Proposition_}.
By \ref{_a_bullet_b_Lambdas_commute_Claim_},
for all $a,b\in A_2$, such that $a\bullet b$,
we have $\eta(a,b)=0$. Therefore,
\ref{_H^*_hyp_Jordan_type_Proposition_} is implied by
the following statement:

\hfill

\lemma \label{_a_bullet_b_dense_in_A_2_x_A_2_Lemma_} 
Let $D\subset A_2\times A_2$ be the set of all pairs
$(a,b)\in A_2\times A_2$ such that $a\bullet b$. Then
$D$ is Zariski dense in $A_2$.

{\bf Proof:} We use the following trivial
statement:

\hfill

\sublemma \label{_a_bullet_b_linear_span_Sublemma_} 
Let $(a,b)$ and $(a',b')$ be two pairs in $A_2$ such that
the linear span of $(a,b)$ is equal to the linear span
of $(a',b')$. Then

\[ a\bullet b\;\;\Leftrightarrow \;\; a'\bullet b'. \]

{\bf Proof:} See
\ref{_a_bullet_b_if_hyperk_exists_Lemma_}.
$\;\;\blacksquare$

\hfill

Denote the Grassmanian of all 2-dimensional planes in $H^2(M,\R)$ by $Gr$.
Let $Gr^\bullet$ be the space of planes generated by $a, b$ with
$a\bullet b$. To prove \ref{_a_bullet_b_dense_in_A_2_x_A_2_Lemma_}
it is sufficient to show that $Gr^\bullet$ is Zariski dense in $Gr$.

\hfill

In Section \ref{_Q_c_defini_Section_}, we defined the period map
$Q_c:\; Comp\arrow G_r$.

\hfill

\claim \label{_D_in_Q_c(Comp)_Claim_} 
The space $Gr^\bullet$ coincides with $Q_c(Comp)$.

{\bf Proof:} Clear from definitions. $\;\;\blacksquare$

\hfill

Since the map $Q_c$ is etale (see the end of Section
\ref{_Q_c_defini_Section_}), its image contains an open set
in $Gr$. Therefore $Q_c(Comp)$ is Zariski dense in $Gr$.
\ref{_a_bullet_b_dense_in_A_2_x_A_2_Lemma_}, and consequently,
\ref{_H^*_hyp_Jordan_type_Proposition_}, is proven.
$\;\;\blacksquare$


\section{Calculation of a zero graded part of the structure Lie
algebra of the cohomology of a hyperkaehler manifold, part I.}
\label{_compu_g_0_part_1_Section_}


In this section, we make steps related to computation of
the grade zero part of the Lie algebra
$\g(A)=\g_{-2}(A)\oplus \g_0(A)\oplus \g_2(A)$.
As in the previous section, $V$ denotes the space $H^2(M,\R)$,
considered as a linear space with the scalar product
$(\cdot, \cdot)_{\c H}$, and $A$ is the Frobenius algebra
$H^*(M,\R)$. We construct an epimorphism
$\g_0(A)\arrow \goth{so}(V)\oplus kH$, where $kH$ is one-dimensional
Lie algebra.

\hfill

With every $I\in Comp$, we associate the semisimple
endomorphism $ad I$ of $A$, defined as follows (see also
Section \ref{hyperk_manif_Section_}). Consider the Hodge decomposition
$H^i(M,\C)=\oplus_{p+q=i}H^{p,q}(M)$. Let
$ad^cI:\; H^i(M,\C) \arrow H^i(M,\C)$ multiply $(p,q)$-forms
by the scalar $(p-q)\1$. One can check that $ad^c I$ commutes
with the standard real structure on $H^i(M,\C)$. Therefore,
$ad^c I$ is a complexification of a (uniquely defined) endomorphism of
$H^i(M,\R)$. Denote this endomorphism by $ad I$. This definition
coinsides with one given in Section \ref{hyperk_manif_Section_}.

Consider the action of $ad I$ on $V=H^2(M, \R)$. Using
Hodge-Riemann relations, we immediately obtain that
$ad I\restrict{V}\in \goth{so}(V)$. Let $\c M_2\subset End(V)$
be the Lie algebra generated by endomorphisms $ad I\restrict{V}$,
for all $I\in Comp$. Clearly, $\c M_2\subset \goth{so}(V)$.
One can show that $\c M_2$ is a Mumford-Tate group of
$(M, I)$ for generic $I\in Comp$, although we never use this
observation.

Let $v:\; B\arrow B$ be an endomorphism of a linear space $B$.
We call the endomorphism

\[ v^{\circ{}}\in End(B),
   \;\;v^{\circ{}}:= v-\frac{1}{\dim B} Tr(B) Id_B
\]
{\bf the traceless part of $v$}. Clearly, the map
$Tl:\; \goth{gl}(V)\arrow \goth{sl}(V)$,
$Tl(v)= v^{\circ{}}$ is a Lie algebra homomorphism.
For all $g\in \g_0(A)$, consider the restriction
$g\restrict{A_2}:\; A_2\arrow A_2$. Let $g^\circ\in End(A_2)$ be the
traceless part of $g\restrict{A_2}$. This defines a Lie
algebra homomorphism\footnote{$V$ is $A_2$ is $H^2(M,\R)$}
$t:\; \g_0(A)\arrow \goth{sl}(V)$, $t(g)= g^\circ$.
Consider the one-dimensional Lie algebra in $End(A)$,
generated by the Hodge endomorphism $H\in End(A)$.
Denote this algebra by $k H$.
Let $s:\; \g_0(A)\arrow k H$ map $g\in \g_0(A)$ to

\[
   -\frac{Tr(g\restrict{A_0})}{\frac{1}{2}\dim_\R M} H\in kH.
\]
The map $s$ is defined in such a way that $s(H)= H$.
Clearly, $s$ is also a Lie algebra homomorphism.

\hfill

\proposition \label{_str_of_g_0_Proposition_} 
The following statements are true:

\hfill

(i) $t(\g_0(A))\subset \c M_2$

(ii) The inclusion $\c M_2\subset \goth{so}(V)$ is an equality:
$\c M_2= \goth{so}(V)$.

(iii) The map $t\oplus s:\; \g_0(A)\arrow \goth{so}(V)\oplus k H$
is an epimorphism.

\hfill

{\bf Proof:} We use the following simple lemma.

\hfill

\lemma \label{_g_0_gener_by_[L,Lambda]_Lemma_} 
(see also \cite{_Lunts-Loo_})
Let $\c A$ be a Lefschetz-Frobenius algebra of Jordan type,
$\g(\c A)=\g_{-2}(\c A)\oplus\g_0(\c A)\oplus\g_2(\c A)$ be its
structure Lie algebra. Then $\g_0(\c A)$ is a linear span
of the elements $[L_a,\Lambda_b]$, where $a,b$ are Lefschetz
elements in $\c A_2$.

{\bf Proof:} Clear.
$\;\;\blacksquare$

\hfill

Let $S\subset A_2$ be the set of all elements of Lefschetz type.
Let $\nu:\; S\times S \arrow \g_0(A)$ be the map
$a,b\arrow [L_a,\Lambda_b]$. \ref{_g_0_gener_by_[L,Lambda]_Lemma_}
implies that $\g_0(A)$ is a linear span of the set
$\nu(S\times S)$. Therefore, to prove \ref{_str_of_g_0_Proposition_}
(i) it is sufficient to show that for all $a,b\in S$,
$t(\nu(a,b))\in \c M_2$. As we have seen, $S$ is a
Zariski open subset in $A_2$. Consider $S$ as an algebraic
manifold with the algebraic structure induced from $A$.
With respect to this algebraic structure, both $t$ and
$\nu$ are algebraic maps. Therefore it is sufficient to
show that for a Zariski dense subset $D\subset S\times S$,
$t(\nu(D))\subset \c M_2$. By
\ref{_a_bullet_b_dense_in_A_2_x_A_2_Lemma_},
the set of all $a,b\in S$, $a\bullet b$ is
Zariski dense in $S$. Therefore, the inclusion
$t(\nu(S\times S))\subset \c M_2$ is implied by the following
statement:

\hfill

\lemma \label{_[L_a,Lambda_b]=I_for_(a,b)_in_Q_c(I)_Lemma_} 
Let $a,b \in S$, $a\bullet b$, where $a$ is not proportional to $b$.
Let $\c L$ be a plane in
$H^2(M,\R)$ spanned by $a$ and $b$. According to
\ref{_D_in_Q_c(Comp)_Claim_}, there exist $\c I\in Comp$ such
that $Q_c(I)= \c L$. Let

\[ \xi_1:=
   \bigg([L_a,\Lambda_b]\bigg|_{{}_{H^2(M,\R)}}\bigg)^\circ
   \in End(H^2(M, \R)
\]
be the traceless part of the restriction of the linear operator
$[L_a,\Lambda_b]\in End(A)$
to $H^2(M, \R)=A_2\subset A$. Let $\xi_2$ be the restriction
of $ad I$ to $H^2(M, \R)$. Then $\xi_1$ is proportional to $\xi_2$.

{\bf Proof:} The space $A_2\cong H^2(M, \R)$ is equipped with the
normalized Hodge-Riemann pairing $(\cdot,\cdot)_{\c H}$.
Let $b,x$ be an orthogonal basis of $\c L$. Clearly,
$x=\lambda a +\mu b$, where $\lambda,\mu\in \R$. This implies that

\begin{equation} \label{_adding_L_in_[L,Lambda]_Equation_}
   [L_x,\Lambda_b] = \lambda H +\mu [L_a,\Lambda_b].
\end{equation}
Since the endomorphism $\lambda H\restrict{A_2}\in End(A_2)$ is proportional
to identity, the traceless parts

\[ \bigg([L_a,\Lambda_b]\bigg)^\circ, \;\;
   \bigg([L_x,\Lambda_b]\bigg)^\circ
\]
are proportional. Therefore, proving
\ref{_[L_a,Lambda_b]=I_for_(a,b)_in_Q_c(I)_Lemma_}, we may
assume that

\[ (a,b)_{\c H}=0, \;\;
   (a,a)_{\c H}= (b,b)_{\c H}>0.
\]

Let $\c H\in Hyp$ be a hyperkaehler structure such that
$\Phi^{hyp}_c(\c H)= \c I$. The relation $\Phi^{hyp}_c(\c H)= \c I$
means that $\c H= (\c I, J, K, (\cdot,\cdot))$ for some $J, K, (\cdot,\cdot)$.
By definition of $Q_c$, the linear span of $P_2(\c H)$, $P_3(\c H)$
coinsides with $\c L$. Using
\ref{_action_SO(3)_on_Hyp_via_periods_Lemma_},
we can easily find a hyperkaehler structure $\c H'\in Hyp$,
$\c H'$ equivalent to $\c H$, such that $P_2(\c H)=\lambda a$,
$P_3(\c H)=\lambda b$ for some $\lambda \in \R$. Now,
\ref{_[L_a,Lambda_b]=I_for_(a,b)_in_Q_c(I)_Lemma_}
is a consequence of the following simple statement:

\hfill

\claim \label{_[L_2,Lambda_3]=adI_Claim_} 
Let $\c H\in Hyp$, $\c H= (I, J, K, (\cdot,\cdot))$,
$a=P_2(\c H)$, $b=P_3(\c H)$.
Then the following endomorphisms of $A$
coinside:

\[ [L_a,\Lambda_b]=ad I  \]

{\bf Proof:} See \cite{_so5_on_cohomo_}. $\;\;\blacksquare$

\hfill

This finishes the proof of \ref{_str_of_g_0_Proposition_} (i).
To prove \ref{_str_of_g_0_Proposition_} (ii), we recall the
following linear-algebraic construction. Let $V$ be a linear space
equipped with non-degenerate symmetric bilinear form $(\cdot,\cdot)$,
and $\c L\subset V$ be a 2-dimensional plane in $V$,
such that the restriction of $(\cdot,\cdot)$ to $\c L$
is non-degenerate. Let $\c L^\bot$ be the orthogonal complement
of $V$ to $\c L$.
Let $T_{\c L}$ be the set of all non-trivial skew-symmetric endomorphisms
of $V$ which vanish on $\c L^\bot$. As we have seen
previously, all elements of $T_{\c L}$ are proportional.

Let $Gr^\circ$ be the space of all 2-dimensional planes
$\c L\subset V$ such that the restriction
$(\cdot,\cdot)_{\c H}\restrict{\c L}$ is non-degenerate.

\hfill

\claim \label{_T_L_generate_SO_Claim_} 
The linear span of the union

\[ \bigcup\limits_{\c L\in Gr^\circ} T_{\c L} \]
coinsides with $\goth{so}(V)$.

{\bf Proof:} Clear. $\;\;\blacksquare$

\hfill

By \ref{_T_L_generate_SO_Claim_}, to prove that $\c M_2=\goth{so}(V)$ it
is sufficient to show that for a Zariski dense set $\c E\subset Gr^\circ$,
we have

\[
  \forall \c L\in \c E, \;\; T_{\c L}\subset \c M_2.
\]
For $\c E$, we take the set $Gr^\bullet$ of
\ref{_D_in_Q_c(Comp)_Claim_}. Since $Gr^\bullet = Q_c(Comp)$
is Zariski dense in $Gr$, it is sufficient to show that for all
$\c L\in Q_c(Comp)$, we have $T_{\c L}\subset \c M_2$. This is
implied by \ref{_[L_2,Lambda_3]=adI_Claim_} and the following
statement.

\hfill

\claim \label{_T_Q(I)_is_[L,Lambda]_Claim_} 
In assumptions of \ref{_[L_2,Lambda_3]=adI_Claim_},
the following two sets coinside:

\[ \bigg\{\lambda [L_a,\Lambda_b],\;\lambda\in\R\backslash 0\bigg\} =
   T_{Q_c(I)}.
\]
{\bf Proof:} Since all elements of $T_{Q_c(I)}$ are proportional
and $[L_a,\Lambda_b]\neq 0$, it is sufficient to show
that $[L_a,\Lambda_b]\in T_{Q_c(I)}$. Let $\c L\subset V$ be the linear
span of $a$ and $b$. Obviously from definition, $Q_c(I)=\c L$.
Let $\c L^\bot$ be the orthogonal complement to $\c L$ in $V$.
We need to show that the restriction of $[L_a,\Lambda_b]$
to $V=H^2(M, \R)$ vanishes on $\c L^\bot$. Consider the
Hodge decomposition on $H^2(M)$ associated with the complex
structure $I$. By definition,

\[ \c L =
   \bigg( H^{2,0}_I(M)\oplus H^{0,2}_I(M)\bigg) \cap H^2(M, \R).
\]
Hodge-Riemann relations imply that
\[
   \c L^\bot= H^{1,1}_I(M)\cap H^2(M, \R)
\]
(see Section
\ref{_Period_and_Hodge_Riemann_Section_} for details). By
defintion of $ad I$,  we have
$ad I(H^{1,1}_I(M))=0$. Since $[L_a,\Lambda_b] = ad I$, we
obtain that $[L_a,\Lambda_b]\in T_{Q_c(I)}$. This proves
\ref{_T_Q(I)_is_[L,Lambda]_Claim_}, and consequently,
proves \ref{_str_of_g_0_Proposition_} (ii). $\;\;\blacksquare$

\hfill

Using the fact that

\[ \{ a,b\in A_2\times A_2 \;\; | \; \; a\bullet b \} \]
is Zariski dense in $A_2\times A_2$, we can derive
from \ref{_T_Q(I)_is_[L,Lambda]_Claim_} the following
corollary:

\hfill

\corollary \label{_[L,Lambda]_in_T_a,b_Corollary_} 
Let $a,b\in A_2$ be elements of Lefschetz type. Let
$\c L\subset A_2$ be a plane generated by $a$ and $b$
and $\c L^\bot$ be its orthogonal complement in $A_2=V$.
Let $[L_a,\Lambda_b]\bigg|_{{}_{V}}^\circ$ be a traceless part of
$[L_a,\Lambda_b]\bigg|_{{}_{V}}\in End(V)$. Then
$[L_a,\Lambda_b]\bigg|_{{}_{V}}^\circ$ vanish on $\c L^\bot$.

{\bf Proof:} Using the argument with Zariski dense sets,
we see that it is sufficient to check
\ref{_[L,Lambda]_in_T_a,b_Corollary_}
for $a,b$ such that $a\bullet b$, and $(b,b)_{\c H}\neq 0$.
When $(a,b)_{\c H} =0$,
\ref{_[L,Lambda]_in_T_a,b_Corollary_}
is a direct consequence of \ref{_T_Q(I)_is_[L,Lambda]_Claim_}.
If $(a,b)_{\c H} \neq 0$, take

\[ x:= a - \frac{(a,b)_{\c H}}{(b,b)_{\c H}} b. \]
Clearly, $(x,b)_{\c H} =0$.
The traceless part of $[L_b, \Lambda_b]\bigg|_{{}_{V}}= H\restrict{V}$ is zero.
Therefore,

\[ [L_a,\Lambda_b]\bigg|_{{}_{V}}^\circ =
   [L_x,\Lambda_b]\bigg|_{{}_{V}}^\circ
\]
This reduces \ref{_[L,Lambda]_in_T_a,b_Corollary_}
to the case $(a,b)_{\c H} =0$.
$\;\;\blacksquare$

\hfill

To prove \ref{_str_of_g_0_Proposition_} (iii), we notice that from
the proof of \ref{_str_of_g_0_Proposition_} (ii) it follows also that
$t:\; \g_0(A)\arrow \goth{so}(V)$ is an epimorphism. Therefore, it is
sufficient to find an element $e\in \g_0(A)$ such that $t(e)=0$,
$s(e)\neq 0$. This element is a Hodge endomorphism
$H\in \g_0(A)\subset End(A)$. We proved
\ref{_str_of_g_0_Proposition_} (iii). $\;\;\blacksquare$


\section{Calculation of a zero graded part of the structure Lie
algebra of the cohomology of a hyperkaehler manifold, part II.}
\label{_compu_g_0_part_2_Section_}


We work in assumptions of the previous section.

\hfill

\theorem \label{_g_0_computed_Theorem_} 
In assumptions of \ref{_str_of_g_0_Proposition_},
consider the homomorphism

\[
   u:\; \g_0(A)\arrow \goth{so}(V) \oplus k H,
\]
$u=t\oplus s$. Then $u$ is an isomorphism.

{\bf Proof:}  We use the following technical result:

\hfill

\proposition \label{_[L_a,Lambda_b]-[Lambda_a,L_b]_Proposition_} 
Let $A$ be a Lefschetz-Frobenius algebra of cohomology of a compact
simple hyperkaehler manifold. Let $a, b\in A_2$ be elements of Lefschetz
type. Then

\[ {{(b,b)_{\c H}}} [L_a, \Lambda_b] -
   {{(a,a)_{\c H}}} [\Lambda_a,L_b] =2(a,b)_{\c H} H,
\]
where $(\cdot,\cdot)_{\c H}$ is a normalized Hodge-Riemann pairing.

{\bf Proof:} First of all, we prove the following lemma.

\hfill

\lemma \label{_Lambda_additive_for_a_bullet_b_Lemma_}
In assumptions of
\ref{_[L_a,Lambda_b]-[Lambda_a,L_b]_Proposition_},
let $a,b\in A_2$, $a\bullet b$.
Then $a+b$ is also of Lefschetz type, and

\[ {(a+b,a+b)_{\c H}}\Lambda_{a+b} =
   {(a,a)_{\c H}}\Lambda_a + {(b,b)_{\c H}}\Lambda_b.
\]

{\bf Proof:} Let $\c H\in Hyp$ be a hyperkaehler structure,
such that there exist $x_i, y_i\in \R$, $i=1,2,3$, such that
\[
   a =\sum x_i\omega_i, \;\; b=\sum y_i\omega_i,
\]
where $\omega_i=P_i(\c H)$, $i=1,2,3$. Such $\c H$ exists
by definition of the relation $a \bullet b$.
Clearly,
\[ (\omega_1,\omega_1)_{\c H}= (\omega_2,\omega_2)_{\c H} =
   (\omega_3,\omega_3)_{\c H}.
\]
Let $c:= {(\omega_1,\omega_1)_{\c H}}$.
We are going to show that

\begin{equation}\label{_Lambda_a_as_lin_combi_Lambda_omega_Equation_}
   \Lambda_a = c\frac{\sum x_i\Lambda_{\omega_i}}{{(a,a)_{\c H}}.}
\end{equation}
First of all, we notice that for all
$a =\sum x_i\omega_i, \;\; a\neq 0,$ the triple

\[
   L_a,\;\; H,\;\; \Lambda :=
   \frac{\sum x_i\Lambda_{\omega_i}}{{\sum x_i^2}}
\]
is a Lefschetz triple.
This can be shown by an easy calculation which uses
\eqref{_so5_relations_Equation_}. On the other hand,
$\omega_i$ are orthogonal with respect to $(\cdot,\cdot)_{\c H}$,
and therefore $(a,a)_{\c H}= c\sum x_i^2$. Therefore,

\[ c\frac{\sum x_i\omega_i}{{(a,a)_{\c H}}} =
   \frac{\sum x_i\omega_i}{{\sum x_i^2}}
\]
Now, \eqref{_Lambda_a_as_lin_combi_Lambda_omega_Equation_}
immediately implies \ref{_Lambda_additive_for_a_bullet_b_Lemma_},
as an easy calculation shows:

\[ {(a+b,a+b)_{\c H}}\Lambda_{a+b} =^*
   c\sum (x_i+y_i)\Lambda_{\omega_i} =
\]

\[   = c\sum x_i\Lambda_{\omega_i}+c\sum y_i\Lambda_{\omega_i}
   =^* {(a,a)_{\c H}}\Lambda_{a} +{(b,b)_{\c H}}\Lambda_{b},
\]
where $=^*$ marks an application of
\eqref{_Lambda_a_as_lin_combi_Lambda_omega_Equation_}.
This proves \ref{_Lambda_additive_for_a_bullet_b_Lemma_}.
$\;\;\blacksquare$

\hfill

Since

\[ \{ (a,b)\in A_a \;\; |\;\; a\bullet b \} \]
is Zariski dense in $A_2$,
\ref{_Lambda_additive_for_a_bullet_b_Lemma_}
has the following corollary:

\hfill

\corollary \label{_Lambda_additive_Corollary_} 
Let $a,b\in A_2$, be the elements of Lefschetz type, such that
$a+b$ is also of Lefschetz type. Then

\[ {(a+b,a+b)_{\c H}}\Lambda_{a+b} =
   {(a,a)_{\c H}}\Lambda_a + {(b,b)_{\c H}}\Lambda_b.
\]

$\blacksquare$

\hfill

\ref{_Lambda_additive_Corollary_} implies the following interesting
result. The set $S\subset A_2$ of Lefschetz elements
is preserved by homotheties.
Therefore we may speak of homogeneous functions from
$S$ to some linear space. Consider the function
$\Lambda:\; S\arrow \g_{-2}$, $x\arrow \Lambda_2$.
Clearly, this map is homogeneous of degree $-1$.
Therefore, the map $\tilde r:\;S\arrow \g_{-2}$,
$x\arrow (x,x)_{\c H} \Lambda_x$ is homogeneous of degree 1.
By \ref{_Lambda_additive_Corollary_}, $\tilde r$ is linear
on $S$. Therefore, $\tilde r$ is a restriction of a linear
map $r:\; A_2 \arrow \g_{-2}$

\hfill

\corollary \label{_g_-2_is_quotie_of_A_2_Corollary_} 
The map $r:\; A_2\arrow \g_{-2}$ is a surjection
of linear spaces.

{\bf Proof:} For all Lefschetz-Frobenius
algebras $C$, $\g_{-2}(C)$ is spanned by $\Lambda_x$
for $x$ of Lefschetz type (\cite{_Lunts-Loo_}).
$\;\;\blacksquare$

\hfill

To prove \ref{_[L_a,Lambda_b]-[Lambda_a,L_b]_Proposition_},
we need only the formula
\eqref{_Lambda_a_as_lin_combi_Lambda_omega_Equation_}.
Since

\[ \{ (a,b)\in A_a \;\; |\;\; a\bullet b \} \]
is Zariski dense in $A_2$, it is sufficient to prove
\ref{_[L_a,Lambda_b]-[Lambda_a,L_b]_Proposition_}, for $a,b$
such that $a\bullet b$. Consider notation introduced
in the proof of \ref{_Lambda_additive_for_a_bullet_b_Lemma_}.
The formula \eqref{_Lambda_a_as_lin_combi_Lambda_omega_Equation_}
together with relations \eqref{_so5_relations_Equation_}
implies that

\begin{equation}\label{_[Lambda_a,L_b]_via_x_y_Equation_}
   [ \Lambda_a, L_b ] =
   \bigg[ \frac{\sum x_i \Lambda_{\omega_i}}{\sum x_i^2},
   \sum y_i L_{\omega_i} \bigg]
   = \frac{1}{\sum x_i^2} \bigg(-\sum x_i y_i H +
   \sum\limits_{i\neq j} x_i y_j [\Lambda_{\omega_i}, L_{\omega_j}]\bigg),
\end{equation}
and

\begin{equation}\label{_[L_a,Lambda_b]_via_x_y_Equation_}
   [ L_a, \Lambda_b ] =
   -\bigg[ \frac{\sum y_i \Lambda_{\omega_i}}{\sum y_i^2},
   \sum x_i L_{\omega_i} \bigg]
   = \frac{1}{\sum y_i^2} \bigg(\sum x_i y_i H -
   \sum\limits_{i\neq j} x_i y_j [L_{\omega_i}, \Lambda_{\omega_j}]\bigg).
\end{equation}
Let $c\in \R$ be the constant defined in the proof
of \ref{_Lambda_additive_for_a_bullet_b_Lemma_}.
We have

\[ (a,a)_{\c H}= c \sum x_i^2,\;\; (b,b)_{\c H}= c \sum y_i^2,\;\;
   (a,b)_{\c H} = c \sum x_i y_i.
\]
Making the corresponding substitutions in
\eqref{_[Lambda_a,L_b]_via_x_y_Equation_},
\eqref{_[L_a,Lambda_b]_via_x_y_Equation_}, we obtain

\begin{equation}\label{_[Lambda_a,L_b]_via_()_H_and_x_i_Equation_}
   [ \Lambda_a, L_b ]
   = \frac{1}{(a,a)_{\c H}} \bigg(-(a,b)_{\c H} H +
   \frac{1}{c}
   \sum\limits_{i\neq j} x_i y_j [\Lambda_{\omega_i}, L_{\omega_j}]\bigg)
\end{equation}
and
\begin{equation}\label{_[L_a,Lambda_b]_via_()_H_and_x_i_Equation_}
   [ L_a, \Lambda_b ]
   = \frac{1}{(b,b)_{\c H}} \bigg((a,b)_{\c H} H +
   \frac{1}{c}
   \sum\limits_{i\neq j} x_i y_j [L_{\omega_i}, \Lambda_{\omega_j}]\bigg).
\end{equation}
A linear combination of these equations yields

\[ {(b,b)_{\c H}}[L_a, \Lambda_b]-{(a,a)_{\c H}}[\Lambda_a, L_b] =
   (a,b)_{\c H} H +
\]
\[ +\frac{1}{c}
   \sum\limits_{i\neq j} x_i y_j
   \bigg([\Lambda_{\omega_i}, L_{\omega_j}] -
         [L_{\omega_i}, \Lambda_{\omega_j}] \bigg)
\]
To prove \ref{_[L_a,Lambda_b]-[Lambda_a,L_b]_Proposition_}
it remains to show that

\[ \sum\limits_{i\neq j} x_i y_j
   \bigg([\Lambda_{\omega_i}, L_{\omega_j}] -
         [L_{\omega_i}, \Lambda_{\omega_j}] \bigg) =0,
\]
which follows from the relation

\[ [\Lambda_{\omega_i}, L_{\omega_j}] =
         [L_{\omega_i}, \Lambda_{\omega_j}], \;\;i \neq j
\]
from \eqref{_so5_relations_Equation_}.
\ref{_[L_a,Lambda_b]-[Lambda_a,L_b]_Proposition_}
is proven. $\;\;\blacksquare$

\hfill

Let $s:\; \g_0\arrow kH$ be a Lie algebra homomorphism
of \ref{_str_of_g_0_Proposition_}.
Either of equations \eqref{_[L_a,Lambda_b]_via_()_H_and_x_i_Equation_}
and \eqref{_[Lambda_a,L_b]_via_()_H_and_x_i_Equation_} implies the
following useful corollary:

\hfill

\corollary \label{_s_of_[L_a,Lambda_b]_Corollary_} 
Let $a,b \in S$.  Then

\[ s([L_a,\Lambda_b]) = \frac{(a,b)_{\c H}}{(b,b)_{\c H}} H. \]

$\blacksquare$

\hfill

Now we can easily prove \ref{_g_0_computed_Theorem_}.
Let $n=\dim V$. Clearly, $\dim \goth{so}(V)=\frac{n(n-1)}{2}$.
Since $u:\; \g_0\arrow \goth{so}(V)\oplus k$ is an epimorphism by
\ref{_str_of_g_0_Proposition_}, it is sufficient to show that
$\dim \g_0\leq \frac{n(n-1)}{2}+1$. The element $H$ belongs to the center
of $\g_0$. Let $\bar \g_0:= \g_0/k H$ be the quotient Lie algebra, and
$q:\; \g_0\arrow \bar\g_0$ be the quotient map. Let $S\subset A_2$ be the
set of all elements of Lefschetz type. By
\ref{_g_0_gener_by_[L,Lambda]_Lemma_}, the space
$\bar \g_0$ is spanned by all vectors
$q([L_a,\Lambda_b])$, where $a,b\in S$.
Denote the map $S\times S\arrow \bar\g_0$,
$a,b \arrow q([L_a,\Lambda_b])$ by
$\tilde\nu:\; S\times S \arrow \bar\g_0$. Let
$\nu:\; S\times S\arrow \bar\g_0$,
\[
   \nu(a,b):= \frac{q([L_a,\Lambda_b])}{(b,b)_{\c H}}.
\]
Consider the Zariski open set $S\subset A_2$ as a space with an
associative commutative group structure, which is defined by rational
maps. In other words, $S$ is equipped with an addition, which
is defined not everywhere, but in a Zariski open subset of
$S$. This addition is induced from $A_2\supset S$, which is a linear space.
By \ref{_Lambda_additive_Corollary_}, the map $\nu$ is bilinear with respect
to this operation. Consider $\nu$ as a rational map from $A_2\times A_2$
to $\bar \g_0$. This rational map is also bilinear. An easy check
shows that a linear rational map of linear spaces is defined everywhere.
Hence, $\nu$ can be uniquely lifted to a bilinear
map

\[
   \nu:\; A_2\times A_2 \arrow \bar\g_0,
\]
such that the square

\[ \begin{array}{ccc}
     S\times S & \hookrightarrow & A_2\times A_2 \\[3mm]
     \bigg\downarrow \nu &&\bigg\downarrow \nu \\[3mm]
     \bar\g_0 &\stackrel{Id}{\arrow} &\bar\g_0
   \end{array}
\]
is commutative.
By \ref{_[L_a,Lambda_b]-[Lambda_a,L_b]_Proposition_},
the bilinear map $\nu:\; A_2\times A_2\arrow \bar\g_0$ is
skew-symmetric. Let $\eta$ be the corresponding linear map from the
exterrior square of $A_2$ to $\bar \g_0$:

\[ \eta:\;
   \bigwedge^2 A_2\arrow \bar\g_0.
\]
Clearly, $\nu(A_2\times A_2) = \eta(\bigwedge^2 A_2)$.
As \ref{_g_0_gener_by_[L,Lambda]_Lemma_} implies,
$\bar\g_0$ is generated by the image of $\nu(A_2\times A_2)$.
Hence, $\eta$ is an epimorphism. Therefore,
$\dim \bar\g_0 \leq \dim \bigwedge^2 A_2 = \frac{n(n-1)}{2}$.
We obtained an upper bound on $\dim\g_0$:

\[
   \dim \g_0\leq \frac{n(n-1)}{2} + 1.
\]
Since the Lie algebra $\goth{so}(V)\oplus k$ has the same dimension,
$\frac{n(n-1)}{2}+1$ and the map $u:\; \g_0\arrow \goth{so}(V)\oplus k$
is an epimorphism as we have seen previously, the map $u$
is an isomorphism. \ref{_g_0_computed_Theorem_} is proven.
$\;\;\blacksquare$

\hfill

Consider the map $r:\; A_2\arrow \g_{-2}$ constructed in
\ref{_g_-2_is_quotie_of_A_2_Corollary_}. Let $r':\; A_2\arrow \g_2$
be a standard isomorphism. Then $\nu$ can be defined
as a composition of $r'\times r:\; A_2\times A_2 \arrow \g_2\times \g_{-2}$
and a commutator map $[,]:\; \g_2\times \g_{-2}\arrow \g_0$.
The map $\nu$ was described in terms of the standard bilinear
map $A_2\times A_2\arrow \Lambda^2A_2$. This description
easily implies that for all $y\in A_2$ there exists
$x\in A_2$ such that $\nu(x,y)\neq 0$. Since
$\nu=r'\times r\circ [,]$, the map $r$ has no kernel.
It is epimorphic by \ref{_g_-2_is_quotie_of_A_2_Corollary_}.
We obtained the following statement:

\hfill

\corollary \label{_g_-2_is_A_2_Corollary_} 
Consider the linear map $r:\; A_2\arrow \g_{-2}$ constructed in
\ref{_g_-2_is_quotie_of_A_2_Corollary_}. Then $r$ is
an isomorphism of linear spaces.

$\;\;\blacksquare$

\hfill

\ref{_g_0_computed_Theorem_} gives a nice insight on the
Hodge structures on $H^*(M)$ corresponding to various
complex structures $I\in Comp$. As elsewhere,
for each $I\in Comp$ we define an endomorphism $ad I\in End(A)$,
$ad I(\omega)= (p-q)\1 \omega$ for all $\omega \in H^{p,q}(M)$.
Let $\goth M\subset End(A)$ be a Lie algebra generated by
$ad I$ for all $I\in Comp$.

\hfill

\lemma \label{_M_acts_by_deriva_Lemma_} 
The Lie algebra $\goth M$ acts on $A$ by derivations:
for all $m\in \goth M$, $a,b \in A$, we have
\[
    m(ab) = m(a) b+ am(b).
\]

{\bf Proof:} For all $I\in Comp$, the operator $ad I$ is a derivation,
as a calculation shows.
A commutator of derivations is a derivation by obvious reasons.
$\;\; \blacksquare$

\hfill

\theorem\label{_g_0_is_Mumf_Tate_Theorem_} 
The following Lie subalgebras of $End(A)$ coinside:

\[
   \g_0(A)\cong \goth M\oplus kH.
\]
{\bf Proof:} The inclusion $\goth M\subset \g_0(A)$ is implied by
\ref{_[L_2,Lambda_3]=adI_Claim_}. The inclusion
$\g_0(A)\subset \goth M\oplus kH$ is proven as follows.
We have seen that $\g_0(A)$ is linearly spanned by $H$ and
endomorphisms $[L_a,\Lambda_b]$, where $a\bullet b$, $(a,b)_{\c H}=0$,
$(a,a)_{\c H}=(b,b)_{\c H}\neq 0$. Let $\c H'$ be a hyperkaehler
structure such that $a\bullet_{\c H'} b$. Applying
\ref{_action_SO(3)_on_Hyp_via_periods_Lemma_}, we
obtain a hyperkaehler structure $\c H$, $\c H$ equivalent to
$\c H'$, such that $a=P_2(\c H)$ and
$b=P_3(\c H)$. Then, \ref{_[L_2,Lambda_3]=adI_Claim_}
implies that $[L_a,\Lambda_b]\in \goth M$. $\;\;\blacksquare$

\hfill

By trivial reasons, $H$ acts on $A$ as a derivation. Therefore,
\ref{_g_0_is_Mumf_Tate_Theorem_} together with
\ref{_M_acts_by_deriva_Lemma_} implies the following:

\hfill

\corollary \label{_g_0_derivatives_Corollary_} 
The Lie algebra $\g_0\subset End(A)$ acts on $A$ by derivations.

$\;\;\blacksquare$


\section[Computing the structure Lie algebra for the
cohomology of a hyperkaehler manifold, part \ II.]
{Computing the structure Lie algebra for the
cohomology of a hyperkaehler manifold, \\part \ II.}
\label{_computing_g_for_hyperk_pt-2_Section_}


In this section, we prove the isomorphism $\g\cong \goth{so}(V, +)$.
This is done as follows. In previous sections, we
have computed dimensions of the components
of $\g\cong \g_{-2}\oplus \g_0\oplus\g_2$. Let $n:= \dim V$.
\ref{_g_-2_is_A_2_Corollary_} implies that $\dim \g_2=n$,
\ref{_g_-2_is_A_2_Corollary_} implies that $\dim\g_{-2}=n$,
and \ref{_g_0_computed_Theorem_} implies that
$\dim\g_0 = \frac{n(n-1)}{2}+1$. Therefore,

\[
   \dim\g= \frac{n(n-1)}{2}+1 +2n = \frac{(n+2)(n+1)}{2}.
\]
A trivial computation yields

\[
   \dim \goth{so}(V,+)= \frac{(n+2)(n+1)}{2}.
\]
We see that dimensions of $\g$ and $\goth{so}(V, +)$ are equal;
it is easy to see that $\g$ is isomorphic to $\goth{so}(V, +)$
as a graded linear space.
We construct an isomorphism of graded linear spaces
$U:\; \g \arrow \goth{so}(V, +)$ and prove that
it commutes with the Lie algebra operation.

\hfill

Let $B= B_0\oplus B_2\oplus B_4$ be the minimal graded Frobenius
algebra associated with $V$, $(\cdot,\cdot)_{\c H}$.
By definition, $\g(B)= \goth{so}(V,+)$. We are going to construct
an isomorphism of linear spaces $U=U_{-2}\oplus U_0\oplus U_2$,

\[
    U:\; \g(A)\arrow \g(B),
    \;\; U_{2i}:\; \g_{2i}(A)\arrow \g_{2i}(B).
\]
We have canonical isomorphisms

\[ u_B:\; \g_0(B)\arrow \goth{so}(V)\oplus kH,\;\;
   u_A:\; \g_0(A)\arrow \goth{so}(V)\oplus kH.
\]
These isomorphisms yield a natural Lie algebra isomorphism
$U_0:\; \g_0(A)\arrow \g_0(B)$. The homomorphism
$U_2:\; \g_2(A)\arrow \g_2(B)$ is provided by the natural
isomorphism $\g_2(A)\cong A_2$ (\ref{_g_2_is_A_2_Corollary_}),
which exists for every Lefschetz-Frobenius algebra of Jordan type.
To construct the isomorphism $U_{-2}:\; \g_{-2}(A)\arrow \g_{-2}(B)$,
we use \ref{_g_-2_is_A_2_Corollary_}. According to this statement,
$\g_{-2}(A)$ is naturally isomorphic to $V$.
The natural isomorphism $\g_{-2}(B)\cong V$ is constructed in
\ref{_calculation_of_g(A)_for_minim_Theorem_}. Composing
these isomorphisms, we obtain $U_{-2}$.
Now, the isomorphism $\g(A)\cong \goth{so}(V, +)$
is implied by the following proposition:

\hfill

\proposition \label{_U_is_Lie_homomo_Proposition_} 
The map $U:\; \g(A)\arrow \g(B)$
is a homomorphism of Lie algebras.

{\bf Proof:} By our construction, the restriction of $U$ to $\g_0(A)$
is a homomorphism of Lie algebras. Therefore, it suffices to
check that

\begin{equation}\label{_U_commu_with_commutator_Equation_}
   U([X,Y])=[U(X),U(Y)]
\end{equation}
in the following cases:

(i) $X\in \g_2$, $Y\in \g_{-2}$

(ii) $X\in \g_0$, $Y\in \g_2$

(iii) $X\in \g_0$, $Y\in \g_{-2}$.

\hfill

We start the proof of \eqref{_U_commu_with_commutator_Equation_}
with (i). We represent $u_A:\; \g_0(A)\arrow \goth{so}(V)\oplus k H$
as a sum $u_A=t_A\oplus s_A$, where $t_A:\; \g_0(A)\arrow \goth{so}(V)$
and $s_A:\; \g_0(A)\arrow k H$ are components of $u_A$. Consider
the similar decomposition of $u_B:\; \g_0(B)\arrow \goth{so}(V)\oplus k H$,
$u_B=t_B\oplus s_B$.
Let $S_A\subset A_2=V$ be the set of all elements of Lefschetz type
in $A$, and $S_B\subset B_2=V$ be the set of all elements of Lefschetz type
in $B$. Let $S:= S_A\cap S_B$. For $x\in S_A$ ($S_B)$, we denote
by $L^A_x$, $\Lambda^A_x$ ($L^B_x$, $\Lambda^B_x$) the corresponding
elements in $\g(A)$ ($\g(B)$). Clearly, $S$ is Zariski open in $V$.
Therefore, to prove \eqref{_U_commu_with_commutator_Equation_}
in case (i) it is sufficient to show that for all $x,y \in S$,

\begin{equation}\label{_U_[L,Lambda]=[U(L),U(Lambda)]_Equation_}
   U([L^A_x, \Lambda^A_y]) =
   [ U(L^A_x), U(\Lambda^A_y) ].
\end{equation}

Checking the definition of $U$, one can easily see that
$U(L_x^A)= L_x^B$ and $U(\Lambda_x^A)= \Lambda_x^B$.
Therefore, \eqref{_U_[L,Lambda]=[U(L),U(Lambda)]_Equation_}
is equivalent to

\begin{equation}\label{_U_[L,Lambda]=[L^B,Lambda^B]_Equation_}
   U([L^A_x, \Lambda^A_y]) =
   [L^B_x, \Lambda^B_y ].
\end{equation}

By definition of $U$, \eqref{_U_[L,Lambda]=[L^B,Lambda^B]_Equation_}
is equivalent to

\[
   u_A([L^A_x, \Lambda^A_y]) = u_B([L^B_x, \Lambda^B_y ]).
\]
Using the decomposition of $u_A$, $u_B$, we obtain that
this equation is implied by the following two relations:

\begin{equation}\label{_t_commu_with_commutato_Equation_}
   t_A([L^A_x, \Lambda^A_y]) = t_B([L^B_x, \Lambda^B_y ]),
\end{equation}
\begin{equation}\label{_t_commu_with_commutato_another_Equation_}
   s_A([L^A_x, \Lambda^A_y]) = s_B([L^B_x, \Lambda^B_y ]).
\end{equation}
The second of these relations is implied by the equation

\[ s_A([L^A_x, \Lambda^A_y]) =
   \frac{(x,y)_{\c H}}{(y,y)_{\c H}} H
\]
(\ref{_s_of_[L_a,Lambda_b]_Corollary_}).
We proceed to prove \eqref{_t_commu_with_commutato_Equation_}.
Consider the action of the operators
$t_A([L^A_x, \Lambda^A_y]),
t_B([L^B_x, \Lambda^B_y ]) \in \goth{so}(V)$
on $V$. Let $\c L$ be the two-dimensional plane
generated by $x,y\in V$. Let $\c L^\bot$ be its
orthogonal complement. By
\ref{_[L,Lambda]_in_T_a,b_Corollary_},
$t_A([L^A_x, \Lambda^A_y])$ acts as zero on $\c L^\bot$.
By \ref{_Lambda_vanish_Corollary_},
$t_B([L^B_x, \Lambda^B_y])$ also vanish on $\c L^\bot$.
The space of skew-symmetric endomorphisms of $V$ which
vanish on $\c L^\bot$ is one-dimensional. Hence,
the operators $t_A([L^A_x, \Lambda^A_y])$ and
$t_B([L^B_x, \Lambda^B_y])$ are proportional.
To prove that they are equal we have to compute
only the coefficient of proportionality.

\hfill

Denote the result of application of $\xi\in \goth{so}(V)$ to $x\in V$ by
$\xi x$. To prove \eqref{_t_commu_with_commutato_Equation_}
it is sufficient to show that

\begin{equation}\label{_L_x_Lambda_y_to_y_Equation_}
   t_A([L_x^A,\Lambda_y^A])y=t_B([L_x^B,\Lambda_y^B])y\neq 0.
\end{equation}
In the case when $x,y$ are collinear,

\[
   t_A([L_x^A,\Lambda_y^A])=t_B([L_x^B,\Lambda_y^B])=0.
\]
(see \eqref{_t(L_y_Lambda_y)=0_Equation_}).
Therefore, in this case \eqref{_L_x_Lambda_y_to_y_Equation_}
is not true. However, \eqref{_t_commu_with_commutato_Equation_}
is vacuously true in this case. Therefore, to prove
\eqref{_t_commu_with_commutato_Equation_} it is sufficient
to prove \eqref{_L_x_Lambda_y_to_y_Equation_} assuming
that $x$ and $y$ are not collinear.

\hfill

Let $x,y\in V$ be two vectors which are not collinear.
We prove \eqref{_L_x_Lambda_y_to_y_Equation_} as follows.
Denote the Hodge operators in $A$, $B$ by $H^A$, $H^B$
respectively. By definition of the Lefschetz triple,

\[
    [L_y^A,\Lambda_y^A])= H^A, \;\;[L_y^B,\Lambda_y^B] = H^B.
\]
This implies that

\begin{equation}\label{_t(L_y_Lambda_y)=0_Equation_}
   t_A([L_y^A,\Lambda_y^A])=t_B([L_y^B,\Lambda_y^B])=0
\end{equation}
Let $\lambda\in R$.
Since the expressions $t_A([L_x^A,\Lambda_y^A])$,
$t_B([L_x^B,\Lambda_y^B])$ are bilinear by $x$,
we have

\[
    t_A([L_{x+\lambda y}^A,\Lambda_y^A])=t_A([L_x^A,\Lambda_y^A])
\]
and

\[
    t_B([L_{x+\lambda y}^B,\Lambda_y^B])=t_B([L_x^B,\Lambda_y^B]).
\]
Therefore, \eqref{_L_x_Lambda_y_to_y_Equation_}
is equivalent to

\[
    t_A([L_{x+\lambda y}^A,\Lambda_y^A])y=
    t_B([L_{x+\lambda y}^B,\Lambda_y^B])y\neq 0.
\]
By \ref{_el-t_with_non_zero_square_Lefschetz_Lemma_},
$z\in S_B$ if and only if $(z,z)_{\c H}\neq 0$. Since
$S\subset S_B$, the number $(y,y)_{\c H}$ is non-zero.

Take $\lambda=\frac{(x,y)_{\c H}}{(y,y)_{\c H}}$.
Then $(x+\lambda y,y)_{\c H}=0$. Replacing $x$
by $x+\lambda y$, we see that it is sufficient
to prove \eqref{_L_x_Lambda_y_to_y_Equation_} in the
case when $(x,y)_{\c H}=0$.

Let $\mu\in \R$, $\mu\neq 0$.
If we replace $x$  by a vector $\mu x$,
both sides of \eqref{_L_x_Lambda_y_to_y_Equation_}
are multiplied by the number $\mu$. Choosing the
appropriate coefficient $\mu$, we may assume that
$(x,x)_{\c H}=(y,y)_{\c H}>0$. We obtained that we may prove
\eqref{_L_x_Lambda_y_to_y_Equation_} under the following set
of assumptions:

\[
  (x,x)_{\c H}=(y,y)_{\c H}>0, \;\; (x,y)_{\c H}=0.
\]
This is implied by the following lemma.

\hfill

\lemma \label{_[L_x,Lambda_y]_to_y_is_-2x_Lemma_} 
Let $x,y\in S$, $(x,x)_{\c H}=(y,y)_{\c H}>0$,
$(x,y)_{\c H}=0$. Then

\hfill

(i) $t_A([L_x^A,\Lambda_y^A])y =-2x$, and

(ii) $t_B([L_x^B,\Lambda_y^B])y =-2x$.

\hfill

{\bf Proof:} (i) Let

\[ T:= \bigg \{ (x,y)\in S \;\; | \;\; (x,x)_{\c H}=(y,y)_{\c H},\;\;
    (x,y)_{\c H}=0 \bigg\}.
\]
Let
\[ T^\bullet:= \bigg \{ (x,y)\in S \;\; | \;\;
    (x,x)_{\c H}=(y,y)_{\c H}, \;\;
    (x,y)_{\c H}=0, \;\;x \bullet y \bigg\}.
\]
A standard argument with periods and comparing dimensions
implies that $T^\bullet$ is Zariski dense in $T$.
Therefore, we may prove (i) assuming that $x\bullet y$.
Let $\tilde {\c H}$ be a hyperkaehler structure such that
$x\bullet_{\tilde{\c H}} y$. Using
\ref{_action_SO(3)_on_Hyp_via_periods_Lemma_}, we replace
$\tilde {\c H}$ by an equivalent hyperkaehler structure
$\c H=(I,J, K, (\cdot,\cdot))$ such that $P_2(\c H)=x$,
$P_3(\c H) =y$. In this case

\[
   [ L^A_x, L^A_y] = ad I
\]
(\ref{_[L_2,Lambda_3]=adI_Claim_}). Let $\Omega:= x+ \1 y$.
By definition of $ad I$, $adI (\Omega) = 2\1 \Omega$.
This immediately implies
\ref{_[L_x,Lambda_y]_to_y_is_-2x_Lemma_} (i).

\hfill

{\bf Proof of \ref{_[L_x,Lambda_y]_to_y_is_-2x_Lemma_} (ii):}
By definition, for all $a, b\in B_2$,

\begin{equation}\label{_ab_in_B_Equation_}
   ab= (a,b)_{\c H}\Omega_B
\end{equation}
for a fixed vector $\Omega_B\in B_4$. Therefore, $L_x^B y=0$.
We obtain that

\begin{equation}\label{_commut_applied_to_y__in_B_Equation_}
   [ L_x^B,\Lambda_y^B]y = L_x^B\Lambda_y^B y.
\end{equation}
Let ${\Bbb I}_B\in B_0$ be the unit in $B$. Then
$\Lambda_y^B y= \Lambda_y^B L_y^B {\Bbb I}_B$. Since
$[L_y^B,\Lambda_y^B] =H^B$ and $\Lambda_y^B {\Bbb I}_B=0$,
we have

\[ \Lambda_y^B L_y^B {\Bbb I}_B =
   - H^B (\Bbb I) = -2 {\Bbb I}_B.
\]
Using \eqref{_ab_in_B_Equation_}, we can easily check that

\[
   [L_x^B,\Lambda_y^B]{\Bbb I}_B=0.
\]
Therefore,
\[
   s_B([L_x^B,\Lambda_y^B])=0.
\]
This implies that the action of $[ L_x^B,\Lambda_y^B]$ on
$B_2$ coinsides with the action of $t_B([ L_x^B,\Lambda_y^B])$ on $V$.
Hence, \eqref{_commut_applied_to_y__in_B_Equation_} implies
\ref{_[L_x,Lambda_y]_to_y_is_-2x_Lemma_} (ii). This finishes
the proof of \eqref{_U_commu_with_commutator_Equation_},
case (i). $\;\;\blacksquare$

\hfill

Relation \eqref{_U_commu_with_commutator_Equation_},
case (ii) is implied by the following explicit
description of the commutators between $\g_0(A)$
and $\g_2(A)$, $\g_{-2}(A)$, which is valid for
many Lefschetz-Frobenius algebras of Jordan type.

\hfill

\proposition \label{_commu_between_g_0_and_g_+-2_Proposition_} 
Let $C=\oplus C_i$ be a Lefschetz-Frobenius algebra of Jordan type.
Let $\xi\in \g_0(C)$. Assume that $\g_0(C)$ acts on $C$ by derivations.
For any $x\in C_2$, denote by
$\xi(x)$ the image of $x$ under an action of $\xi:\; C_2\arrow C_2$.
Then $[\xi, L_x]= L_{\xi(x)}$ for all $x\in A_2$, $\xi\in\g_0$.

{\bf Proof:} By definition of a derivation, $\xi(xa)= \xi(x) a +x\xi(a)$.
By definition, $\xi L_x(a)= \xi(xa)$ and $L_x \xi (a)=x\xi(a)$.
Substracting one from another, we obtain $[\xi, L_x](a)= L_{\xi(x)}(a)$.
$\;\;\blacksquare$

\hfill

According to \ref{_g_0_derivatives_Corollary_},
$\g_0(A)$ acts on $A$ by derivations.
To show that $\g_0(B)$ acts on $B$ by derivations, we notice
the following. Let $C$ be an associative algebra over a field
and $\g$ be a Lie algebra, $\g\subset End(C)$. Consider the
corresponding Lie group $G\subset End(C)$. Then
$\g$ acts on $C$ by derivations if and only if $G$ acts
on $C$ by algebra automorphisms. Now, $\g_0(B)\cong \goth{so}(V)\oplus kH$.
The algebra $kH$ acts on $B$ by derivations for obvious reasons.
On the other hand, the Lie group $SO(V)$ acts on $B$
by automorphisms, as follows from definition of
$B$. Therefore, \ref{_commu_between_g_0_and_g_+-2_Proposition_}
can be applied to $A$ and $B$. Let $g\in \goth{so}(V)\oplus kH$.
Let $g^A$, $g^B$ be the elements of $\g_0(A)$, $\g_0(B)$
which correspond to $g$. Then,
\eqref{_U_commu_with_commutator_Equation_},
case (ii) is equivalent to

\[
   U^A([g^A, L_x^A]) = U^B([g^B, L_x^B]),\;\; \forall x\in V.
\]
By \ref{_commu_between_g_0_and_g_+-2_Proposition_},
this is equivalent to

\[ g^A(x)=g^B(x), \]
which is clear from definitions. We proved
\eqref{_U_commu_with_commutator_Equation_},
case (ii). $\;\;\blacksquare$

\hfill

It remains to prove \eqref{_U_commu_with_commutator_Equation_}
in case (iii). Consider the action of $\g_0(A)$, $\g_0(B)$ on
$\g_{-2}(A)$, $\g_{-2}(B)$ by commutators.
The action of $\g_0(B)$ on $\g_{-2}(B)$ consides
with that on $B_2\cong V$ (we use the standard isomorphism
$\g_{-2}(B)\cong V$ which is apparent from the explicit
description of $\g(B)$). This means that $k H\subset \g_0(B)$
acts on $\g_{-2}(B)$ trivially, and $\goth{so}(V)\subset \g_0(B)$
acts on $\g_{-2}(B)\cong V$ in a standard fashion. Denote this
action by $\rho_1:\; \goth{so}(V)\arrow End(V)$.
Similarly, $kH\subset \g_0(A)$ acts trivially on
$\g_{-2}(A)$. It remains to compare the action
of $\goth{so}(V)\subset \g_0(A)$ on $\g_{-2}(A)$
to $\rho_1$. Consider the isomorphism $r:\; A_2\arrow \g_{-2}(A)$
constructed in \ref{_g_-2_is_A_2_Corollary_}. The action
of $\g_0(A)$ on $\g_{-2}(A)$ defines an action
of $\goth{so}(V)\subset \g_0(A)$ on $A_2$. Using
the isomorphism $A_2\cong V$, we write this action
as a homomorphism $\rho_2:\;\goth{so}(V)\arrow End(V)$.
In this notation, the equation
\eqref{_U_commu_with_commutator_Equation_},
case (iii) is equivalent to the following statement:

\hfill

\lemma \label{_rho_1_is_rho_2_Lemma_} 
The representations $\rho_1$, $\rho_2$ coinside.

{\bf Proof:} Let $I\in Comp$.
Consider the endomorphism $ad I\in \g_0(A)$.
We identify $\g_0(A)$ and $\g_0(B)$ using the isomorphism
$U_0$. Let $t:\; \g_0(A)\arrow \goth{so}(V)$ be a projection
on a summand.
Let $C\subset \goth{so}(V)$ be the union of $t(ad I)$
for all $I\in Comp$. As we have seen, $C$ is Zariski
dense in $\goth{so}(V)$. Therefore it is sufficient
to show that $\rho_1$, $\rho_2$ coinside on $C$.
Let $x\in A_2$. By definition, $\rho_1(t(ad I))(x)= ad I(x)$,
and $\rho_2(t(ad I))x = r^{-1}[ad I, r(x)]$, where
$r:\; V\arrow \g_{-2}(A)$ is a map of \ref{_g_-2_is_A_2_Corollary_}.
To prove \ref{_rho_1_is_rho_2_Lemma_}
we have to show the following:

\begin{equation} \label{_r_commu_w_so(V)_Equation_}
   r(ad I(x)) = [ad I, r(x)].
\end{equation}

Both sides of \eqref{_r_commu_w_so(V)_Equation_} are
linear by $x$. Therefore it is sufficient to check
\eqref{_r_commu_w_so(V)_Equation_} in two cases:

\hfill

(i) $x\in H^{1,1}_I(M)$

(ii) $x\in H^{2,0}_I(M)\oplus H^{0,2}_I(M)$,

\hfill

\hspace{-6mm}where $H^{p,q}_I(M)$ is Hodge decomposition
associated with the complex structure $I\in Comp$.
In case (i), $ad I(x)=0$, so \eqref{_r_commu_w_so(V)_Equation_}
is equivalent to

\begin{equation} \label{_adI_commu_w_r(x)_for_x_in_H^1,1_Equation_}
   [ad I, r(x)] =0.
\end{equation}

Since elements of Lefschetz type are Zariski dense in $A_2$, it is
sufficient to prove \eqref{_adI_commu_w_r(x)_for_x_in_H^1,1_Equation_}
assuming that $x$ is of Lefschetz type. In this case,
$r(x)=(x,x)_{\c H}\Lambda_x$. Therefore,
\eqref{_adI_commu_w_r(x)_for_x_in_H^1,1_Equation_}
follows from the equation $[ad I, \Lambda_x] =0$. Since
$x\in H^{1,1}_I(M)$, the operator $L_x$ preserves weights
of the Hodge decomposition. An easy linear algebraic check insures that
in this case $\Lambda_x$ also preserves Hodge weights.
Therefore, by definition of $ad I$, we have $[ad I, \Lambda_x] =0$.
We proved \eqref{_r_commu_w_so(V)_Equation_}, case (i).

\hfill

Consider a non-zero holomorphic symplectic form $\tilde \Omega$
over the complex manifold $(M,I)$. Let $\Omega\in H^2(M)$ be
cohomology class represented by $\Omega$. Then $Im(\Omega)$,
$Re(\Omega)$ constitute basis in two-dimensional space
$H^{2,0}_I(M)\oplus H^{0,2}_I(M)$. Let
$\c H = (I, J, K, (\cdot,\cdot))$ be a hyperkaehler
structure such that $P_2(\c H)= Re(\Omega)$,
$P_3(\c H)= Im(\Omega)$. Such $\c H$ exists by
Calabi-Yau theorem. Since \ref{_r_commu_w_so(V)_Equation_}
is linear by $x$, we may check
\ref{_r_commu_w_so(V)_Equation_} case (ii) only for $x_2=P_2(\c H)$,
$x_3=P_3(\c H)$.
Clearly from definitions,
\[
    ad I(x_2) = 2 x_3, \; ad I(x_3) = -2 x_2.
\]
Let $c=(x_2,x_2)_{\c H}=(x_3,x_3)_{\c H}$.
By definition, $r(x_i)=c\Lambda_{x_i}$, $i=2,3$.
Therefore, \eqref{_r_commu_w_so(V)_Equation_} case (ii)
is equivalent to the following pair of equations:

\[ 4c \Lambda_{2x_3} = [ad I, c\Lambda_{x_2}], \]

and

\[ -4 c\Lambda_{2x_2} = [ad I, c\Lambda_{x_3}] \]

Since $\Lambda_{2a}=1/2\Lambda_a$, these two equations can
be rewritten as

\begin{equation}\label{_adI_on_Lambdas_Equation_}
   2\Lambda_{x_3} = [ad I, \Lambda_{x_2}], \;
   -2\Lambda_{x_2} = [ad I, \Lambda_{x_3}]
\end{equation}
\ref{_[L_2,Lambda_3]=adI_Claim_} implies that in notation of
\eqref{_so5_relations_Equation_}, $ad I = K_{23}$.
Therefore \eqref {_adI_on_Lambdas_Equation_}
is a consequence of \eqref{_so5_relations_Equation_}.
This proves \ref{_rho_1_is_rho_2_Lemma_}, and consequently,
\eqref{_U_commu_with_commutator_Equation_}
case (iii). Proof of \ref{_U_is_Lie_homomo_Proposition_}
is finished. This also finishes the proof
of \ref{_g(A)_for_hyperkae_Theorem_}, which spanned
four sections of this paper. $\;\;\blacksquare$


\section{The structure of the cohomology ring for
compact hyperkaehler manifolds.}
\label{_cohomolo_compu_Section_}


Let $M$ be a compact simple hyperkaehler manifold and $A=H^*(M)$
be its cohomology ring. Let $\c V=H^2(M,\R)$ considered as a linear
space with non-degenerate symmetric pairing $(\cdot,\cdot)_{\c H}$.
Applying \ref{_g(A)_for_hyperkae_Theorem_} and
\ref{_all_alg_with_so_are_^dA_Theorem_} to $A$, we immediately
obtain the following statement.

\hfill

\theorem\label{_cohomo_of_hyperk_are_^dA(V)_Theorem_} 
Let $A^r\subset H^*(M)$ be the subalgebra of $A=H^*(M)$
generated by $A_0$, $A_2$. Then $A^r\cong {}^dA(\c V)$, where
${}^dA(\c V)$ is a Frobenius algebra considered in
\ref{_all_alg_with_so_are_^dA_Theorem_}

$\blacksquare$

\hfill

We see that $A^r\subset H^*(M)$ can be expressed purely in
linear-algebraic terms. In this section, we make this
description as explicit as possible. Fix a linear space
$V$ with a non-degenerate bilinear symmetric form $B:\; S^2V\arrow k$,
where $k$ is a ground field. We express ${}^dA(V)$ in terms of
$V$ and $B$ as follows.

Let

\[
   C = C_0\oplus C_2\oplus ... \oplus C_{4d}
\]
be a graded linear space, with

\[ C_{2i} = S^iV, \;\;i\leq 2d, \]

\[ C_{2i} = S^{2i-d}V, \;\;i\geq 2d. \]

We describe a multiplicative structure on $C$ using some
classical results of linear algebra. Let $V$ be a linear space
equipped with non-degenerate bilinear symmetric product
$B:\; V\otimes V\arrow k$. Consider the Lie group $SO(V)$
associated with $V$ and $B$. This group naturally acts
on the symmetric powers $S^nV$ for all $n$. The representation
$SO(V)\arrow End(S^nV)$ is not irreducible. Its irreducible
decomposition is a classical result of linear algebra.
We describe this decomposition explicitely, and
define $SO(V)$-invariant multiplicative structure on $C$
in terms of this decomposition.

Let $\goth V:= \{ x_1,...., x_n\}$ be a basis in $V$.
We represent the vectors from $S^nV$ by the polynomials

\[
   \sum \alpha_{i_1,...,i_m} x_{i_1},..., x_{i_m},
\]
where $\alpha_{i_1,...,i_m}\in k$ and $x_{i_j}\in \goth V$, $k$
is a ground field. Consider an $SO(V)$-invariant vector $r\in S^2 V$
represented by the polynomial

\[
   r:= \sum_{i,j} B(x_i,x_j) x_i x_j.
\]
Let $L_r:\; S^nV\arrow S^{n+2}V$ be the map multiplying
the polynomial $P$ by $r$. Since the product in $S^*V$
commutes with the $SO(V)$-action and $r$ is $SO(V)$-invariant,
the map $L_r$ is a homomorphism of $SO(V)$-representations.

The scalar product $B$ on $V$ can be extended to an $SO(V)$-invariant
scalar product $(\cdot,\cdot)_{V^{\otimes_n}}$
on the space of $n$-tensors $\otimes^n V$ by the law

\[ (x_1\otimes x_2\otimes ... \otimes x_n,
    y_1\otimes y_2\otimes ... \otimes y_n)_{V^{\otimes_n}}
   = B(x_1,y_1) B(x_2,y_2) ... B(x_n,y_n).
\]
The space $S^nV\subset \otimes ^n V$ is $SO(V)$-invariant.
Using Schur's lemma, one can see that the restriction
of $(\cdot,\cdot)_{V^{\otimes_n}}$ to $S^nV$ is non-degenerate.
Denote this scalar product by $(\cdot,\cdot)_{S^n V}$.
For all $n>1$, the map $L_r:\; S^{n-2}V\arrow S^n V$
is an embedding. Let $R^n V\subset S^n V$ be the orthogonal
complement to the image of $L_r\; S^{n-2}V\arrow S^n V$.
Using the embedding $L_r:\; S^{i-2}V\arrow S^i V$
for different $i$,we obtain a decomposition

\begin{equation}\label{_decompo_of_S^nV_Equation_}
    S^n V\cong R^n V\oplus R^{n-2} V \oplus ... \oplus R^{n\;
\mbox{\tiny mod}\; 2} V,
\end{equation}
where $R^0V= k$ and $R^1 V = V$.

\proposition \label{_R^nV_is_irredu_Proposition_} 
For all $n\in \Bbb N$, the $SO(V)$-representation $R^n V$ is irreducible.

{\bf Proof:} See \cite{_Weyl_}. $\;\;\blacksquare$

\hfill

We obtain that \eqref{_decompo_of_S^nV_Equation_}
is an irreducible decomposition of an $SO(V)$-representation
$S^n V$. By definition, the representation

\[ S^nV / R^n V \oplus R^{n-2}V \oplus...
   \oplus R^{n-2i +2}V
\]
is canonically isomorphic to $S^{n-2i}V$. Denote the corresponding
quotient map by $B^n_{n-2i}:\; S^n V\arrow S^{n-2i}V$. Using the
maps $B^n_{n-2i}$ we define the multiplicative structure on $C$
as follows. Let $S^* V= \oplus S^n V$ be the algebra
of symmetric tensors over $V$. Let $\phi:\; S^* V\arrow C$,
$\phi=\oplus \phi_i$, where $\phi_i:\; S^iV \arrow C_{2i}$ is the following
map. For $2i\leq 2d$, $C_{2i}\cong S^i V$. For such $i$,
$\phi_i$ is defined as identity map. For $2d<2i\leq 4d$,
we have $C_{2i}= S^{2d-i}V$. Let $\phi_i: S^i V\arrow C_{2i}$
be equal to $B^i_{2d-i}:\; S^i V\arrow S^{2d-i}V$. For
$2i>4d$, $C_{2i}=0$ and we take $\phi_i=0$.
Clearly, the map $\phi$ is onto.

\hfill

\lemma \label{_ker_phi_ideal_in_S^*V_Lemma_} 
Let $I$ be a kernel of $\phi:\; S^*V\arrow C$. Then
$I$ is an ideal  in $S^*V$.

{\bf Proof:}
Let $x\in S^l V$, $x\in I$, and $y\in S^m V$. We have to show that
$xy\in I$. This relation is
vacuously true except in case when $d<l<l+m <2d$.
Let $\Lambda_r:\; S^n V\arrow S^{n-2} V$
be equal to $B^n_{n-2}$ for all $n = 2,3,...$. Clearly,
$B^n_{n-2i}= \Lambda^i_r$, where $\Lambda^i_r$ is $\Lambda_r$
to the power of $i$.
Therefore $x\in I$ is equivalent to $\Lambda^{l-d}_r(x)=0$,
and $xy\in I$ is equivalent to $\Lambda^{l+m-d}_r(xy)=0$.
Therefore, \ref{_ker_phi_ideal_in_S^*V_Lemma_}
is a special case of the following statement.

\hfill

\lemma \label{_Lambda_r_to_multi_Lemma_} 
Let $l,m,n$ be positive integer numbers, $x\in S^l V$,
$y\in S^m V$. Assume that $\Lambda^{n}_r(x) =0$.
Then $\Lambda^{n+m}_r(xy) =0$.

{\bf Proof:} Consider an element of $S^n V$ as polynomial function
on $V$ considered as an affine space. Let $\Delta:\; S^n V\arrow S^{n-2}V$
be the Laplace operator associated with the metric structure
$B$ on $V$.  By definition,

\[ \Delta(P) = \sum B(x_i,x_j)
   \frac{\partial^2}{\partial x_i \partial x_j} P,
\]
where $x_1,..., x_n$ is a basis in $V$. This definition
is independent of the choice of basis in $V$. The operator
$\Delta$ commutes with an action of $SO(V)$.
Checking that $\Delta(S^n V)$ contains $y^{n-2}$ for all $y\in V$,
we obtain that the map $\Delta$ is onto. This imples that
$\ker(\Delta^i:\; S^n V \arrow S^{n-2i} V)$ coinsides
with $R^n V\oplus R^{n-2} V\oplus ... \oplus R^{n-2i+2}V\subset S^n V$.
Therefore, $\ker(\Delta^i)=\ker(\Lambda_r^i)= \ker(B^n_{n-2i})$.
We obtain that \ref{_Lambda_r_to_multi_Lemma_} is equivalent
to the following statement:

\hfill

\lemma \label{_laplace_multi_Lemma_} 
Let $x\in S^l V$, $y\in S^m V$. Assume that $\Delta^i x =0$.
Then $\Delta^{i+m}(xy)=0$.

{\bf Proof:} We prove \ref{_laplace_multi_Lemma_}
using induction by $l$, $m$. We denote by $\bf L(l_0, m_0)$
the statement of \ref{_laplace_multi_Lemma_} applied to $l=l_0$,
$m=m_0$. Clearly,

\begin{equation}\label{_Laplace_of_produ_Equation_}
   \Delta(xy) = \Delta(x) y + x\Delta(y)
   + 2\sum \frac{\partial x}{\partial x_i}
     \frac{\partial y}{\partial x_j} B(x_i,x_j).
\end{equation}

By $\bf L(l-1, m)$, we have

\[ \Delta^{i+m -1}(\Delta(x) y) =0. \]
By $\bf L(l, m-2)$,
\[ \Delta^{i+m -1}(\Delta(x) y) =0. \]
Laplacian commutes with partial derivatives, and therefore
$\Delta^{i}(\frac{\partial x}{\partial x_i})=0$. Hence, by the
virtue of $\bf L(l, m-1)$,

\[ \Delta^{i+m -1} \bigg(\sum \frac{\partial x}{\partial x_i}
     \frac{\partial y}{\partial x_j} B(x_i,x_j)\bigg) =0
\]
Therefore, $\Delta^{i+m -1}$ applied to the right hand side
of \eqref{_Laplace_of_produ_Equation_} is zero.
This finishes the proof of \ref{_laplace_multi_Lemma_}.
We proved \ref{_Lambda_r_to_multi_Lemma_} and
\ref{_ker_phi_ideal_in_S^*V_Lemma_}. $\;\;\blacksquare$

\hfill

By definition, $C= S^*V / I$. Since $I$ is an ideal in $S^* V$,
the space $C$ inherits a canonical ring structure. Denote this
algebra by ${}^dC(V)$. The following
theorem characterizes the $H^2(M)$-generated subring of
cohomology ring of a simple compact hyperkaehler manifold
in terms of $C$.

\hfill

\theorem \label{_^dA(V)_is_C_Theorem_} 
Let $V$ be a linear space equipped with bilinear symmetric pairing $B_V$.
Then the algebra ${}^dA(V)$ is naturally isomorphic to ${}^dC(V)$.

{\bf Proof:} We prove \ref{_^dA(V)_is_C_Theorem_} as follows.
We consider ${}^dC(V)$ and ${}^dA(V)$ as graded linear spaces
with an action of the group $SO(V)$. These spaces are isomorphic
as graded $SO(V)$-representations. We notice that ${}^dA(V)$
is by definition a quotient of $SO(V)$ by the $SO(V)$-invariant
ideal $J$. We show that there is a unique graded $SO(V)$-invartiant
ideal $J$ in $S^*V$ such that $S^*V/J$ is isomorphic as a graded
$SO(V)$-representation to ${}^dC(V)\cong S^*(V)/I$. This implies
that $I$ coinsides with $J$, which proves \ref{_^dA(V)_is_C_Theorem_}.

In Section \ref{_minimal_Fro_Section_} we considered the graded Lie algebra
$\goth{so}(V,+)$. By definition, $\goth{so}(V,+)$ is a Lie algebra
of skew-symmetric endomorphisms of the space $V\oplus \goth H$.
Denote $V\oplus \goth H$ by $W$. The minimal Frobenius algebra
$A(V)\cong {}^1A(V)$ is isomorphic to $W$ as
$\goth{so}(W)$-representation. Therefore ${}^dA(V)$ is an
irreducible $\goth{so}(W)$-subrepresentation of $S^d W$ generated
by ${\Bbb I}\otimes{\Bbb I}\otimes...\otimes{\Bbb I}$. Consider the action
of the group $SO(W)$ on ${}^dA(V)$ which corresponds to
this Lie algebra action. We immediately obtain the following:

\hfill

\claim \label{_^dA(V)_is_R^d(W)_Claim_} 
Let $V$ be a linear space equipped with a non-degenerate bilinear
symmetric form, and ${}^dA(V)$ be a Frobenius
algebra defined in Section \ref{_^dA(V)_Section_}.
Let $\goth H$ be the two-dimensional
space with the hyperbolic metric, and $W:= V\oplus \goth H$.
As shown above, there is a natural action of $SO(W)$ on
the space ${}^dA(V)$. Earlier in this section, we defined the
$SO(W)$-representation $R^d W$. Then
${}^dA(V)$ is isomorphic to $R^d W$
as a representation of $SO(W)$.

$\blacksquare$

\hfill

Consider a natural embedding $SO(V)\subset SO(W)$
which corresponds to the decomposition $W=V\oplus \goth H$.
Consider $R^d W$ as a graded space, with the grading inherited
from ${}^dA(V)$. We proceed to demonstrate that, as a graded
$SO(V)$-representation, $R^d W\cong {}^dC(V)$. This is done
as follows. Let $r_V\in S^2 V$, $r_W\in S^2 W$ be $SO(V)$-invariant
(respectively, $SO(W)$-invariant) polynomials of degree 2defined earlier in
this
section. Denote the scalar product in $W$ by $B_W(\cdot,\cdot)$.
Earlier we denoted the scalar product in $V$ by $B_V(\cdot,\cdot)$.
Let $x_1, ..., x_n$ be a basis in $V$. Consider the vectors
${\Bbb I}$, $\Omega\in W$ which were introduced when we gave definition
of $A(V)\cong W$. By definition, all vectors $x_i$ are orthogonal
to ${\Bbb I}$ and $\Omega$, and $B_W({\Bbb I}, \Omega)=1$,
while $B_W({\Bbb I}, {\Bbb I})=0$ and $B_W(\Omega, \Omega)=0$.
Clearly, the vectors ${\Bbb I}, \Omega, x_1,...,x_n$ form
a basis in $W$. By definition,

\[
   r_W = {\Bbb I}\Omega - r_V.
\]
Consider a linear map $\gamma:\; S^d W\arrow S^d W$ which maps a monomial
$P = {\Bbb I}^i \Omega^j T$
to

\[ \gamma(P) :=
   \begin{array}{l}
   r_V^i \Omega^{j-i},
   \;\; j\geq i \\[5mm]
   {\Bbb I}^{i-j} r_V^j,
   \;\; j< i,
   \end{array}
\]
where $T= x_1^{\alpha_1} x_2^{\alpha_2}... x_n^{\alpha_n}$,
$\sum \alpha_i= d-i-j$.
Clearly, $\ker \gamma = r_W S^{d-2}W$. Therefore, the image of
$\gamma$ is naturally isomorphic as a linear space to $R^dW$. Consider a
multiplicative grading on $S^d W$ defined as follows:
$gr(\Bbb I)=0$, $gr(x_i)=2$, $i=1,2,..., n$, $gr(\Omega)=4$.
Clearly, this grading induces the standard one on
the space ${}^dA(V)\cong R^dW\subset S^d V$. By definition, the map
$\gamma:\; S^d W\arrow S^d W$ preserves this grading.
Since $\Bbb I$, $r_V$ and $\Omega$ are $SO(V)$-invariant,
the map $\gamma$ commutes with an action of $SO(V)$
on $S^d W$. Therefore, $\gamma(S^d W)$ is isomorphic to
$R^d W$ as a graded representation of $SO(V)$.
On the other hand, $\gamma(S^d W)$ is isomorphic
to ${}^d C(V)$ (again, as a graded representation of $SO(V)$
as the following argument shows.

For $2i\leq 2d$, the grading-$2i$ subspace
$\bigg( \gamma(S^d W)\bigg)_{2d}\subset \gamma(S^d W)$
is a linear span of monomials

\[ {\Bbb I}^{d-i} x_1^{\alpha_1} x_2^{\alpha_2} ... x_n^{\alpha_n},
   \;\; \sum \alpha_i = i.
\]
Similarly, for $2i>2d$, the space $\bigg( \gamma(S^d W)\bigg)_{2d}$
is a linear span of monomials

\[ \Omega^{i-d} x_1^{\alpha_1} x_2^{\alpha_2} ... x_n^{\alpha_n},
   \;\; \sum \alpha_i = 2d-i.
\]

Therefore, the grade $2i$ part of $R^dW\cong {}^d(W)$ is
$S^iV$ for $i<d$ and is $S^{2d-i}V$ for $i>d$. We proved the
following statement:

\hfill

\lemma \label{_C_iso_to_A(V)_as_SO(V)_repre_Lemma_} 
The spaces ${}^d C(V)$ and ${}^dA(V)$ are isomorphic as graded
representations of $SO(V)$.

$\blacksquare$

\hfill

By construction, ${}^d C(V)$ and ${}^dA(V)$ are quotient
algebras of $S^*V$. The canonical epimorphisms
$S^*V \arrow {}^d C(V)$ and $S^*V \arrow {}^d A(V)$
are $SO(V)$-invariant. Therefore, \ref{_^dA(V)_is_C_Theorem_}
is a consequence of \ref{_C_iso_to_A(V)_as_SO(V)_repre_Lemma_}
and the following lemma.

\hfill

\lemma \label{_quotie_alge_of_S^*V_iso_as_repre_iso_as_alge_Lemma_} 
Let $V$ be a linear space equipped with a scalar product.
Let $D$, $E$ be quotients of $S^*V$ by graded ideals
$I$, $J\subset S^*V$. Consider $D$, $E$ as graded algebras,
with the grading inherited from $S^*V$. Assume that $I$, $J$ are
$SO(V)$-invariant subspaces in $S^*V$ and $D$ is
isomorphic to $E$ as a graded representation of $SO(V)$.
Then $D$ is isomorphic to $E$ as an algebra.

{\bf Proof:} Consider the irreducible decomposition of $S^i V$
given by \eqref{_decompo_of_S^nV_Equation_}.
The summands of this decomposition are pairwise non-isomorphic.
Therefore, by Schur's lemma every $SO(V)$-invariant subspace
of $S^nV$ is a direct sum of several components of the
decomposition \eqref{_decompo_of_S^nV_Equation_}.
Let $I_n$, $J_n\subset S^n V$
be the $n$-th grade components of $I$, $J$. Since $I_n$, $J_n$
are $SO(V)$-invariant, these subspaces are direct sum of
several components of \eqref{_decompo_of_S^nV_Equation_}.
The quotient spaces $D_n =S^n V/I_n$, $E_n =S^n V/J_n$ are isomorphic
as $SO(V)$-representations. These spaces can be identified with the direct
sums of those components of the decomposition
\eqref{_decompo_of_S^nV_Equation_}
which don't appear in the decomposition of $I_n$, $J_n$.
Since $E_n$ is isomorphic to $D_n$, the spaces $I_n$ and
$J_n$ are isomorphic (as representations of $SO(V)$).
Since the components of \eqref{_decompo_of_S^nV_Equation_}
are pairwise non-isomorphic, the spaces $I_n$ and $J_n$ coinside.
This proves \ref{_quotie_alge_of_S^*V_iso_as_repre_iso_as_alge_Lemma_}.
\ref{_^dA(V)_is_C_Theorem_} is proven. $\;\;\blacksquare$


\section{Calculations of dimensions.}
\label{_calcu_dimensi_Section_}


We obtain easy numerical lower bounds on the dimensions
\[ \dim H^{p,q}(M)\] of Hodge components of cohomology spaces
of compact hyperkaehler manifolds. We have computed the
part of cohomology generated by $H^2(M)$.
The dimension of $\dim H^{p,q}(M)$ cannot be
lower than the dimension of the space $\bar H^{p,q}(M)$ of all
$(p,q)$ cohomology classes which are generated
by $H^2(M)$. In this section, we
compute dimensions of $\bar H^{p,q}$ for all $p$, $q$.

Let $\bar H^*(M)= \oplus_{p,q}\bar H^{p,q}\subset H^*(M)$ be the subring of
$H^*(M)$ generated by $H^2(M)$. Clearly,

\[
   \dim H^{p,q}(M) \geq \dim \bar H^{p,q}(M).
\]

By $p(n,m)$, we denote dimension of the space of
homogeneous polynomials of degree $m$ of $n$ variables.
This number is known from combinatorics as partition number.
It is given by the following formal serie, which was discovered
by Euler:

\[ \sum p(n,m) s^n t^m = \prod\limits_{i=1}^\infty
   \bigg(\frac{1}{1-t^is}\bigg)
\]

Consider the ring $S_2$ of polynomials of $n+1$ variables,
where $n$ variables are assigned degree 1 and one variable is
assigned degree 2. Let $p_2(n,m)$ be the space of homogeneous
polynomials of degree $m$ in $S_2$. Clearly,

\[
   p_2(n,m) = \sum_{i=0}^{i=[\frac{m}{2}]} p(n,m-i).
\]

\hfill

\theorem \label{_dimens_of_bar_H_^pq_in_terms_of_p_Theorem_} 
Let $M$ be a simple compact hyperkaehler manifold, $\dim_\R M =4d$,
$b_2(M)=n$ (we denote by $b_2$ the second Betti number).  Then

\hfill

(i) $\dim \bar H^*(M) = p(n,d) - p(n, d-2)$.

\hfill

(ii) $\dim \bar H^{2i}(M)= \begin{array}{l} p(n,i),\;\; i\leq d\\[2mm]
					    p(n,2d-i),\;\; i\geq d
			   \end{array}$

\hfill

(iii) $\dim\bar H^{p,q}(M)=0$ for $p+q$ odd.

\hfill

(iv) $\dim \bar H^{p,q}(M) = \dim \bar H^{p,2d-q}(M) =$ \\
\centerline{$=\dim \bar H^{2d-p,q}(M)=\dim \bar H^{2d-p,2d-q}(M)$.}

\hfill

(v) For $p+q\leq 2d$, $p\leq q$,

\[ \dim \bar H^{p,q}(M) = p_2(n-2,p) \]

\hfill

{\bf Proof:}

\hfill

(i) Follows from \ref{_^dA(V)_is_R^d(W)_Claim_}.

\hfill

(ii) \ref{_^dA(V)_is_C_Theorem_}

\hfill

(iii) Clear

\hfill

(iv) See \cite{_so5_on_cohomo_}

\hfill

(v) Let $V= H^2(M)$. \ref{_^dA(V)_is_C_Theorem_}
implies that $\bar H^{2m}(M) \cong S^m V$. Clearly, the Hodge decomposition
on $\bar H^{2m}(M)= S^mV$ is induced from that on
$V= H^{2,0}(M)\oplus H^{1,1}(M)\oplus H^{0,2}(M)$. The spaces
$H^{2,0}(M)$ and $H^{0,2}(M)$ are one-dimensional. Let $z\in H^{2,0}(M)$,
$\bar z \in H^{0,2}(M)$ be the non-zero vectors, and
$z,\bar z, x_1,...,x_{n-2}$ be the basis in $H^2(M)\cong V$.
Then the space $\bar H^{p,q}(M)$, $p+q\leq 2d$, $p\leq q$,
is a linear span of the monomials

\[ T_{a,b,\alpha_1,...\alpha_{n-2}} =
   z^a \bar z^b x_1^{\alpha_1} x_2^{\alpha_2}... x_{n-2}^{\alpha_{n-2}}
\]
where $b-a=q-p$, $\sum \alpha_i = p+q-(b+a)$.
Let $\Theta=z \bar z\in S^2 V$. Then

\[ T_{a,b,\alpha_1,...\alpha_{n-2}} = \Theta^{a} \bar z^{b-a}
   x_1^{\alpha_1} x_2^{\alpha_2}... x_{n-2}^{\alpha_{n-2}},
\]
where $a+b + \sum \alpha_i = p+q$. Since $b-a= q-p$,

\[ T_{a,b,\alpha_1,...\alpha_{n-2}} = \Theta^{a} \bar z^{p-q}
   x_1^{\alpha_1} x_2^{\alpha_2}... x_{n-2}^{\alpha_{n-2}},
\]
where $2a + (q-p) +\sum \alpha_i = p+q$. Translating $q-p$ to the
right hand side, we obtain that $T_{a,b,\alpha_1,...\alpha_{n-2}}$
is numbered by the different combinations of $a, \alpha_1,...\alpha_{n-2}$,
which satisfy the condition $2a +\sum \alpha_i = 2p$.
This number is by definition $ p_2(n-2,p)$. $\;\;\blacksquare$

\hfill

{\bf Acknowledgements:}
I am very grateful to my advisor David Kazhdan for warm support
and encouragement. This paper owes its existence to F. A. Bogomolov
and A. Todorov, whose studies of hyperkaehler manifolds made this
work possible. The last parts of this paper were inspired by
joint work with Valery Lunts.

I am extremely grateful to
Pierre Deligne, who was most kind and helpful. I owe Deligne several
important corrections in the final version of this manuscript.

Thanks due to B. A. Dubrovin, P. Etingof,
R. Bezrukavnikov, L. Positsel'sky, A. Todorov,
A. Polishchuk, T. Pantev and M. Bershadsky for stimulating
discussions. Prof. Y.-T. Siu and Prof. J. Bernstein kindly answered the
questions vital for the development of this work.

I am also grateful to MIT math department
for allowing me the use of their computing facilities.

\hfill


\begin{thebibliography}{20}

\bibitem[Beau]{_Beauville_}
Beauville, A. Varietes K\"ahleriennes dont la pere classe de Chern est
nulle. // J. Diff. Geom. 18, p. 755-782 (1983).

\bibitem[Bes]{_Besse:Einst_Manifo_}
Besse, A., Einstein Manifolds. // Springer-Verlag, New York (1987)

\bibitem[Bog] {_Bogomolov_}Bogomolov, F. A. , Hamiltonian Kaehler
manifolds. // Dokl. Akad. Nauk SSSR (Mat.), 243, 1101-1104 (1978).
English translation: Soviet Math. Dokl., 19, 1462-1465 (1978).

\bibitem[GH]{_Griffiths_Harris_}
Griffits, Ph. and Harris, J. Principles of algebraic geometry. //
Wiley-Interscience, New York (1978).


\bibitem[KS]{_Kodaira_Spencer_}
Kodaira K., Spencer D. C. On deformations of complex structures II
// Kodaira K., Collected Works vol. II, Princeton Univ Press (1975).

\bibitem[LL]{_Lunts-Loo_} Looijenga, E., Lunts, V. // preprint, 1994.


\bibitem[Spr]{_Springer_}  Springer T. A.,
     Jordan algebras and algebraic groups. //
     in series ``Ergebnisse der Mathematik und ihrer Grenzgebiete'',
     Bd. 75, Springer-Verlag, 1973.



\bibitem[Tod] {_Todorov_}
Todorov, A. Moduli of Hyper-K\"ahlerian manifolds I, II. // Preprint
MPI (1990)

\bibitem[Tod-Symp]{_Todorov:Symplectic_is_Kaehler_} Todorov, A. Every
holomorphic symplectic manifold admits a Kaehler metric. //
Preprint MPI.




\bibitem[V-$\protect\goth{so}$(5)]{_so5_on_cohomo_}
Verbitsky, M. On the action of a Lie algebra SO(5) on the cohomology
of a hyperk\"ahler manifold. // Func. Analysis and Appl. 24(2)
p. 70-71 (1990).

\bibitem[V-Bun]{_Verbitsky:Hyperholo_bundles_}
Verbitsky, M., Hyperholomorphic bundles.  //
alg-geom electronic preprint 9307008, 43 pages, LaTeX

\bibitem[V-Sym]{Verbitsky:Symplectic_I_}
Verbitsky, M., Hyperkaehler embeddings and
holomorphic symplectic geometry I. // alg-geom electronic preprint 9307009,
12 pages, LaTeX.

\bibitem[V-Sym II]{Verbitsky:Symplectic_II_}
Verbitsky, M., Hyperkaehler embeddings and holomorphic
symplectic geometry II. // alg-geom electronic preprint 9403006,
14 pages, LaTeX.

\bibitem[Wi]{_Weil_} Weil, Andre,
 Introduction a l'etude des varietes kahleriennes.//
 Paris, Hermann, 1958.


\bibitem[Wy]{_Weyl_}
   Weyl, H. The classical groups, their invariants and representations.
// Princeton Univ. Press, New York (1939)


\bibitem[Y]{_Yau:Calabi-Yau_}
Yau, S. T. On the Ricci curvature of a compact K\"ahler manifold
and the complex Monge-Amp\`ere equation I. // Comm. on Pure and Appl.
Math. 31, 339-411 (1978).

\end{thebibliography}
\end{document}